\documentclass[11.5pt]{amsart}

\usepackage{txfonts}
\usepackage[T1]{fontenc}

\usepackage[utf8]{inputenc} 
\usepackage{geometry} 
\usepackage{graphicx} 
\usepackage{booktabs} 
\usepackage{array} 
\usepackage{paralist} 
\usepackage{verbatim} 
\usepackage{subfig}
\usepackage{tabularx}
\usepackage{enumitem}
\usepackage{amsmath,amsthm,amscd}
\usepackage{amssymb,amsfonts}
\usepackage{fancyhdr} 
\usepackage{pgf}
\usepackage{caption}
\usepackage{rotating}
\usepackage{xcolor}
\usepackage{lscape}
\usepackage{xspace}
\usepackage{xfrac}
\usepackage{blindtext}
\usepackage{bm}
\usepackage[ruled,vlined]{algorithm2e}

\usepackage{soul}

\usepackage{tikz-cd}
\usepackage{tikz}
        \usetikzlibrary{patterns,hobby}
        \usepackage{pgfplots}
        \pgfplotsset{compat=1.6}

\usepackage{faktor} 
\usetikzlibrary{arrows.meta}
\usetikzlibrary{decorations.markings}
\usetikzlibrary{patterns}

\usepackage{color}   
\usepackage{hyperref}
\hypersetup{
    colorlinks=true, 
    linktoc=all,     
    linkcolor=blue,  
    citecolor=blue,
    filecolor=blue,
    urlcolor=blue
}

\usepackage[T1]{fontenc}
\usepackage{hyperref}
\usepackage{nameref}

\newcommand{\Ima}{\mathrm{Im\,}}
\newcommand{\Rea}{\mathrm{Re\,}}

\newcommand{\C}{\mathbb{C}}

\DeclareMathOperator{\res}{\textbf{res}}
\DeclareMathOperator{\per}{\textbf{per}}
\newcommand{\Bl}{\operatorname{Bl}}
\newcommand{\val}{\operatorname{val}}

\newcommand\CC{{\mathbb C}}
\newcommand\PP{{\mathbb P}}

\newcommand{\calB}{\mathcal B}
\newcommand\calH{{\mathcal H}}
\newcommand\calR{{\mathcal R}}

\newcommand\calG{\mathcal G}
\newcommand\calP{{\mathcal P}}
\newcommand{\calT}{\mathcal T}

\newcommand\Z{{\mathbb Z}}
\newcommand\ZZ{{\mathbb Z}}
\newcommand\calF{{\mathcal F}}

\newcommand\R{{\mathbb R}}

\DeclareMathOperator\MS{MS}
\DeclareMathOperator\Res{Res}

\newcommand{\cIH}{{\mathcal I\mathcal H}}

\newcommand{\RAMS}{\operatorname{RAMS}}
\newcommand{\PRes}{\textbf{PRes}}

\usepackage[
	    style=alphabetic,
            isbn=false,
            doi=false,
            url=false,
            backend=bibtex,
            maxnames = 5,
            sorting=nty,
           ]{biblatex}
\defbibheading{bibliography}[\bibname]{%
   \section*{References}%
}

\DeclareFieldFormat[article]{title}{{\it #1}}
\DeclareFieldFormat{journaltitle}{{\rm #1}}
\renewbibmacro{in:}{%
  \ifentrytype{article}{}{\printtext{\bibstring{in}\intitlepunct}}}
  \bibliography{biblio}

\newtheorem{thm}{Theorem}[section]
\newtheorem{cor}[thm]{Corollary}
\newtheorem{prop}[thm]{Proposition}
\newtheorem{lem}[thm]{Lemma}
\newtheorem{opprob}[thm]{Open Problem}

\theoremstyle{definition}
\newtheorem{defn}[thm]{Definition}
\theoremstyle{remark}
\newtheorem{rmk}[thm]{Remark}
\theoremstyle{definition}

\theoremstyle{definition}

\theoremstyle{definition}
\newtheorem{ex}[thm]{Example}
\theoremstyle{definition}

\theoremstyle{definition}

\numberwithin{equation}{section}
\title{The residual map as a stratified fibration over the resonance arrangement}

\author{Dawei Chen} 
\address[Dawei Chen]{Department of Mathematics, Boston College, Chestnut Hill, MA 02467, USA}
\email{dawei.chen@bc.edu}

\author{Scott Mullane} 
\address[Scott Mullane]{School of Mathematics and Statistics, University of Melbourne, Australia}
\email{mullanes@unimelb.edu.au}

\author{Nick Salter} 
\address[Nick Salter]{Department of Mathematics, University of Notre Dame, 255 Hurley Building,
Notre Dame, IN 46556}
\email{nsalter@nd.edu}

\author{Guillaume Tahar}
\address[Guillaume Tahar]{Department of Mathematics with Computer Science, Guangdong Technion - Israel Institute of Technology (GTIIT), Shantou, Guangdong, China}
\email{guillaume.tahar@gtiit.edu.cn}

\newcommand{\para}[1]{\noindent \textbf{#1.}}

\newcommand{\into}{\hookrightarrow}

\date{\today}
\keywords{
Residual map, Meromorphic differentials, Translation structures, Resonance arrangement, Braid arrangement, Residue-aligned multi-scale differentials}

\usepackage{pinlabel}
\begin{document}
\begin{abstract}
Given a stratum $\mathcal{H}(\mu)$ of meromorphic one-forms on a Riemann surface with prescribed orders of zeros and poles, the residual map assigns to each differential the collection of its residues at the poles. In this paper, we prove that the residual map is a locally trivial topological fibration over each open flat of the stratification induced by the resonance arrangement.
\par
To give an algebro-geometric proof of this result, we introduce a refined compactification of the projectivized stratum $\mathbb{P}\mathcal{H}(\mu)$, namely the moduli space of residue-aligned multi-scale differentials $\operatorname{RAMS}(\mu)$. We prove that the projectivized residual map extends to $\operatorname{RAMS}(\mu)$ as a stratified fibration over the projectivized resonance arrangement, thereby yielding a resolution of the projectivized residual map. This approach naturally extends to the analogue of the residual map for strata of $k$-differentials.
\par
Complementing the algebro-geometric approach, we also develop a flat-geometric point of view. In genus zero, each isoresidual fiber carries a canonical higher-dimensional analogue of a translation structure. These fibers admit a metric completion whose boundary has a modular interpretation in terms of the WYSIWYG compactification, for which we provide explicit local models. Isoresidual fibers admit a polyhedral cell decomposition whose geometry varies linearly with the underlying configuration of residues in period coordinates. Degenerations of these cells occur precisely when resonance conditions are satisfied, providing a geometric proof that the residual map is a stratified fibration. Finally, certain isoresidual fibers admit natural interpretations as Hurwitz spaces of branched covers, yielding unexpected homeomorphisms between moduli spaces of covers with prescribed ramification profiles.
\end{abstract} 
\maketitle
\setcounter{tocdepth}{1}
\tableofcontents

\section{Introduction}

For any collection of non-negative integers $a_{1},\ldots,a_{n}$ and any collection of negative integers $-b_{1},\ldots,-b_{p}$ forming a partition $\mu$ of $2g-2$, we write $\mathcal{H}(\mu)$ for the moduli space of pairs $(X,\omega)$ formed by a compact Riemann surface $X$ of genus $g$ and a meromorphic one-form $\omega$ on $X$ with labeled zeros and poles of respective orders $a_{1},\dots,a_{n}$ and $b_{1},\dots,b_{p}$, considered up to biholomorphic equivalence. For technical reasons, we introduce the possibility of using marked points, which then appear as zeros of order $0$. Each moduli space $\mathcal{H}(\mu)$ carries the structure of a complex-analytic orbifold of dimension $2g+n+p-2$ (and it has no orbifold points in the case of $g=0$ and $n+p > 2$ as $\mathbb{CP}^1$ with at least three labeled points admits only the trivial automorphism). 

\begin{defn}\label{defn:ResidualSpace}
For $p \geq 2$, the \textit{residual space} $\calR_{p}$ is the complex vector space consisting of configurations $\lambda = (\lambda_{1},\dots,\lambda_{p})$ satisfying $\sum_{j=1}^{p} \lambda_{j} = 0$.
\end{defn}

Each moduli space $\mathcal{H}(\mu)$ is endowed with a \textit{residual map} $$\res\colon \mathcal{H}(\mu) \to \calR_{p}$$
$$\omega \mapsto (\Res_{1}\omega,\dots,\Res_{p}\omega)$$ that assigns to each differential $\omega$ the sequence of its residues at the poles.
\par
Work combining flat geometry and intersection theory on moduli spaces has led to a precise description of the fibers of $\res$ for strata of differentials on $\mathbb{CP}^{1}$ with $n=1$ and $n=2$ (see \cite{GeTaIso, ChPrIso, CGPT}). In these beginning cases, the residual map appears to be singular only over a hyperplane arrangement $\mathcal{A}_{p} \subset \calR_{p}$ in the residual space, called the \textit{resonance arrangement} (see Definition~\ref{defn:resonanceintro}).
\par
In this paper, we show that for a moduli space $\mathcal{H}(\mu)$ with an arbitrary number $n$ of zeros, the residual map is a \textit{stratified fibration} in the sense that, for every stratum $B$ of the stratification of $\calR_{p}$ induced by the resonance arrangement $\mathcal{A}_{p}$, the restriction of the residual map to $\res^{-1} (B)$ is locally trivial, hence defines a topological fibration over $B$.

\subsection{The stratified fibration}

We first recall the definition of the resonance arrangement.

\begin{defn}\label{defn:resonanceintro}
Let $I$ be a nonempty proper subset of $[p] = \{1, \dots, p\}$. The \textit{resonance hyperplane} $A_{I}$ is the subspace defined in the residual space $\calR_{p}$ by the equation $\lambda_{I}=0$ where $$\lambda_{I}=\sum\limits_{i \in I} \lambda_{i}.$$ 
We refer to the set of all the resonance hyperplanes as the \textit{resonance arrangement} $\mathcal{A}_{p}$ in $\calR_{p}$.
\par
For any family $E$ of resonance hyperplanes in $\mathcal{A}_{p}$, let $B(E)$ denote the subset of $\calR_{p}$ consisting of configurations of residues belonging to all hyperplanes in $E$, and to no resonance hyperplane outside $E$. We refer to $B(E)$ as an \textit{open flat}.

\par
The open flats $B(E)$, as $E$ ranges over families of resonance hyperplanes in $\mathcal{A}_{p}$, endow $\calR_{p}$ with a structure of \textit{stratified space} we refer to as $(\calR_{p},\mathcal{A}_{p})$.
\end{defn}

Given a configuration $\lambda$ of residues in $\calR_{p}$, we refer to the subset $\mathcal{F}_{\lambda}(\mu) = \res^{-1}(\lambda)$ in $\mathcal{H}(\mu)$ as an \textit{isoresidual fiber}. The residual map induces a stratification of the moduli space $\mathcal{H}(\mu)$ by pullback of the stratification of the residual space. For any open flat $B$ of $(\calR_{p},\mathcal{A}_{p})$, we write $\mathcal{H}(\mu,B)$ for the preimage subset $\res^{-1} (B)$ in $\mathcal{H}(\mu)$. 

\begin{thm}\label{thm:MAIN}
For any stratum $\mathcal{H}(\mu)$ of meromorphic one-forms and any open flat $B$ in $(\calR_{p}, \mathcal{A}_{p})$, the restriction of the residual map $\textnormal{\textbf{res}}$ to $\mathcal{H}(\mu, B)$ is a locally trivial fibration over $B$. In other words, $\textnormal{\textbf{res}}$ is a stratified fibration of $\mathcal{H}(\mu)$ over $(\calR_{p},\mathcal{A}_{p})$.
\end{thm}

In Section~\ref{sec:RAMS}, we give an \textit{algebro-geometric proof} of Theorem \ref{thm:MAIN}, relying on the properties of the projectivized residual map in the moduli space of residue-aligned multi-scale differentials, a refined version of the multi-scale compactification.
\par
In the particular case of moduli spaces of meromorphic one-forms on $\mathbb{CP}^{1}$, we give a \textit{flat-geometric proof}, consisting of endowing each isoresidual fiber with a geometric structure that can be regarded as a higher-dimensional analogue of a translation surface, decomposing it into polyhedra, and showing that its deformations lead to degenerations of this geometric structure only upon crossing a resonance hyperplane.
\par
This second proof, although longer and currently limited to the genus-zero setting, involves several intermediate results that are of independent interest, as they provide an explicit description of the isoresidual fibers and their deformations. In the remainder of this introductory section, we outline this strategy and the main results obtained along the way.

\subsection{Higher-dimensional translation structures}\label{SS:higherdimtrans}

In any stratum $\mathcal{H}(\mu)$ of meromorphic or holomorphic one-forms, whether in genus zero or higher genus, fixing the periods of the differential along closed loops, i.e. the absolute periods, defines the \textit{isoperiodic foliation}. In genus zero, this amounts to fixing a connected component of an isoresidual fiber. Given an isoperiodic leaf $\mathcal{L}$ in $\mathcal{H}(\mu)$, the remaining degrees of freedom in the period coordinates (see Section~\ref{sub:PeriodAtlas}), after fixing the residues and the absolute periods, are given by the relative periods, namely the integrals of the differential along arcs connecting the $n$ zeros.
\par
These relative periods $(Z_{ij})_{1 \leq i, j \leq n}$ naturally assemble into a matrix-valued system of local coordinates, yielding a local identification of $\mathcal{L}$ with a space of \textit{exact matrices} $\mathcal{E}_n(\mathbb{C})$ (see Section~\ref{sub:ExactMatrices}). This latter space also canonically identifies with the quotient space $V_{n} = \mathbb{C}^n / \Delta$, where $\Delta$ is the diagonal line. 
\par
Different choices of connecting arcs produce different coordinate systems, yet as these differ only by absolute periods, any isoperiodic leaf $\mathcal{L}$ is naturally equipped with a translation atlas. The isoperiodic leaf $\mathcal{L}$ inherits from $V_{n}$ a \textit{Hermitian flat metric} which allows us to define the \textit{metric completion} $\mathcal{L}^{\mathrm{MET}}$ of the leaf (see Section~\ref{sub:Translation}).
\par
The boundary of the metric completion consists of limit objects resulting from collisions of zeros. These objects belong to the {\em WYSIWYG} (What You See Is What You Get) compactification, obtained from multi-scale limit differentials by retaining the top level and contracting the lower levels (see Section~\ref{subsec:local}).
\par
Collision-type multi-scale limit differentials are detected by the vanishing of the corresponding relative periods. Section~\ref{subsec:local} compares their local metric completion with the operation of contracting lower levels in WYSIWYG. 

\begin{thm}\label{thm:MetricCompletion}
In every stratum $\mathcal{H}(\mu)$, any isoperiodic leaf $\mathcal{L}$ can be endowed with a canonical translation structure that is locally isomorphic to the flat Hermitian space $V_{n} = \mathbb{C}^n / \Delta$.
\par
Any nonzero point $\omega$ of its metric completion $\mathcal{L}^{\mathrm{MET}}$ admits a neighborhood that is locally isomorphic to a branched cover of $V_{n}$, whose branch locus is contained in the \textit{collision arrangement} $\mathcal{B}_{n}$ (also called the braid arrangement modulo translations; see Section~\ref{sub:ExactMatrices}).
\end{thm}

Theorem~\ref{thm:MetricCompletion} is proved in Section~\ref{subsec:local} through Propositions~\ref{prop:translation} and \ref{prop:local}.

\begin{rmk}
There are several obstacles to obtaining a flat-geometric proof of Theorem~\ref{thm:MAIN} for isoresidual fibers in moduli spaces of differentials of higher genus. First, the isoresidual fibers would no longer be translation manifolds, but would instead be endowed, like $\mathcal{H}(\mu)$, with an integral-affine structure induced by the period atlas. Second, since there is no Poincaré--Bendixson theorem in higher genus, critical leaves of the vertical foliation may accumulate. Consequently, there is no straightforward analogue of the critical graph that would provide a combinatorial decomposition of the fiber; see Section~\ref{sub:IntroDeformation}.
\end{rmk}

\begin{rmk}
The stratified Euclidean structure of the isoperiodic leaves in the space of stable forms (i.e., sections of the dualizing bundle) was considered by McMullen in \cite{McMullen}. In contrast, our approach explicitly defines a translation structure, extends it to moduli spaces of meromorphic one-forms with arbitrary poles, and gives a modular interpretation of its metric closure in terms of the WYSIWYG compactification.
\end{rmk}

\subsection{Deformations of translation manifolds}\label{sub:IntroDeformation}
In genus zero, the translation structure on a completed isoresidual fiber $\overline{\mathcal{F}_{\lambda}}(\mu)$ admits an explicit description as a polyhedral cell complex. Here we give an overview of this theory, as developed in Sections \ref{sec:Polyhedral} and \ref{sec:Period}. 
\par
Our analysis is built around a study of the directional foliations of meromorphic differentials in genus zero. Adjusting via the $\C^*$-action if necessary, we can assume throughout that our residues satisfy the following genericity condition.

\begin{defn}[Vertically generic, cf. Definition \ref{def:vgen}]
A configuration of residues $\lambda = (\lambda_1, \dots, \lambda_p) \in \calR_p$ (or a differential $\omega \in \mathcal{F}_{\lambda}(\mu)$) is said to be {\em vertically generic} if for all subsets $I \subset [p]$ such that $\Rea(2 \pi \mathrm{i} \lambda_I) = 0$, then moreover $\lambda_I = 0$. A subset $U \subset \calR_p$ is likewise said to be vertically generic if each $\lambda \in U$ is.
\end{defn}

For $\lambda$ vertically generic, the vertical directional foliation of any differential $\omega \in \mathcal{F}_{\lambda}(\mu)$ is particularly simple (cf. Lemma \ref{lemma:vnrcrit}). There are no periodic or minimal leaves: every leaf is compact and flows from one singularity to another. The structure of the foliation is thus encoded by a \textit{critical graph} consisting of all critical leaves of the foliation (i.e., leaves incident to a zero of the differential).
The case of a boundary element in the metric completion $\overline{\mathcal{F}_{\lambda}}(\mu)$ is analogous, except that, since it represents a singular curve, one must consider the critical graph on a disjoint union of possibly several irreducible components, each isomorphic to $\mathbb{CP}^{1}$. 
\par
The decomposition of $\overline{\mathcal{F}_{\lambda}}(\mu)$ according to the possible critical graphs yields a polyhedral cell structure: each critical graph carries a finite set of {\em distinguished saddle connections} (Definition \ref{def:distsaddle}), whose periods furnish coordinates on the set of differentials with a given critical graph. The periods of distinguished saddles are not completely independent; they satisfy {\em cycle equations} (Definition \ref{def:cyclespace}) that record when a collection of saddles determines a closed loop, whose period is therefore the sum of the enclosed residues. These conditions, in combination with a positivity constraint arising from a preferred orientation on each distinguished saddle, lead to an unbounded polyhedral region of geometrically feasible periods. 
\par
For a given critical graph $\Gamma$ and configuration of residues $\lambda$, we let $\calP^\lambda_\Gamma$ denote this polyhedron. The key technical result of Sections \ref{sec:Polyhedral} and \ref{sec:Period}, Theorem \ref{lem:Pproduct}, shows that the combinatorial type of $\calP^\lambda_\Gamma$ is locally invariant along an open flat.

\begin{thm}\label{lem:Pproduct}
   Let $\Gamma$ be any critical graph, let $B \subset \calR_p$ be an open flat, and let $U\subset B$ be connected and vertically generic. Then the combinatorial isomorphism type of $\calP_\Gamma^\lambda$ is independent of $\lambda \in U$. Moreover, the equations defining $\calP_\Gamma^\lambda$ are all of the form $\sum \epsilon_i \Rea(s_i) \ge \Rea(2 \pi \mathrm{i}\lambda_I)$ and $\sum \epsilon_i\Ima(s_i) \ge 0$, where the variables $s_i$ denote periods of a certain finite set of distinguished saddle connections and $\epsilon_i \in \{-1,0,1\}$.
\end{thm}

We obtain the flat-geometric proof of Theorem \ref{thm:MAIN} by showing that $\calF_\lambda(\mu)$ is homeomorphic to a polyhedral complex $\calP^\lambda$ constructed by gluing the polyhedra $\calP^\lambda_\Gamma$ that parametrize the differentials realizing a given critical graph (see Theorem~\ref{thm:periodmaphomeo}), where Theorem~\ref{lem:Pproduct} is the key ingredient in establishing the local triviality of the residue map.

\subsection{Topological complexity of isoresidual fibers}

One application of the polyhedral decomposition and the correspondence between polyhedral cells and critical graphs is a polynomial upper bound for the Betti numbers of isoresidual fibers in genus zero. It is worth noting that the constants $C_{n,p}$ can grow very quickly as functions of $n$ and $p$ (they count certain graphs with $n+p$ labeled vertices). The content of Theorem \ref{thm:BettiEstimate} is that for fixed $n$ and $p$, the Betti numbers grow at most polynomially as a function of the {\em orders} of the zeros and poles.

\begin{thm}\label{thm:BettiEstimate}
Given a partition $\mu = (a_1,\ldots, a_n, -b_1,\ldots, -b_p)$ of $-2$, set $B\coloneqq\sum_{j=1}^p b_j$. There exists a constant $C_{n,p}$, depending only on $n$ and $p$ such that for any configuration of residues $\lambda \in \calR_{p}$, the Betti numbers of the isoresidual fiber $\mathcal{F}_{\lambda}(\mu)$ satisfy
\[
\operatorname{Betti}_i(\mathcal{F}_{\lambda}(\mu))
\leq C_{n,p}(2B - 2p + 1)^{n + p - 3}
\quad \text{for } 0\leq i\leq n-1
\]
and moreover,
\[
\operatorname{Betti}_i(\mathcal{F}_{\lambda}(\mu)) = 0
\quad \text{for } i\geq n.
\]
The Betti numbers of the metric completion satisfy the same bound
\[
\operatorname{Betti}_i(\overline{\mathcal{F}_{\lambda}}(\mu))
\leq C_{n,p}(2B - 2p + 1)^{n + p - 3}
\quad \text{for } 0\leq i\leq 2n-2
\]
with vanishing in the range
\[
\operatorname{Betti}_i(\overline{\mathcal{F}_{\lambda}}(\mu)) = 0
\quad \text{for } i\geq 2n - 1.
\]
\end{thm}

Theorem~\ref{thm:BettiEstimate} is proved in Section~\ref{sec:Betti}.

\subsection{Generalization to $k$-differentials}

For $k \geq 2$ and $\mu$ being a partition of $k(2g-2)$, $\mathcal{H}^{k}(\mu)$ denotes the moduli space of $k$-differentials of type $\mu$, parameterizing pairs $(X,\xi)$ where $X$ is a compact Riemann surface of genus $g$ and $\xi$ is a possibly meromorphic section of the line bundle $K_{X}^{\otimes k}$ whose zero and pole orders are  prescribed by $\mu$.
\par
The study of moduli spaces of $k$-differentials is a relatively new subject, but one that has seen rapid progress in recent years; see, for instance, the classification of the connected components of strata of $k$-differentials in \cite{chge,apisa2026components}.
\par
An analogue of the isoresidual fibration can also be defined for strata of $k$-differentials. There are, however, several notable differences. First, only poles whose order is a multiple of $k$ carry an additional local invariant besides their order, namely the $k$-residue. Second, there is no analogue of the residue theorem for $k\geq 2$. As a consequence, for a stratum with $p$ poles of order divisible by $k$, the residual space $\calR^{k}_{p}$ has dimension $p$ rather than $p-1$. Moreover, the resonance conditions are not expressed directly in terms of the $k$-residues. Rather, they correspond to the vanishing of partial sums of $k^{th}$ roots of these $k$-residues. They do not define linear subspaces of the residual space, but they nevertheless define a stratification.
\par
For moduli spaces where the isoresidual fibers are discrete ($k$-differentials on $\mathbb{CP}^{1}$ where every singularity except two is a pole of order divisible by $k$), the cardinality of isoresidual fibers over each stratum of the stratification of the residual space has been computed in \cite{chen2025finite}.
\par
In the present paper, by applying the same methods as in the algebro-geometric proof of Theorem~\ref{thm:MAIN}, we prove in Theorem~\ref{thm:kRAMS-submersion} that for every moduli space $\mathcal{H}^{k}(\mu)$, the $k$-residual map is a stratified fibration over the stratification of the $k$-residual space $\calR^{k}_{p}$ induced by the resonance subvarieties.

\subsection{Relations with other works}

\subsubsection{Isoresidual fibration: state of the art}

Recent progress in the study of the isoresidual fibration has been driven by two main approaches, sometimes developed independently and sometimes combined: flat-geometric methods and intersection-theoretic methods.
\begin{itemize}
\item a characterization of empty isoresidual fibers in every connected component of every stratum of meromorphic one-forms, see \cite{getaab};
\item a complete description of discrete isoresidual fibers, including explicit formulas for their cardinalities on each open flat of the resonance stratification, see \cite{GeTaIso,ChPrIso};
\item for strata with one-dimensional isoresidual fibers (corresponding to the case $n = 2$), the translation structure of the fibers was thoroughly investigated in \cite{CGPT}. In particular, the authors established a complete correspondence between degenerations of parametrized differentials and the local invariants of singularities arising in the fibers;
\item as an application of the results on one-dimensional isoresidual fibers, Theorem~1.7 of \cite{CGPT} provides a classification of the connected components of generic fibers for an arbitrary number $n$ of zeros.
\end{itemize}
Unfortunately, little is known about isoresidual fibers of higher dimension (corresponding to the case $n \geq 3$). The present work provides foundational results towards the construction of an appropriate compactification of such strata.
\par
In the present work, we prove that the residual map defines a topological fibration away from the resonance arrangement. A natural and intriguing question is then to understand its monodromy. For $n=1$ and $p \leq 3$, a complete description of the monodromy groups was obtained in \cite{GeTaIso}, revealing in particular the exotic embedding $\mathfrak{S}_{5} \subset \mathfrak{S}_{6}$ as the monodromy group of the isoresidual fibration of $\mathcal{H}(6, -3, -3, -2)$.
\par
For $n \geq 2$, the monodromy is no longer described by permutations but rather by elements of the mapping class group. For this reason, it is natural to replace each isoresidual fiber by a suitable homological object and to study the induced monodromy action in the form of a Gauss--Manin connection. In Section~7 of \cite{CGPT}, several structural results and explicit examples were provided in this direction for the case $n=2$. In particular, the appropriate homological object turns out to be the one naturally associated with a translation surface, namely the first relative homology group of the compactified fiber punctured at the poles, relative to the zeros of the meromorphic one-form defining the translation structure.

\subsubsection{Applications to moduli theory}

More recently, the geometry of the isoresidual fibration has found applications in moduli theory. For example, the curve class given by the general fiber of the isoresidual fibration with exactly two zeros was used to show that the pseudo-effective cone of divisors of certain moduli spaces of pointed stable rational curves are non-polyhedral, and hence these spaces fail to be Mori Dream Spaces; see \cite{mullane2025isoresidual}.
\par
In Section~\ref{sec:isoresidual-elliptic-surface}, we illustrate the interaction between flat, elliptic, and birational geometry for the signature $\mu=(2, 2, -2, -2, -2)$. The canonical translation structure on a generic compactified isoresidual curve in this case has signature $(4^4, -2^8)$. In the associated rational elliptic surface, the zeros form a degree-four multisection with nontrivial monodromy around the resonance locus, while the poles form eight sections. The elliptic surface has Mordell--Weil rank two, giving infinitely many $(-1)$-curves and hence a non-polyhedral pseudo-effective cone. Hence the singularities of the translation structure have a concrete interpretation as curves controlling the birational geometry of the total space. 
\par
Another area of applications arises from the multi-scale compactification of strata of abelian differentials; see \cite{BCGGM2}. In this setting, the boundary of a stratum consists of nodal curves equipped with meromorphic differentials whose residues at the poles are required to satisfy certain resonance conditions.
\par
These boundary objects are parametrized by products of \textit{generalized strata}, namely loci in strata of meromorphic differentials where the poles are partitioned into subsets on which the corresponding residues sum to zero. The connected components of these generalized strata have recently been classified by Lee and collaborators; see \cite{lee2024connected, lee2025connected, lee2025one}.
\par
A long-term objective of this program is the computation of the cohomology ring of the boundary complex. To date, nontrivial classes in the third cohomology group of compactified strata have been constructed using, among other ingredients, results stemming from the study of the isoresidual fibration; see \cite{chen2026non}.

\subsubsection{Translation manifolds}

In isoresidual fibers of strata with $n \geq 2$ zeros, relative periods between the zeros form a local system of coordinates defined up to translations in some complex vector space. In other words, they are (open) translation manifolds.
\par
Defining a reasonable notion of a translation structure on a compact complex manifold requires specifying which singularities are admissible. At present, no consensus has emerged on such a notion. An interesting candidate is given by polyhedral K\"{a}hler manifolds, introduced in \cite{Pan}, which exhibit good geometric properties. In particular, recent works of Borbon and Panov \cite{de2025polyhedral,  de2026miyaoka} show that they provide a natural framework for establishing a Miyaoka–Yau inequality for hyperplane arrangements in $\mathbb{CP}^{n}$, generalizing classical results of Hirzebruch on complex algebraic surfaces; see \cite{hirzebruch1983arrangements}.
\par
However, polyhedral Kähler manifolds are constructed by gluing flat simplices and do not seem to correspond naturally to some complex-analytic data. In contrast, isoresidual fibers directly arise as moduli spaces of complex-analytic objects. It would be highly desirable to have a natural category of translation manifolds that would keep the two-faceted quality of translation surfaces: a complex-analytic definition (complex structure with a differential) and an equivalent geometric definition (translation atlas extending to some well-behaved singularities). It is natural to expect that any reasonable category of translation manifolds should include suitably compactified isoresidual fibers, at least for generic configurations of residues.
\par
Within such a framework, it should be possible to compute the topological invariants of the complex manifold (notably Hodge and Chern numbers) directly from the combinatorics of the singular loci, much as one recovers the Euler characteristic of a Riemann surface from the orders of the zeros and poles of a translation structure. This would, for instance, allow one to locate the isoresidual fibers corresponding to strata with $n=3$ zeros within the Enriques–Kodaira classification of complex surfaces.

\subsubsection{Isoperiodic foliation}
The foliation obtained by fixing the absolute periods of a holomorphic differential within a stratum has been studied extensively, both from a topological and a dynamical perspective, in \cite{McMullen, calsamiglia2023transfer}. It plays an important role in the study of totally geodesic submanifolds of $\mathcal{M}_{g,n}$; see, for instance, \cite{wright2020totally}. Connections with Teichm\"{u}ller dynamics have been established, as well as analogues of Ratner's theorems. Numerous partial ergodicity results have also been obtained for these foliations on the connected components of strata; see, for instance, \cite{Hamenstadt2018,winsor2021dynamics}.
\par
In contrast, the study of the isoperiodic foliation on strata of meromorphic differentials has so far been carried out only for very particular strata (see \cite{calsamiglia2025isoperiodic, FaTaZa}). The problem can be formulated as follows. 

\begin{opprob}
Given an isoresidual fiber $\mathcal{F}_{\lambda}(\mu)$ in a stratum $\mathcal{H}(\mu)$ of meromorphic one-forms, find conditions on $\lambda$ and $\mu$ under which the isoperiodic foliation is ergodic in $\mathcal{F}_{\lambda}(\mu)$.
\end{opprob}

Among the isoperiodic foliations, a particularly remarkable subcase is provided by real-normalized differentials, namely differentials whose absolute periods are all real. The Riemann bilinear relations imply, in particular, that such differentials necessarily have poles. They arise naturally in the theory of integrable systems; see, for instance, \cite{grushevsky2013real} for an application to the spectral theory of the elliptic Calogero--Moser system.
\par
Real-normalized differentials also play an important role in the algebraic geometry of moduli spaces. Using them, Krichever proved in \cite{krichever2012real} the Arbarello conjecture concerning the intersection of compact complex subvarieties of $\mathcal{M}_g$ with the Weierstrass loci. They also led Grushevsky and Krichever \cite{grushevskykrichever2009} to a new proof of the Diaz theorem on the upper bound of the dimensions of complete subvarieties in the moduli space of curves (see \cite{diaz1984bound}).
\par
To establish ergodicity results for the isoperiodic foliation in moduli spaces of real-normalized differentials, Krichever, Lando, and Skripchenko introduced in \cite{KLS} a combinatorial model in the case of differentials with a single double pole. Using the translation structure introduced in the present paper, such a combinatorial model can be constructed in full generality. Indeed, on the one hand, the horizontal directional flow and its decomposition into irreducible components provide a combinatorial classification of these differentials into distinct types. On the other hand, the imaginary parts of the relative periods between pairs of zeros give a family of globally defined harmonic functions. Their zero loci decompose the isoperiodic leaves into chambers within which the combinatorial type remains constant. Ergodicity of the isoperiodic foliation for real-normalized differentials can then be reformulated as a question about the incidence relations between these chambers.

\subsubsection{Stratified braid groups}\label{subsub:Stratified}

The third-named author came to be interested in the theory of isoresidual fibers by way of a connection with braids and polynomials, the theory of so-called {\em stratified braid groups}. These were introduced in \cite{SaBraid} and studied further in \cite{SaBraid2, SaHuxBraid}. The objective here is to understand a stratification on the classical configuration space $\operatorname{Conf}_n(\C)$ of $n$-tuples of distinct points in $\C$. Viewing $\operatorname{Conf}_n(\C)$ as the space of monic squarefree polynomials of degree $n$, it carries the {\em equicritical stratification} induced by tracking the multiplicities of the critical points of the polynomials. There are $n-1$ critical points in total, so that strata are indexed by a partition $\kappa$ of $n-1$; we let $\operatorname{Conf}_n(\C)[\kappa]$ denote the corresponding stratum.

The connection with the isoresidual fibration arises as follows: associating $f \in \operatorname{Conf}_n(\C)[\kappa]$ to its logarithmic derivative $\frac{df}{f}$ yields a holomorphic map $d\colon \operatorname{Conf}_n(\C)[\kappa] \to \calF_\lambda(\kappa, -1^{n+1})$, where the $n+1$ simple poles correspond to the roots of $f$ and $\infty$, with residue $\lambda_i = 1$ at each root of $f$ by the argument principle. The fibers of $d$ correspond to orbits of $f(z)$ under the action of the group of affine transformations $az + b$ of $\C$, realizing an equicritical stratum as the total space of an (orbifold) fiber bundle over a particular isoresidual fiber.

The works \cite{SaBraid, SaBraid2, SaHuxBraid} are concerned with a group-theoretic question arising from this, and fit more broadly into the enterprise of understanding fundamental groups and monodromy groups of moduli spaces. Tracking the $n$ roots and $k$ critical points, there is an embedding $i_\kappa\colon\operatorname{Conf}_n(\C)[\kappa] \to \operatorname{Conf}_{n+k}(\C)$, and hence a homomorphism $i_{\kappa,*}\colon \pi_1(\operatorname{Conf}_n(\C)[\kappa]) \to B_{n+k}$, where $B_{n+k}$ is the classical braid group on $n+k$ strands. The paper \cite{SaBraid2} gives a description of the image of this map as a ``framed braid group'', a certain subgroup of $B_{n+k}$ consisting of braids ``compatible'' with the flat structure on $\C$ arising from the logarithmic differential. In a complementary direction, \cite{SaHuxBraid} shows that in the two critical point case $\kappa = \{a, b\}$, the map $i_{\kappa,*}$ is not injective.

Analogous questions can be formulated for any isoresidual fiber, by considering instead the inclusion map $i_{\lambda,\mu}\colon \calF_\lambda(\mu) \to \calH(\mu)$ and taking the induced map on $\pi_1$. Note that $\calH(\mu)$ can be identified with a certain $\C^*$-bundle over a configuration space of points in $\mathbb{CP}^{1}$ (tracking the locations of zeros and poles) and so its fundamental group is closely related to a braid group. Theorem \ref{thm:MAIN} provides important structural information about $i_{\lambda,\mu}$, though by itself it does not seem to be sufficient to address the central question of when $i_{\lambda,\mu,*}\colon \pi_1(\calF_\lambda(\mu)) \to \pi_1(\calH(\mu))$ is an injection. 

\subsubsection{Hurwitz spaces}

In a more general setting than the one considered in Section~\ref{subsub:Stratified}, for certain Hurwitz spaces of rational maps with prescribed ramification orders, logarithmic differentials describe an isoresidual fiber in which the configuration of residues is prescribed by the ramification orders. By establishing that the residual map is a stratified fibration, Theorem~\ref{thm:MAIN} then shows that the homeomorphism type of these Hurwitz spaces is entirely determined by the resonance conditions satisfied by the configuration of residues of the corresponding isoresidual fiber. This bridge between Hurwitz spaces and isoresidual fibers is established in Section~\ref{sec:isoresidual-hurwitz-fibers}, and more precisely in Corollary~\ref{cor:Hurwitz-homeomorphisms}. This also provides a solution to the Hurwitz realization problem for some of these spaces (i.e., proving that they are nonempty), as shown in Proposition~\ref{prop:Hurwitz-realization}.
\subsection{Organization of the paper}

\begin{itemize}
    \item In Section~\ref{sec:Background}, we review the basic material on translation structures, meromorphic differentials, and period coordinates on strata of differentials.
    \item In Section~\ref{sec:RAMS}, we discuss the residue-aligned multi-scale compactification and give the algebro-geometric proof of Theorem~\ref{thm:MAIN}. We also generalize these constructions to $k$-differentials.
    \item In Section~\ref{sec:wysiwyg-metric-completion}, we introduce the higher-dimensional translation structure on isoresidual fibers (in genus zero) and isoperiodic leaves (in higher genus). We then define their metric completion and provide a modular interpretation of the boundary points in terms of the WYSIWYG compactification. We describe their local geometry in terms of ramified covers of a braid arrangement, thereby proving Theorem~\ref{thm:MetricCompletion}.
    \item In Section~\ref{sec:Polyhedral}, we introduce a polyhedral cell structure on isoresidual fibers in genus zero and prove the invariance of the combinatorial type of each polyhedron as the configuration of residues varies in a non-resonant subset of an open flat of the resonance arrangement (Theorem \ref{lem:Pproduct}).    
    \item In Section~\ref{sec:Period}, we give the global construction of the polyhedral period domain in genus zero and describe its topology. This provides a flat-geometric proof of Theorem~\ref{thm:MAIN} in genus zero.
    \item In Section~\ref{sec:Betti}, drawing on the polyhedral cell decomposition, we obtain upper bounds on the Betti numbers of isoresidual fibers and prove Theorem~\ref{thm:BettiEstimate}.
    \item In Section~\ref{sec:isoresidual-hurwitz-fibers}, we interpret certain isoresidual fibers as Hurwitz-type spaces and deduce from Theorem~\ref{thm:MAIN} conditions under which two of them are homeomorphic.
    \item In Section~\ref{sec:isoresidual-elliptic-surface}, we study an interesting interplay between flat, elliptic, and birational geometry for the isoresidual fibration with the signature $\mu = (2,2,-2,-2,-2)$. 
\end{itemize}
\par 
\paragraph{\bf Acknowledgements.}
During the preparation of this work, D.C. was supported by NSF grant DMS-2301030 and Simons Travel Support for Mathematicians; S.M. was supported by DECRA Grant DE220100918 from the Australian Research Council; N.S. was supported by NSF grant DMS-2338485; G.T. was supported by the French National Research Agency under the project TIGerS (ANR-24-CE40-3604). Part of the work was carried out when S.M. visited Boston College in February 2026 and when D.C. visited BIMSA in August 2026. The authors thank these institutions for their hospitality.   
\smallskip 
\paragraph{\bf AI usage disclosure.} Some of the authors used ChatGPT for literature review, exposition, and proofreading. The ideas and content are the authors' own. 

\section{Meromorphic differentials and translation structures on  Riemann surfaces}\label{sec:Background}

By a slight abuse of terminology, we refer to a pair $(X,\omega)$, where $X$ is a Riemann surface of genus $g$ and $\omega$ is a meromorphic one-form, as a \textit{meromorphic differential of genus $g$}.

\subsection{Translation atlas on the regular locus}\label{sub:TransRegular}

Given a pair $(X,\omega)$ formed by a compact Riemann surface $X$ with a nonzero meromorphic one-form $\omega$, we denote by:
\begin{itemize}
    \item $Z_{\omega}$ the set of the zeros of $\omega$;
    \item $P_{\omega}$ the set of the poles of $\omega$;
    \item $X^{\ast}$ the complement of $Z_{\omega} \cup P_{\omega}$ in $X$, also called the \textit{regular locus} of $(X,\omega)$.
\end{itemize}
Primitive functions of~$\omega$ are locally injective on $X^{\ast}$ and therefore provide an atlas of holomorphic local coordinates. Since local primitives of the same differential differ by a constant, the transition maps of this atlas are translations of the complex plane. We refer to such an atlas as a \textit{translation atlas} and $(X,\omega)$ as a \textit{translation surface}.
\par
Geometric notions that make sense on the Euclidean plane $(\mathbb{C},dz)$ globalize on a translation surface. In particular, $X^{\ast}$ is endowed with a flat metric $|\omega|$. Conversely, any translation atlas defines a complex structure on $X^{\ast}$, and the pullbacks of the one-form $dz$ on $\mathbb{C}$ by the local charts globalize to a holomorphic differential $\omega$ on $X^{\ast}$.
\par
In a translation surface $(X,\omega)$, for each slope $\theta$ in the circle $\mathbb{S}^{1}$ of directions, there is a foliation ${\mathrm{Fol}}_{\theta}$ of the regular locus $X^{\ast}$, where the leaves are locally conjugated to oriented lines of slope $\theta$ in the translation charts.

\subsection{Local models for zeros}\label{sub:LocalModelsZeros}

Up to a holomorphic change of variables, a zero of order $a$ in a holomorphic one-form is given as $z^{a}dz$. The flat metric of the translation structure induced by the differential is obtained by pulling back the standard Euclidean metric of the complex plane via a cyclic cover of degree $a+1$ ramified over a point. Thus, the flat metric induced by the differential extends to the zero as a conical singularity of angle $(a+1)2\pi$.

\subsection{Local models for simple poles}\label{sub:LocalModelsSimplePoles}

Up to a biholomorphic local change of variables, in a neighborhood of a simple pole, the differential $\omega$ is of the form $\dfrac{rdz}{z}$, where $r$ is the residue at the pole. Its flat-geometric interpretation is that of the end semi-infinite cylinder, foliated by a one-parameter family of closed geodesics having a period equal to $2\pi {\mathrm{i}} r$.

\subsection{Local models for poles of higher order}\label{sub:LocalModelsPolesHigherOrder}

In the language of translation surfaces, the standard Euclidean plane is written as $(\mathbb{C},dz)$. It can be compactified as $\hat{\mathbb{C}} \cong \mathbb{CP}^{1}$ where $dz$ has a residueless double pole at infinity.
\par
A \emph{flat cone of order $a \geq 0$} is the translation surface defined by the differential $z^{a}dz$ on $\hat{\mathbb{C}} \cong \mathbb{CP}^{1}$. It has a unique zero of order a and a unique pole of order $a+2$. This surface can be obtained as a cyclic cover of degree $a+1$ of the Euclidean plane, branched at a single point. Neighborhoods of the point at infinity in these cones provide flat-geometric local models for a residueless pole of order $a+2$.
\par
To construct local models for poles of order $b \geq 2$ with a nontrivial residue, we consider a flat cone of order $b-2$ and remove an $\epsilon$-neighborhood of a semi-infinite line emanating from the conical singularity, together with a neighborhood of the singularity itself. We then identify the resulting boundary segments via an isometry. By suitable rotation and rescaling, this produces a pole of order $b$ with any prescribed nonzero residue. This construction is explained in detail in \cite{Bo}.
\par
Just as the flat-geometric interpretation of a zero is a conical singularity, the flat-geometric interpretation of a pole of order at least two is a neighborhood of infinity in a cone, which we may call an \emph{anti-conical singularity}, following the convenient terminology introduced in \cite{BCT25}.

\subsection{Relative homology and saddle connections}

In a translation surface $(X,\omega)$, a smooth path $\gamma$ joining two zeros of $\omega$ in the punctured surface $X \setminus P_{\omega}$ represents a \textit{relative homology class} $[\gamma]$ in $H_{1}(X \setminus P_{\omega},Z_{\omega}; \mathbb Z)$.
\par
A \textit{saddle connection} is an arc joining two zeros, whose interior contains no singularity and is locally conjugated to a segment of constant slope in the translation charts. In the flat metric induced by $\omega$, a saddle connection is a geodesic segment. A \textit{closed saddle connection} is a saddle connection whose endpoints coincide.

\subsection{Period atlas on strata of differentials}\label{sub:PeriodAtlas}

The integrals of $\omega$ along a basis of the $\mathbb{Z}$-module $H_{1}(X \setminus P_{\omega}, Z_{\omega}; \mathbb Z)$ provide local complex-analytic coordinates for the ambient stratum, which are called {\em period coordinates}.
\par
After passing locally to a homologically marked cover, the period coordinates take values in
$$
H^1(X\setminus P_\omega, Z_\omega; \mathbb C).
$$
Their restriction to
$$
H^1(X\setminus P_\omega; \mathbb C)
$$
gives the {\em absolute periods}. Fixing the absolute periods defines the local isoperiodic leaves. By the exact sequence of the pair,
$$
\ker\left(
H^1(X\setminus P_\omega,Z_\omega; \mathbb C)
\to
H^1(X\setminus P_\omega; \mathbb C)
\right)
\cong
H^0(Z_\omega;\mathbb C)/\mathbb C.
$$
Thus the tangent directions along an isoperiodic leaf are precisely the relative periods between the zeros. In genus zero, the absolute homology of $X\setminus P_\omega$ is generated by peripheral loops around the poles, whose periods are $2\pi \mathrm{i}$ times the residues. Hence fixing the residues fixes all absolute periods.
\par
For a stratum $\mathcal{H}(a_{1},\dots, a_{n}, -b_{1},\dots, -b_{p})$ with $\sum\limits_{i=1}^n a_{i} - \sum\limits_{j=1}^p b_{j} = 2g-2$ and $n, p \geq 1$, the \textit{period atlas} defines a structure of a complex-analytic orbifold of complex dimension $2g + n + p-2$.

\section{Residue-aligned multi-scale differentials and the resolution of the residual map}\label{sec:RAMS}

We introduce a residue-aligned refinement of the multi-scale compactification constructed for all genera and signatures in \cite{BCGGM} and \cite{BCGGM2}. Its role is to resolve the projectivized residual map at the residueless loci, while retaining a modular interpretation analogous to ordinary multi-scale differentials.

We keep the notation of the preceding sections. For a moduli space $\mathcal{H}(\mu)$ of meromorphic one-forms with $p$ marked poles $p_1,\ldots,p_p$, the residual map is 
\[
  \res\colon \mathcal H(\mu)\to \calR_{p}, \quad (X,\omega)  \mapsto (\Res_{p_1}\omega,\ldots,\Res_{p_p}\omega),
\]
and its projectivization is the rational map
\[
  \PRes\colon \PP \calH(\mu)\dashrightarrow \PP \calR_{p}.
\]

Throughout this section we exclude the case in which $\res$ is identically zero on $\mathcal H(\mu)$. In that case the unprojectivized residual map is the constant map to the origin of $\calR_{p}$, and the whole stratum is a single isoresidual fiber. This is the only case in which there is no nonzero residue function to projectivize; for instance it occurs when there is only one marked pole (i.e., $p = 1$). Indeed, when there are at least two marked poles, the residue function at each pole cannot be identically zero on $\calH(\mu)$, as each residue contributes to the local period coordinates subject to the residue theorem as the only constraint. 

\subsection{Review of multi-scale differentials}\label{subsec:ordinary-MS-review}

We recall the ingredients of the multi-scale compactification that will be refined below; see \cite{BCGGM} and \cite{BCGGM2} for the full construction. A {\em multi-scale differential} is represented by a stable marked curve with an {\em enhanced level graph} $\Gamma$. The vertices of $\Gamma$ correspond to the irreducible components and carry levels in a totally ordered set, which we normalize as $0,-1,\ldots,-N$. Edges joining vertices on the same level are called {\em horizontal}, and edges joining different levels are called {\em vertical}. Each component $X_v$ carries a possibly meromorphic differential $\eta_v$, with the prescribed orders at the marked points and compatible orders at the branches of every node. At a horizontal node the residues at the two branches add to zero. 

At a vertical node, the lower branch is a pole for the lower-level differential; the {\em prong-matching} structure records the cyclic gluing choice of local horizontal directions. The {\em rescaling ensemble} records the relative scales of the levels. In local multi-scale coordinates these scales are encoded by {\em adjusting parameters} $h_v$, defined up to units, such that $h_v^{-1}\omega$ has leading term $\eta_v$ on $X_v$. Multiplying all differentials on one level by a scalar gives the {\em level-rotation equivalence}, and the projectivized {\em multi-scale compactification} $\MS(\mu)$ also quotients by overall scaling \cite[Sections 9.2 and 13.3]{BCGGM2}.

Conversely, given an enhanced level graph $\Gamma$ and differentials $\eta_v$ on the components $X_v$ satisfying the aforementioned conditions, the {\em global residue condition} (GRC) determines whether it can be smoothed into the respective stratum of differentials $\calH(\mu)$, namely, whether it can be an actual limit when differentials in $\calH(\mu)$ degenerate from smooth curves to nodal curves. (If we do {\em not} impose the GRC, the remaining structure is called a {\em generalized} multi-scale differential; see \cite[Section 3]{CGHMS}.) 

More precisely, let $L$ be a level and let $Y$ be a connected component of the subgraph $\Gamma_{>L}$ lying strictly above $L$. If $Y$ contains a marked pole lying above $L$, no GRC equation is imposed. Otherwise, if $E(Y,L)$ denotes the set of vertical edges from $Y$ to vertices at level $L$, then the GRC is
\[
  \sum_{e\in E(Y,L)} \Res_{q_e^-}(\eta_{v_-(e)}) = 0,
\]
where $q_e^-$ is the lower end of $e$. Equivalently, in the bottom-up plumbing construction this is the condition for the existence of the {\em modification differential} that opens the nodes between $Y$ and level $L$ \cite[Sections 9 and 10]{BCGGM2}.

There is an analogous compactification for loci of differentials with a set $\Lambda$ of {\em prescribed linear residue conditions}. In that setting the ordinary GRC in the above is replaced by a generalized GRC, also called the $E_\Lambda$-GRC, which records the boundary specializations of the chosen linear residue equations together with the vertical-edge residues appearing at the relevant level cuts. This is the smoothability condition for the corresponding residue-constrained boundary strata; see \cite[Proposition~4.2]{CMZ} and \cite[Theorem~3.2]{CGPT}.

For {\em strata of $k$-differentials}, the analogous multi-scale compactification can be described after passing to the {\em canonical cyclic cover}: the chosen $k^{th}$-root one-form satisfies ordinary matching and residue conditions on the cover, while the structures descend by equivariance under the deck transformation to the stratum of $k$-differentials. The resulting compactification imposes the corresponding conditions for prong-matching, scaling, and $k$-GRC; see \cite[Sections~2, 4, and 5]{BCGGM3} and \cite[Sections~2--4]{CMZarea}.

The space $\MS(\mu)$ is stratified by the boundary strata with fixed types of enhanced level graphs. Its local coordinates are horizontal smoothing parameters, level scaling parameters, lower side residues at vertical edges, and (perturbed) periods on each level, subject to the matching, prong, level-rotation, and GRC conditions. The refined residue-aligned compactification below will keep the same multi-scale structure and add only the scaling information needed to remember the first nonzero scale of the residue function at each marked pole. 

\subsection{Alignment of residueless marked poles to lower levels}

We first give the modular description of the residue-aligned multi-scale space. Let $\Gamma$ be the enhanced level graph of a multi-scale differential, with normalized levels $0,-1,\ldots,-N$. We treat each marked point as a leg (i.e., a half-edge) in $\Gamma$, and call it a {\em marked polar (or zero) leg} if the marked point corresponds to a pole (or zero). 

A {\em residue-aligned level graph} $\widetilde \Gamma$ is obtained from $\Gamma$ by allowing the endpoint of each marked polar leg to lie on a level not higher than the level of the vertex to which the leg is attached. These levels may be original levels of $\Gamma$ or new levels inserted by the residue alignment. The (possibly lower) endpoint of a marked polar leg records the scale of the corresponding residue function; the marked pole remains attached to the original component containing it. Polar legs carry no prong-matching structure, as prong-matching is  needed only for vertical nodes.

Suppose the marked polar leg corresponding to $p_i$ is attached to a vertex $v$. If the residue of the multi-scale differential $\eta_v$ at $p_i$ on the component $X_v$ is nonzero, the endpoint remains on the level of $v$, and this residue is recorded. If this residue is zero, then a residue alignment must assign the marked polar leg a level 
strictly lower than that of $v$. For a one-parameter smoothing transverse to the zero-residue locus, the new level of the endpoint of the polar leg records the scale at which the scaled residue has a finite nonzero limit. (This is similar to scaling periods that shrink to zero so that the corresponding subsurfaces are placed at  lower levels.) The nonzero leading coefficient of the residue function is denoted by $r_{p_i}$ and is called the {\em aligned residue} at $p_i$. 

Note that the aligned endpoint of a marked pole may lie on an existing level, on a new residue-only level, or below all pre-existing levels. In particular, if the residue function of a marked pole is identically zero on the stratum  (e.g. as the unique pole, or considered on a linear submanifold with residueless constraint), then its endpoint is aligned to the lowest level (lower than all other poles with non-identically-zero residues) and its aligned residue value remains zero. 

Equivalently, the endpoint level of $p_i$ is determined by comparing the residue function at $p_i$ with the local adjusting parameters of the multi-scale differential. Thus every marked pole placed on a residue level carries a nonzero aligned residue. Marked poles below the highest nonzero residue level nevertheless contribute zero under the (resolved) projectivized residual map defined in Section~\ref{subsec:extension} below, since $\PRes$ records only the leading projective residue directions. 

For a vertex $v$, the residue theorem for the differential $\eta_v$ involves only the residues at the level of $v$. If a marked pole attached to $v$ is aligned to a lower level, its residue at the level of $v$ is zero; its nonzero aligned residue enters instead in the residue equation for the construction of the modification differential at the lower aligned level.

\begin{ex}[Alignment of a residueless differential]
\label{ex:residueless}
    Suppose a family of differentials has three marked poles with residues given by 
    $$ t^2-2t^3, \quad -t^2, \quad 2t^3, $$
    respectively, as $t\to 0$. The limit differential at $t = 0$ is residueless and $[0,0,0]$ is not a valid image point of $\PRes$ in the space of projectivized residue tuples. However, in the setting of residue-aligned differentials, at $t = 0$, 
    we align the first and second marked polar legs to a lower level (according to the adjusting parameter $t^2$), with aligned residues 
    $1$ and $-1$, respectively, and we align the third marked polar leg to a further lower level (according to the adjusting parameter $t^3$) with aligned residue $2$; see the following figure:  
    \begin{center}
\begin{tikzpicture}[scale=1,
    vertex/.style={circle,fill=black,inner sep=1.8pt},
    every node/.style={font=\small}]
  \node[vertex,label=above:$v$] (v) at (0,2.5) {};

  \draw[dashed] (v) -- (-1.7,1.15);
  \draw[dashed] (v) -- ( 0,1.15);
  \draw[dashed] (v) -- ( 2.25,-0.25);

  \node[anchor=north east] at (-1.7,1.15) {$p_1$};
  \node[anchor=north] at (0,1.15) {$p_2$};
  \node[anchor=north west] at (2.25,-0.25) {$p_3$};
  \node[anchor=south east] at (-1.55,1.15) {$1$};
  \node[anchor=south west] at (0.12,1.15) {$-1$};
  \node[anchor=south west] at (2.10,-0.25) {$2$};
\end{tikzpicture}
\end{center} 

    Then we define the image of $\PRes$ at $t = 0$ to be 
    $$[1,-1,0]$$
    using the higher-level aligned residues from the first two marked poles. Note that as projective tuples, for $t\neq 0$, we have 
    $$ 
    [t^2-2t^3, -t^2, 2t^3] = [1-2t, -1, 2t]. 
    $$
    It has limit $[1,-1,0]$ as $t\to 0$, which coincides with the   extension of $\PRes$ using the aligned residues. 
    \end{ex}

\begin{defn}[Residue-aligned multi-scale differential]
\label{def:residue-aligned-MS}
A \emph{residue-aligned multi-scale differential} is a multi-scale differential, with its usual prong-matchings, rescaling ensemble, and level-rotation equivalence, together with a residue-aligned level graph and aligned residues at the marked poles, subject to the aligned global residue condition (A-GRC) defined in~\eqref{eq:A-GRC} below.
\end{defn}

If A-GRC is omitted, the remaining structures define \emph{generalized} residue-aligned multi-scale differentials as in \cite[Section~3]{CGHMS}; they are smoothable into the interior of the stratum if and only if the relevant A-GRC holds. 
\par 
Let $T_\Gamma$ be the ordinary level-rotation torus with the twist group quotient used in the multi-scale construction, and let $R(\widetilde{\Gamma})$ be the set of new residue-only levels. Locally the enlarged action of $T_{\widetilde{\Gamma}}$ is obtained by adjoining one $\CC^*$-factor for every level in $R(\widetilde{\Gamma})$. If the endpoint of a marked polar leg $p_i$ lies on level $L$, rotation of $L$ by $c\in\CC^*$ sends $r_{p_i}$ to $cr_{p_i}$. In particular, when the component carrying $p_i$ lies above $L$, rotation by that upper component level does not change $r_{p_i}$. If $L$ is residue-only, its rotation acts only on the residue level scale and on the aligned residues with endpoints on $L$; it does not act on any vertex or prong-matching. Overall scaling is according to the same projectivization convention as in the ordinary multi-scale space.

Two residue-aligned level graphs are equivalent if there is an isomorphism of their underlying residue-aligned level graphs preserving the enhanced structures, together with an order-preserving bijection of their enlarged level sets. Two residue-aligned multi-scale differentials are isomorphic when their underlying multi-scale differentials are isomorphic, their aligned graphs are equivalent in this sense, and their aligned residues correspond under the action of the enlarged level-rotation torus.

\subsection{The aligned global residue condition}
\label{subsec:A-GRC}

Let $\widetilde\Gamma$ be a residue-aligned level graph. For a level $L$, let $\widetilde\Gamma_{>L}$ be the subgraph lying strictly above $L$, and let $Y$ be a connected component of $\widetilde\Gamma_{>L}$. If $Y$ contains the endpoint of a marked polar leg lying strictly above $L$, no aligned global residue condition is imposed for $(Y,L)$. This is the ordinary GRC exception: a pole whose scaled residue is already above the cut level can absorb the residues of lower-level modification differentials without changing its leading aligned term. With prescribed linear residue conditions involving such a pole, one imposes the corresponding generalized A-GRC as in the usual generalized GRC or $E_\Lambda$-GRC formalism \cite[Section~4.1]{CMZ} and \cite[Section~3.2]{CGPT}.

Otherwise, let $E(Y,L)$ be the set of vertical edges from $Y$ to vertices at level $L$, and let $P(Y,L)$ be the set of marked polar legs whose adjacent vertices lie in $Y$ and whose aligned endpoints lie on level $L$. Edges and marked polar legs with lower endpoints strictly below $L$ do not appear in the residue condition at level $L$, since their contributions vanish after scaling by the adjusting parameter of $L$. Without further ado, we state the \emph{aligned global residue condition}, abbreviated A-GRC, as 
\begin{equation}\label{eq:A-GRC}
  \sum_{e\in E(Y, L)} \Res_{q_e^-}(\eta_{v_-(e)})
  - \sum_{p_i\in P(Y, L)} r_{p_i} = 0.
\end{equation}
Here $q_e^-$ denotes the lower branch of the node corresponding to the edge $e$, and $\eta_{v_-(e)}$ is the multi-scale differential on the lower-level component containing $q_e^-$.

The minus sign in \eqref{eq:A-GRC} reflects the {\em boundary orientation}. The residue of a vertical edge is measured at the lower-side pole of the node. If a marked pole $p_i$ is aligned to the lower level $L$, then its residue at the original  level of its adjacent vertex is zero, and $r_{p_i}$ is the leading coefficient obtained after division by the adjusting parameter of $L$. In the construction of the modification differential at level $L$, a vertical edge with lower-side residue $\rho_e$ contributes $-\rho_e$ on the upper side, whereas an aligned marked pole contributes $r_{p_i}$. The residue theorem for the modification differential gives
\[
  \sum_{p_i\in P(Y,L)} r_{p_i}
  - \sum_{e\in E(Y,L)} \rho_e = 0,
\]
which is equivalent to \eqref{eq:A-GRC}. 

\begin{ex}[A residue-aligned level graph and A-GRC]
\label{ex:V-shaped-AGRC} Consider a  residue-aligned level graph with two top vertices $v_1,v_2$ and one bottom vertex $v_0$, where $v_i$ is joined to $v_0$ by a vertical edge $e_i$. Suppose  $v_i$ carries a marked pole $p_i$ whose aligned endpoint lies on the level of $v_0$, and $v_0$ carries a marked pole $p_0$ with nonzero residue at its own level. Let $\rho_i$ be the lower-side residue of $e_i$. 

\begin{center}
\begin{tikzpicture}[
    x=1.3cm, y=1.25cm,
    vertex/.style={circle,fill=black,inner sep=1.7pt},
    endpoint/.style={circle,fill=black,inner sep=1.25pt},
    every node/.style={font=\small}
]
  \node[vertex,label=above:$v_1$] (v1) at (-1.35,1.5) {};
  \node[vertex,label=above:$v_2$] (v2) at ( 1.35,1.5) {};
  \node[vertex,label=below:$v_0$] (v0) at (0,0) {};

  \draw (v1) -- node[pos=.48,above right] {$e_1$} (v0);
  \draw (v2) -- node[pos=.48,above left] {$e_2$} (v0);

  \node[anchor=east] at (-0.23,0.23) {$\rho_1$};
  \node[anchor=west] at ( 0.23,0.23) {$\rho_2$};
  \draw[dashed] (v1) -- (-2.25,0);
  \draw[dashed] (v2) -- ( 2.25,0);
  \node[endpoint] at (-2.25,0) {};
  \node[endpoint] at ( 2.25,0) {};
  \node[anchor=south east] at (-2.25,0.05) {$p_1$};
  \node[anchor=south west] at ( 2.25,0.05) {$p_2$};
  \draw[dashed] (v0) -- (0.9,0);
  \node[endpoint] at (0.9,0) {};
  \node[anchor=south] at (0.9,0.05) {$p_0$};
\end{tikzpicture}
\end{center}

The A-GRC on the two upper branches gives
\[
  \rho_1 - r_{p_1} = 0,\quad \rho_2 - r_{p_2} = 0.
\]
The residue theorem on the bottom vertex gives
\[
  \rho_1 + \rho_2 + r_{p_0} = 0.
\]
Eliminating $\rho_1$ and $\rho_2$ gives
\[
  r_{p_1} + r_{p_2} + r_{p_0} = 0.
\]
The last equation appears as the limit of the residue theorem equation on nearby smooth differentials after dividing all residues by the scale of the bottom level. Without scaling, the residues of $p_1, p_2, p_0$ all vanish in the limit.
\end{ex}

The A-GRC follows from the same argument for obtaining the ordinary GRC. On a nearby smoothing, the vanishing cycles corresponding to the edges in $E(Y,L)$ and to the aligned polar legs in $P(Y,L)$, together with vanishing cycles descending to lower levels whose contributions are negligible at the scale of $L$, bound a subsurface. Dividing the boundary integral by the scale of $L$ and passing to the limit gives \eqref{eq:A-GRC}. 
\par 
We introduce the following lemma to treat the converse.

\begin{lem}
\label{lem:modifyingdiff}
Fix a level graph $\widetilde{\Gamma}$ for residue-aligned multi-scale differentials.  Assume that the A-GRC equation holds for every level cut $(Y,L)$ for which it is imposed. Then the bottom up smoothing construction admits a compatible family of modifying differentials whose upper side residue at a vertical edge $e$ is $\rho_{e}$ and whose residue at a marked pole $p_i$ aligned to $L$ is $r_{p_i}$. The modifying differential at a given level depends only on the data prescribed on lower levels.
\end{lem}

\begin{proof}
We follow the construction of modifying differentials in \cite[Proposition 9.3]{BCGGM2}, adding the aligned marked poles to the puncture set at the level where their endpoints lie. The required residues at the vertical edges and at these marked poles specify the corresponding periods. For every connected upper region, the period extension criterion \cite[Proposition 9.2]{BCGGM2}
 states that this residue assignment extends to the required period homomorphism exactly when
 $$\sum_{e\in E(Y,L)}\rho_e - \sum_{p_i\in P(Y,L)}r_{p_i} = 0,$$
 that is, the A-GRC equation. With the same choices of the auxiliary subspace and nested multicurves as in that construction, \cite[Proposition 9.3]{BCGGM2} produces the desired modifying differential.  Then \cite[Corollary 9.4]{BCGGM2} shows its level $L$ component depends only on the lower level differential and scaling parameters. Applying this recursively from the lowest level up thus  produces a compatible family at all level cuts. The added terms have at worst simple poles at the aligned marked poles and hence supply the prescribed residue at the required scale without  changing the prescribed higher-order principal part. 
 \end{proof}

The following lemma shows that residue-aligned multi-scale structures with A-GRC are compatible with  specialization. 

\begin{lem}
\label{lem:contractcompat}
Let $\widetilde{\Gamma}\to\widetilde{\Gamma}'$ be a contraction of residue-aligned level graphs of multi-scale differentials. Then a residue-aligned multi-scale structure satisfying the A-GRC on $\widetilde{\Gamma}$ specializes to one satisfying the A-GRC on $\widetilde{\Gamma}'$. 
\end{lem}

\begin{proof}
A graph contraction forgets intermediate levels and components of the underlying curve.  Consider a cut $(Y',L')$ of $\widetilde{\Gamma}'$. If an aligned polar endpoint lies strictly above $L'$, no A-GRC is imposed for this level cut, just as in the ordinary GRC. Otherwise, suppose before the contraction $Y'$ is represented by a union of components separated by intermediate cuts. Add the corresponding A-GRC equations, together with the ordinary residue theorem on any components that are contracted between these cuts. The residue of an internal vertical edge then appears once from each side, with opposite signs, and hence cancels. What remains are exactly the residues of the vertical edges crossing the cut $(Y',L')$, together with the aligned marked pole residues at level $L'$, with the same sign convention as in \eqref{eq:A-GRC}. Thus the resulting equation is precisely the A-GRC for $(Y', L')$. 
\end{proof}

Summarizing the above discussions, the A-GRC is both necessary and sufficient for differentials in the interior of a stratum degenerating to a residue-aligned multi-scale differential. We denote by
\[ 
\RAMS(\mu)
\]
the resulting {\em moduli space of residue-aligned multi-scale differentials}. The above description of $\RAMS(\mu)$ implies that its local structure is similar to that of $\MS(\mu)$, with the only addition from the aligned residue coordinates satisfying the A-GRC. In particular, $\RAMS(\mu)$ enjoys the same property as that of $\MS(\mu)$: it is a {\em complex orbifold with normal crossing boundary}.   

\subsection{The blowup description}

We next relate the above modular description to a global blowup of $\MS(\mu)$. Recall first the local functions entering the blowup description of $\MS(\mu)$. Let $\Gamma$ be the enhanced level graph of a boundary stratum of $\MS(\mu)$. Locally every vertex $v\in V(\Gamma)$ has an adjusting parameter $h_v$, characterized by the property that $h_v^{-1}\omega$ has a nonzero limit differential on the component corresponding to $v$. The parameter $h_v$ is well defined up to multiplication by a unit. The global blowup description of $\MS(\mu)$ is obtained by the weighted principalization of the product of ideals $(h_v,h_{v'})$, as in the logarithmic and multi-scale comparison construction of \cite[Section~7]{CGHMS}, which globalizes the local blowup construction of \cite[Section~14]{BCGGM2}.

For each marked pole $p_i$, let $f_i$ be the residue function at $p_i$. For a vertex $v$, let 
$$ w(v) = 2g(v) -2 + \val(v), $$
where the valence $\val(v)$ counts half-edges and marked legs adjacent to $v$. In a local chart, let 
\begin{equation}\label{eq:Jres}
  J_{\res}(\Gamma) = \left(\prod_{1\leq i < j\leq p}(f_i, f_j)\right)
  \left(\prod_{i=1}^p\prod_{v\in V(\Gamma)}(f_i, h_v)^{w(v)}\right).
\end{equation}
Note that ideal factors that are already principal, or in which one generator is a unit, do not affect the blowup of the total product ideal.

\begin{lem}
\label{lem:Jres-glues}
The local ideals $J_{\res}(\Gamma)$ glue to a global ideal sheaf $J_{\res}$ on $\MS(\mu)$.
\end{lem}

\begin{proof}
The local ideals are unchanged by the choices of adjusting parameters, since these choices multiply the generators by units. It remains to check compatibility under the degeneration of level graphs. As in the global blowup description of the space of multi-scale differentials \cite[Section~7]{CGHMS}, it suffices to check the contraction of one edge, in which case the proof is identical to \cite[Proof of Lemma 7.2]{CGHMS}.
\end{proof}

The blowup above is the residue-aligned analogue of the orderly
blowup construction for the moduli space of ordinary multi-scale differentials
\cite[Theorem~14.8]{BCGGM2}. In a multi-scale local chart, the
adjusting parameters $h_v$ are already totally ordered by divisibility.
The additional factors of type 
$$
(f_i, f_j),\quad (f_i, h_v)
$$
simply insert the residue functions $f_i$ into this divisibility
order.  The resulting blowup therefore
records, for every marked pole, the first scale at which its residue function has a nonzero leading term.  This is precisely the residue-aligned structure in Definition~\ref{def:residue-aligned-MS}. The weights in \eqref{eq:Jres} are used to globalize this local construction under contractions of level graphs, as in the ordinary multi-scale case; see \cite[Section~7]{CGHMS}. 

\begin{prop}[Modular meaning of the blowup]
\label{prop:RAMS-blowup}
There is a natural isomorphism
\[ \RAMS(\mu)\cong \widetilde{\Bl}_{J_{\res}}\MS(\mu),\]
where the latter stands for the normalization of the blowup. In particular, the coarse moduli space of $\RAMS(\mu)$ is projective.
\end{prop}

\begin{proof}
The proof follows the same argument using the universal property of the blowup description of the space of multi-scale differentials as in \cite[Section~7]{CGHMS}. We explain the main points.

Let $T$ be a family of residue-aligned multi-scale differentials. Forgetting the positions of the marked pole endpoints gives a family of ordinary multi-scale differentials, hence a morphism $T\to\MS(\mu)$. On a residue-aligned local chart, the scaling functions for the components and residues are ordered by divisibility. Thus each of the ideals $(f_i, f_j)$ and $(f_i, h_v)$ occurring in \eqref{eq:Jres} pulls back to a principal ideal. Hence $J_{\res}$ pulls back to an invertible ideal, and the universal property of the blowup gives the morphism from $\RAMS(\mu)$ to the blowup, factoring through its normalization using the normality of $\RAMS(\mu)$. 

Conversely, on the normalized blowup, pull back the universal multi-scale differential from $\MS(\mu)$. The blowup of the product of the ideals in \eqref{eq:Jres} makes the residue functions and the level adjusting parameters become comparable by divisibility. For each marked pole $p_i$, this divisibility order determines the suitable scale at which the scaled residue has a finite nonzero leading coefficient. Aligning the endpoint of $p_i$ at that scale's level and recording this coefficient gives the aligned residue $r_{p_i}$. Hence one obtains a family of residue-aligned multi-scale differentials over the normalized blowup. (Here the normality is needed for the existence of level adjusting parameters, cf. \cite[Proposition 11.13]{BCGGM2}.)  
\par 
The two constructed maps in the above are clearly inverse to each other, which shows the desired isomorphism. Finally, projectivity follows from the blowup description. The coarse space of $\MS(\mu)$ is projective, the blowup of a projective variety along a global ideal sheaf is projective, and the normalization is finite over the blowup which preserves projectivity.
\end{proof}

\subsection{Extension of the residual map and properties of \texorpdfstring{$\RAMS(\mu)$}{RAMS(mu)}}
\label{subsec:extension}

It remains to define the extension of the residual map over the residueless loci. Let $x\in \RAMS(\mu)$, and let $L_{\max}$ be the highest level of its residue-aligned graph containing at least one marked pole endpoint with nonzero aligned residue. Define
\[
  \PRes(x) = [s_1: \cdots : s_p]\in \PP \calR_{p},
\]
where
\[
  s_i=\begin{cases}
    r_{p_i}, & \text{if the aligned endpoint of }p_i\text{ lies on level }L_{\max},\\
    0, & \text{if the aligned endpoint of }p_i\text{ lies on a level below }L_{\max}.
  \end{cases}
\]
Recall that the aligned residues $r_{p_i}$ are the leading scaled limits of the ordinary residue coordinates of nearby smooth differentials. Since the marked poles aligned below $L_{\max}$ have residues vanishing faster than those on $L_{\max}$, the residue theorem applied to a nearby smooth differential and divided by the scale of $L_{\max}$ gives $\sum_{i=1}^p s_i=0$. Hence the above expression defines a morphism
\begin{align}\label{eq:residual-extension}
  \PRes\colon \RAMS(\mu)\to \PP \calR_{p}.
\end{align}
On the open stratum of smooth differentials with at least one nonzero residue,  it agrees with the usual projectivized residual map. On the boundary it records the degeneration directions of the projective residue tuples after comparing the residue functions between themselves (see Example~\ref{ex:residueless}). 

\begin{rmk}\label{rmk:minimal-resolution}
In general, $\RAMS(\mu)$ is a refined resolution for the residual map. To achieve a minimal resolution, it suffices to retain the aligned residues on $L_{\max}$ and ignore those in deeper levels. In terms of the blowup description, one can alternatively use the ideal factor $(f_1, \ldots, f_p)$ in $J_{\res}$ in~\eqref{eq:Jres}. This comparison is analogous to the \emph{incidence variety compactification} \cite{BCGGM} versus the (ordinary) multi-scale compactification (also called \emph{partially aligned} versus \emph{fully aligned} in logarithmic geometry), where in the former the lower-level vertices are only partially ordered. 
\end{rmk}

The following statement summarizes some properties of $\RAMS(\mu)$ that will be used later. Its analogue for $k$-differentials will be given in Proposition~\ref{prop:kRAMS}.

\begin{prop}[Properties of $\RAMS(\mu)$]
\label{prop:RAMS-local-properties}
The moduli space of residue-aligned multi-scale differentials $\RAMS(\mu)$ satisfies the following properties: 
\begin{enumerate}[label=\textup{(\roman*)}]
\item The projectivized residual map extends to the morphism~\eqref{eq:residual-extension} defined above.
\item $\RAMS(\mu)$ admits local coordinates similar to those of $\MS(\mu)$ together with aligned residues of marked poles. 
\item $\RAMS(\mu)$ admits a Whitney stratification given by 
locally closed boundary strata with fixed residue-aligned level graphs. 
\end{enumerate}
\end{prop}

\begin{proof}
The morphism in (i) was constructed above from the leading nonzero 
coefficients of the aligned residue functions. More precisely, the blowup in Proposition~\ref{prop:RAMS-blowup} compares among the residue functions of the marked poles. Hence, each marked pole has a well-defined scale at which its aligned residue has a nonzero leading coefficient (see e.g. Example~\ref{ex:residueless}). These leading coefficients as aligned residues on the highest residue-alignment level thus extends the projectivized residual map as in \eqref{eq:residual-extension}.  

The coordinate and stratification statements in (ii) and (iii) then follow from the same descriptions for the multi-scale local charts as in \cite[Sections~9.2 and~13.3]{BCGGM2} and \cite[Section 7]{CGHMS}. In particular, normal crossing boundary implies Whitney. 
\end{proof}

\subsection{Locally trivial isoresidual fibrations}

In this section we prove local triviality of the residual map over any open flat of the resonance arrangement, for any moduli space $\mathcal{H}(\mu)$ in any genus, using the residue-aligned multi-scale compactification introduced above. This resonance stratification is the one already appearing in the study of isoresidual fibrations; see \cite[Sections~2.7 and~3.3]{CGPT}. 
\par
By Proposition~\ref{prop:RAMS-local-properties}, the projectivized residual map extends to a morphism
\[
  \PRes\colon \RAMS(\mu)\to \PP \calR_{p},
\]
where
\[
  \calR_{p} = \left\{(r_1, \ldots, r_p)\in\mathbb C^p\ \middle|\ \sum_{i=1}^p r_i = 0\right\}
\]
is the residual space of the marked poles. Locally on $\RAMS(\mu)$,  the boundary strata possess the residue coordinates satisfying the A-GRC equations. For the argument below, we use only the residue part of the aligned multi-scale coordinate system.

By an abuse of notation that will be clear from the context, we denote by $B = B(E)$ an open flat of $\calR_{p}$ (see Definition~\ref{defn:resonanceintro}) and its projectivization in $\mathbb{P}\calR_{p}$, where $E$ is a family of resonance hyperplanes. 

\begin{thm}[Local triviality of the isoresidual fibration]
\label{thm:RAMS-stratified-submersion}
For every signature $\mu$ and every open flat $B\subset \PP \calR_{p}$ of the projectivized resonance arrangement, the restriction
\[
  \PRes\colon \PRes^{-1}(B)\to B
\]
is a stratified submersion, where the stratification in the source is induced by the loci with fixed residue-aligned level graphs. In particular, the compactified projectivized isoresidual fibers in $\RAMS(\mu)$ are locally topologically trivial over every open flat $B$ of the resonance arrangement. Moreover, the same conclusion holds for the restriction of $\PRes$ to the complement of any union of closed boundary strata in $\RAMS(\mu)$.
\end{thm}

\begin{proof}
Write $B = B(E)$ and fix a residue-aligned multi-scale differential 
\[
x\in D_{\widetilde\Gamma}\cap\PRes^{-1}(B),
\]
where $D_{\widetilde\Gamma}$ is the stratum with fixed residue-aligned level graph.  Let $L$ be the highest level to which some marked pole is aligned. Marked poles aligned to $L$ map to their aligned residues under $\PRes$, while those aligned to lower levels map to zero.

Since $\PRes(x)\in B(E)$, each residual coordinate forced to vanish in this way vanishes identically on $E$; otherwise $\PRes(x)$ would lie on an additional resonance hyperplane not in $B(E)$. Since there are no aligned residues above $L$, it suffices to consider the A-GRC at the cut level $L$ together with the residue theorem for the vertices at level $L$.

For a connected component $Y$ above $L$, rewrite the A-GRC as 
\[
  \sum_{p_i\in P(Y,L)} s_i 
  -  \sum_{e\in E(Y,L)} \rho_e = 0,
\]
where $\rho_e$ is the residue on the lower side of the vertical edge $e$. Take a connected block of the part of the level graph consisting of the vertices at level $L$, the components $Y$ above $L$, and the vertical and horizontal edges joining them. Adding the corresponding A-GRC and residue theorem equations cancels the residue of every vertical edge; the two residues at each horizontal edge cancel by the local matching-residue condition. Hence we obtain
\[
  \sum_{i\in I} s_i = 0
\]
for the marked poles belonging to this connected block.

On the other hand, no other relation among the $s_i$ can arise. Indeed, suppose a linear combination of these equations cancels all residues at the nodes. At each vertical node, and at each horizontal node joining two vertices, the residue occurs with opposite signs in the two corresponding equations, so those two equations must have the same coefficient. Since the concerned block is connected, all coefficients are equal. Thus, up to a scalar, the only such combination is the sum of all the equations, and the only resulting relation among the $s_i$ is the partial sum above. 

We thus conclude that every relation imposed on the $s_i$ is a resonance relation. Since $\PRes(x)\in B(E)$, each of these relations already vanishes identically on $E$; otherwise $\PRes(x)$ would lie on an additional resonance hyperplane. Hence there is no further condition on the residue coordinates along $B(E)$.

The local period and residue coordinates on $D_{\widetilde\Gamma}$ then imply the surjectivity of 
$$  d\PRes\colon  T_xD_{\widetilde\Gamma} \to T_{\PRes(x)}B(E) $$ on the tangent spaces.  
 Hence $\PRes$ is a stratified submersion over $B$.

By the local normal crossings structure of $\RAMS(\mu)$, the stratification by residue-aligned level graphs is Whitney. Moreover, $\PRes^{-1}(B)\to B$ 
is proper, since the extended residual map \eqref{eq:residual-extension} is proper with both domain and target being compact. Thom's first isotopy lemma~\cite[Proposition~11.1]{Mather} thus gives, for every $b\in B$, a sufficiently small neighborhood $U\subset B$ and a stratum preserving homeomorphism
\[
\PRes^{-1}(U)\cong  \PRes^{-1}(b)\times U.
\]
Hence the compactified  projectivized isoresidual fibers are locally topologically trivial over $B$.

Finally, since the trivialization preserves the boundary stratification, it restricts to the  complement of any union of closed boundary strata.
\end{proof}

\subsection{The isoresidual fibration for \texorpdfstring{$k$}{k}-differentials} \label{subsec:k-diff}

Theorem~\ref{thm:RAMS-stratified-submersion} extends from meromorphic one-forms to $k$-differentials. Rather than using directly the hypersurface equations defining the $k$-global residue condition on the original stratum, we pull  back to the canonical $k$-cover, where the relevant $k$-residues have chosen $k^{th}$ roots. There the $k$-global residue conditions become ordinary linear global residue conditions,  as in \cite{BCGGM3} and \cite{chen2025finite}, and the proof thus reduces to the case $k=1$.

Fix $k\geq 2$ and a moduli space $\mathcal{H}^{k}(\mu)$ of meromorphic $k$-differentials of type $\mu$ that are not global powers of differentials of lower order (they are called \textit{primitive differentials}). Let $p$ be the number of marked poles whose orders are {\em divisible by $k$}. These and only these marked poles carry possibly nonzero $k$-residues. Let 
\[
  \calR_{p}^k\cong \mathbb C^p, \quad R = (R_1,\ldots, R_{p}),
\]
for the $k$-residue space. Unlike the abelian case, in general, no residue theorem imposes a linear relation among the coordinates $R_i$ for primitive $k$-differentials with $k\geq 2$. Let 
\[
\mathbf{Res}^k\colon \mathcal{H}^{k}(\mu) \to \calR_{p}^k
\]
 be the \emph{$k$-residual map} and $\PRes^k$ its projectivization. As before, we consider the case when the image of $\mathbf{Res}^k$ is not $\{ 0 \}$ so the projectivization is well defined. 

We will use the following analogue of the residue-aligned multi-scale compactification for the stratum $\mathcal{H}^{k}(\mu)$ of $k$-differentials. It can be constructed by applying residue alignment on the canonical cyclic $k$-cover, where the chosen $k^{th}$-root differential is a (normalized) abelian differential, and then descending by the deck group of the canonical cover, in parallel with the construction for multi-scale $k$-differentials as in \cite{BCGGM3} and \cite{CMZarea}.

\begin{prop}[Properties of $\RAMS^k(\mu)$]
\label{prop:kRAMS}
The construction of the residue-aligned space $\RAMS(\mu)$ generalizes to $k$-differentials. It gives a compactification
\[
  \RAMS^k(\mu)
\]
of the stratum of $k$-differentials of type $\mu$ 
with the following properties: 
\begin{enumerate}[label=\textup{(\roman*)}]
\item The projectivized $k$-residual map  extends to a morphism
\[
  \PRes^k\colon \RAMS^k(\mu) \to \mathbb{P} \calR_{p}^k.
\]
\item $\RAMS^k(\mu)$ admits local coordinates similar to those of $\MS^k(\mu)$ together with chosen $k^{th}$-root aligned residues of marked poles. 
\item $\RAMS^k(\mu)$ admits a Whitney stratification given by locally closed boundary strata with fixed residue-aligned $k$-multi-scale structures. 
\end{enumerate}
\end{prop}

\begin{proof}
The construction for $\RAMS(\mu)$ applies to the multi-scale compactification of strata of $k$-differentials $\xi$ after passage to the canonical cyclic cover $\pi$. On that cover the chosen $k^{th}$-root differential $\widehat\omega$ is an abelian differential satisfying $\widehat\omega^k=\pi^*\xi$.  The multi-scale compactification for $k$-differentials is obtained by imposing the deck-equivariant matching, prong, and $k$-GRC data \cite[Sections~2, 4, 5]{BCGGM3} and \cite[Sections~2--4]{CMZarea}.

After passing to the $k^{th}$-root residue cover, the step of aligning residues is identical to the abelian case. If the chosen $k^{th}$-root residue of a marked pole has zero leading term at a given level, the marked  polar leg is aligned to the lower level at which the leading $k^{th}$-root residue is nonzero (or to the lowest level with zero aligned residue value if the corresponding $k$-residue function at this pole is identically zero on the entire stratum). Equivalently, one aligns the leading term of a chosen $k^{th}$-root residue with the level adjusting parameters. The same blowup and level-recording procedure gives the extension of $\PRes^k$ after taking $k^{th}$ powers on the original stratum.

The aligned $k$-GRC can be formulated on the canonical cover. Once the branch of roots as well as the criss-cross and abelian-type level graph description for nodal canonical covers have been fixed, the residue equations are linear in the chosen $k^{th}$-root residue variables. With the sign convention of \eqref{eq:A-GRC}, a marked pole aligned to a lower level contributes with the sign opposite to that of the $k^{th}$-root residue at the lower side of a vertical edge. Thus, for a level $L$ and a connected component $Y$ of the subgraph above $L$, such that $Y$ is of non-abelian type and does not contain a pole of order at least $k$ (see \cite[Definition 1.4]{BCGGM3}), the local aligned $k$-GRC equations have the form
\[
  \sum_{e\in E(Y,L)} \epsilon_e q_e - \sum_{p_i\in P(Y,L)} \epsilon_i r_i = 0,
  \quad {\rm for}~\epsilon_e,\epsilon_i\in\mu_k,
\]
up to the fixed coefficients determined by the deck group action. These equations correspond to the residue theorem conditions for constructing  the modification differentials on the covering curve, so they give the smoothability criterion.  
\end{proof} 

If $p = 1$, the projectivized target of the $k$-residual map is a point,  and Theorem~\ref{thm:kRAMS-submersion} below holds trivially. From now on we treat only the case $p\geq 2$. Let
\[
  \widetilde{\calR}_{p} = \mathbb{C}^{p}, \quad
  \varpi\colon \widetilde{\calR}_{p}\to \calR_{p}^{k}, \quad
  (r_1,\ldots,r_p)\mapsto (r_1^k,\ldots,r_p^k)
\]
be the finite \emph{$k^{th}$-root residue cover}. We also use $\varpi\colon \mathbb{P} \widetilde{\calR}_{p}\to \PP \calR_{p}^k$ to denote its projectivization. 
\par 
On $\mathbb{P}\widetilde{ \calR}_{p}$ consider the \emph{$k^{th}$-root resonance arrangement}
\begin{align*}
  \widetilde{A}_{k,p} =\left\{
  H_w = \left\{r\in\widetilde{\calR}_{p}\ \middle|\ \sum_{i=1}^p w_i r_i = 0\right\}
  \ \middle|\
  w = (w_1,\ldots, w_p)\in (\{0\}\cup\mu_k)^p\setminus\{0,\ldots, 0\}
  \right\},
\end{align*}
where $\mu_k$ is the set of $k^{th}$-roots of unity. 

Let $\widetilde{E}$ be an open flat with respect to the $k^{th}$-root resonance arrangement in $\widetilde{A}_{k,p}$. Define 
\[
B(\widetilde E) \coloneqq \varpi (\PP \widetilde E)
\]
to be a projectivized {\em $k$-resonance stratum} in $\PP \calR_{p}^k$. Changing the choices of $k^{th}$-roots of unity preserves the $k^{th}$-root resonance arrangement $\widetilde A_{k,p}$. Since an open flat is determined by the resonance hyperplanes containing any of its points, the inverse image of $B(\widetilde E)$ under $\varpi$ is precisely the union of the open flats obtained from $\widetilde E$ by changing these roots in the residue coordinates. 

\begin{thm}[Local triviality of the $k$-isoresidual fibration]
\label{thm:kRAMS-submersion}
Let $B\subset \PP \calR_{p}^k$ be a projectivized $k$-resonance stratum. 
Then
\[
  \PRes^k\colon (\PRes^k)^{-1}(B)\to  B
\]
is a stratified submersion, where the domain is stratified by the  residue-aligned $k$-multi-scale boundary strata as in Proposition~\ref{prop:kRAMS}.

Consequently, the compactified projectivized $k$-isoresidual fibers are locally topologically trivial over $B$. More precisely, for every $b\in B$ there is a sufficiently small neighborhood $U\subset B$ and a stratified homeomorphism
\[
  (\PRes^k)^{-1}(U)\cong (\PRes^k)^{-1}(b) \times U.
\]
The same conclusion holds after replacing $\RAMS^k(\mu)$ by any open subset obtained by removing a union of closed boundary strata.
\end{thm}

\begin{proof}
We only explain what changes from the case $k=1$. Let $\widetilde E$ be an open flat lying over $B$. Choose $k^{th}$-root residue coordinates $s_i$ on a local chart. Fix a domain stratum over $\widetilde E$ and a point $x$ on it, and let $L$ be the highest level to which some marked pole is aligned in the level graph of $x$.

Marked poles aligned to levels below $L$ map to zero under the extended residual map. Since the image of $x$ lies in the open flat $\widetilde E$, all such forced vanishing relations already hold on $\widetilde E$. Since there are no aligned residues above $L$, only the residue equations at level $L$ can impose conditions on the varying target residues.

On the canonical cover, the aligned $k$-GRC at level $L$ is the A-GRC for the corresponding multi-scale abelian differential.  The proof of Theorem~\ref{thm:RAMS-stratified-submersion} for $k=1$ therefore applies without change at this level: combine the A-GRC
equations above $L$ with the ordinary residue theorem on the vertices at level $L$. The residues of vertical nodes cancel between these equations, while the residues of horizontal nodes cancel by the matching-residue condition. Hence the remaining relations among the marked pole residues are sums of residues on the canonical cover. By the same argument as in the case $k=1$, there are also no other relations on the target residues.

In the chosen $k^{th}$-root residue coordinates, each such relation has the form 
\[
\sum_{i\in I}\zeta_i s_i = 0, \quad \zeta_i\in\mu_k.
\]
These are exactly the $k^{th}$-root resonance hyperplanes defining the $k$-resonance stratification. Since $x$ maps to $\widetilde B = B(\widetilde E)$, every relation that occurs must already vanish identically on
$\widetilde E$; otherwise $x$ would lie on an additional
$k^{th}$-root resonance hyperplane that is not in $\widetilde E$. Thus there is no additional condition on the tangent directions to $\widetilde E$. The local period and residue coordinates show that the lifted residual map is a submersion on every domain boundary stratum.

The domain and target strata are finite quotients of these $k^{th}$-root residue charts, as in the construction of $\RAMS^k(\mu)$. On each open flat, the map taking $k^{th}$-powers is locally unramified, since the vanishing residue coordinates are fixed. Hence this gives the desired stratified submersion for
\[
  \PRes^k\colon (\PRes^k)^{-1}(B)\to B. 
\]

The rest is the same as in the case $k=1$. The stratification given by the boundary strata is Whitney, and the restriction of the $k$-residual map over $B$ is proper. Thom's first isotopy lemma thus gives the desired local product over $B$.  Since the trivialization preserves the boundary strata, it also restricts to the complement of any union of closed boundary strata.
\end{proof}

\section{Translation structures and metric completions}
\label{sec:wysiwyg-metric-completion}

In this section, we show how isoresidual fibers, and more generally {\em leaves of the isoperiodic foliation}, can be endowed with holomorphic coordinates for which the transition functions are translations (Proposition \ref{prop:translation}). The resulting flat path metric has an intrinsic completion. Proposition~\ref{prop:local} identifies this completion locally with certain WYSIWYG compactification (cf. \cite{WYSIWYG}) and gives its finite holomorphic map to the model space $V_n$ of relative periods (cf. Definition~\ref{def:arrangement}).

\subsection{Exact matrices and collision arrangement}\label{sub:ExactMatrices}

As we will see in Section~\ref{sub:Translation}, isoperiodic loci are locally modeled on configuration spaces of $n$ points in the complex plane up to translation. We introduce a class of matrices that naturally encodes this structure.

\begin{defn}[Exact matrix]
\label{defn:ExactMatrix}
Fix an integer $n \ge 1$. An $n \times n$ matrix $M = (M_{ij})$ with coefficients in $\mathbb{C}$ is called \emph{exact} if there exists a vector $a = (a_1,\dots,a_n) \in \mathbb{C}^n$ such that
\[
M_{ij} = a_j - a_i \quad \text{for all } 1 \le i,j \le n.
\]
We denote by $\mathcal{E}_n(\mathbb{C})$ the set of exact $n \times n$ matrices with coefficients in $\mathbb{C}$.
\end{defn}

The set $\mathcal{E}_n(\mathbb{C})$ is a vector subspace of $\mathbb{C}^{n \times n}$. Moreover, the map
$$
\mathbb{C}^n \to \mathcal{E}_n(\mathbb{C}), \quad a \mapsto (a_j - a_i)_{1 \le i,j \le n}
$$
is linear and has kernel consisting of constant vectors. Hence, $\dim \mathcal{E}_n(\mathbb{C}) = n - 1$.

\begin{defn}[Collision arrangement]
\label{def:arrangement}
For $n \geq 2$, let $\Delta \subset \mathbb{C}^n$ denote the {\em diagonal subspace}
\[
\Delta = \{(z,\dots, z) \in \mathbb{C}^n \mid z \in \mathbb{C}\}.
\]
In the quotient space 
\[
V_{n} \coloneqq \mathbb{C}^n / \Delta,
\]
the \emph{collision arrangement} is induced by the \emph{polydiagonal hyperplanes}
\[
B_{ij} = \{(\zeta_1,\dots,\zeta_n) \in \mathbb{C}^n \mid \zeta_i = \zeta_j\}, \quad 1 \le i < j \le n.
\]
Since each $B_{ij}$ contains $\Delta$, it descends to a well-defined hyperplane arrangement in $V_{n}$, still denoted by $$\mathcal{B}_{n}=\bigcup\limits_{1 \le i < j \le n}B_{ij} \subset V_{n}.$$ 
This arrangement can also be described as the \emph{braid arrangement modulo translations}.

Moreover, the space $V_{n}$ admits a natural identification with the space $\mathcal{E}_n(\mathbb{C})$ of exact matrices. Indeed, the linear map
\[
\mathbb{C}^n \to \mathcal{E}_n(\mathbb{C}), \quad a \mapsto (a_j - a_i)_{1 \le i, j \le n}
\]
factors through the quotient by $\Delta$ and induces a linear isomorphism
\[
V_{n} \cong \mathcal{E}_n(\mathbb{C}).
\]
Under this identification, the hyperplane $B_{ij}$ corresponds to the linear condition $M_{ij} = 0$, so that the arrangement $\mathcal{B}_{n}$ depends only on the differences $\zeta_i - \zeta_j$ for $1 \le i < j \le n$. 
\end{defn}
The arrangement $\mathcal{B}_{n}$ in $V_{n}$ is {\em essential} (i.e.\ the intersection of all hyperplanes in $\mathcal{B}_{n}$ is $\{0\}$). It is the classical hyperplane arrangement of type $A_{n-1}$.

\subsection{Translation structures on isoperiodic leaves}\label{sub:Translation}

Consider a moduli space $\mathcal{H}(\mu)$ of meromorphic one-forms for a partition $\mu=(a_{1},\dots,a_{n},-b_{1},\dots,-b_{p})$ of $2g-2$. The \textit{isoperiodic foliation} is defined in $\mathcal{H}(\mu)$ by fixing the (absolute) periods of the differentials along closed loops on the surface punctured at the poles; the set of differentials realizing a fixed set of absolute periods is called an {\em isoperiodic leaf}, and will be denoted $\mathcal{L}$. Note that in genus zero, every absolute period is a sum of residues, and so in this case an isoperiodic leaf is the same as a connected component of an isoresidual fiber. The result below shows that the system of relative periods endows any isoperiodic leaf with the structure of a {\em translation manifold}.
\par
\begin{prop}[Translation structure on an isoperiodic leaf]
\label{prop:translation}
    Let $\mathcal{L} \subset \calH(\mu)$ be an isoperiodic leaf. Then there exists an atlas of holomorphic charts on $\mathcal{L}$ for which the transition functions are given by translations.
\end{prop}

\begin{proof}
    Recall that by convention, the zeros (and poles) of any differential $\omega \in \calH(\mu)$ are labeled: $z_{i}$ denotes the zero of order $a_{i}$ for $i=1,\dots,n$. Given $(X,\omega) \in \mathcal L$, let $\gamma$ be a homotopy class of paths on $X$ punctured at the poles of $\omega$. Such $\gamma$ has at least one geodesic representative, and so any homotopy class of paths connecting a pair $z_i, z_j$ of zeros admits a representative given by a chain of saddle connections. Consequently, there exists a tree $\Gamma$ spanning $z_{1},\dots,z_{n}$ that is formed by $n-1$ (oriented) saddle connections.
\par
In a neighborhood of $\omega$ in $\mathcal{L}$, the periods $\int_{\gamma} \omega$, where $\gamma$ are edges of $\Gamma$, give local complex coordinates for $\mathcal{L}$.
Consider any path $\alpha$ joining a pair of zeros $z_{i},z_{j}$ in $X$ punctured at the poles of $\omega$. Let $\beta$ be the broken geodesic from $z_{i}$ to $z_{j}$ contained in $\Gamma$. Then the difference $\int_{\alpha} \omega - \int_{\beta} \omega$ is the period of $\omega$ along the absolute homology class $[\alpha]-[\beta]$. It follows that these two local coordinates differ by a constant.
\par
Denote by $Z_{ij}$ the period of $\omega$ along an arbitrary path $\alpha$ from $z_{i}$ to $z_{j}$. The differential $dZ_{ij}$ does not depend on the choice of $\alpha$ and is thus a globally defined nowhere vanishing holomorphic one-form on $\mathcal{L}$.
For any $i, j, k$, we have $dZ_{ik} = dZ_{ij} + dZ_{jk}$, so that the differentials $(dZ_{ij})$ assemble into an  $\mathcal{E}_n(\mathbb{C})$-valued one-form $dZ$. Moreover, as relative periods $Z_{ij}$ are locally injective holomorphic functions on $\mathcal{L}$ and define a system of local coordinates on $\mathcal{L}$, the form $dZ$ has neither zeros nor poles in $\mathcal{L}$ and is of \textit{maximal rank} $n-1$ at every $[\omega] \in \mathcal{L}$.
\par
Conversely, a \textit{local primitive} of $dZ$ in $\mathcal{L}$ is an  $\mathcal{E}_n(\mathbb{C})$-valued injective holomorphic map. Such maps form an atlas of $\mathcal{L}$ where transition maps are translations of complex vector space $\mathcal{E}_n(\mathbb{C})$. Thus the matrix-valued differential $dZ$ defines a (higher-dimensional) translation structure on $\mathcal{L}$. 
\end{proof}

\begin{rmk}\label{rmk:translation-orbifold}
When a stratum $\mathcal H(\mu)$ is an orbifold, its period coordinates are initially defined on finite local covers of the stratum. They nevertheless descend to the underlying isoperiodic leaves.  Indeed, the stabilizer of an orbifold point fixes every labeled zero, and hence acts trivially on
$$
\ker\left( H^1(X\setminus P, Z; \mathbb C)
\to H^1(X\setminus P; \mathbb C) \right)
\cong H^0(Z;\mathbb C)/\mathbb C.
$$
Thus its action on the corresponding isoperiodic leaf has
trivial linear part.  Since the stabilizer is finite, this affine action is in fact trivial.  Consequently the relative  period coordinates descend to genuine translation charts on the underlying isoperiodic leaf.   
\end{rmk}

\subsection{The metric completion and WYSIWYG}
\label{subsec:local}

The above result endows every isoperiodic leaf $\mathcal{L}$ with a flat path metric induced by $V_{n} = \mathbb{C}^n / \Delta$. We denote the \emph{metric completion of $\mathcal L$} by $\mathcal{L}^{\mathrm{MET}}$. 

We next describe another completion through the contraction of multi-scale limit differentials when the marked zeros collide in the boundary of $\mathcal L$. Given a multi-scale differential, let the lower part be the union of the components strictly below the top level. Following the  construction of~\cite{MW17, WYSIWYG}, we contract (in the topological sense) each connected component $C_\alpha$ of the lower part to a point $Q_\alpha$, retain the top-level differential and every marked point which remains on the top level, and decorate $Q_\alpha$ by the set of all labeled zeros and poles carried by $C_\alpha$. We call the resulting quotient space \emph{WYSIWYG}. In terms of flat geometry, for a WYSIWYG element, the top level differential becomes identically zero on lower levels, so lower level subsurfaces disappear under the top level flat metric, and we fill in any punctures thus created with marked points. (In the genus-zero case the objects parameterized in WYSIWYG also appear later as {\em arboreal curves and differentials} in Definition~\ref{def:arboreal} for our flat-geometric approach.) 
\par 
Let $x$ be a limit multi-scale differential of an isoperiodic leaf $\mathcal L$. We first show that the lower part of $x$ discarded by the WYSIWYG contraction does not contribute to the actual periods.

\begin{lem}
\label{lem:lowerexact}
Every differential on any lower level of $x$ is exact.
\end{lem}

\begin{proof}
We work on a multi-scale local chart containing $x$. Let $L$ be a lower level in the level graph of $x$, let $t_L\to 0$ be its scaling parameter, and let $\eta_L$ be the scaled differential on level $L$. For an absolute cycle $\gamma$ on a vertex of level $L$, transport $\gamma$ to a
cycle $\gamma_t$ on the nearby smooth fiber. The construction of the multi-scale space gives
$$
\int_{\gamma_t}\omega_t =  t_L \int_\gamma\eta_L + o(t_L).
$$
The left-hand side is an absolute period and is therefore constant
along $\mathcal L$, whereas the right-hand side tends to zero.
Hence this constant is zero. Dividing by $t_L$ and passing to the limit
gives
$$
\int_\gamma\eta_L=0.
$$
The same argument can be applied to peripheral cycles around the lower level marked poles and to the vanishing cycles associated with vertical nodes. It follows that all residues of
$\eta_L$ vanish. Thus all absolute periods of $\eta_L$ vanish, and
$\eta_L$ is exact.
\end{proof}

Following Lemma~\ref{lem:lowerexact}, we say that a limit multi-scale differential $x$ of $\mathcal L$ is of \emph{collision type} if its level graph has no horizontal edges, where the additional requirement of no horizontal edges on the top level is due to our consideration of collisions of zeros rather than zeros moving infinitely far apart. In other words, in the metric completion $\mathcal L^{\mathrm{MET}}$ we only consider the collisions of marked zeros when some relative periods become arbitrarily short rather than arbitrarily long. 

Fix a nonzero collision-type multi-scale limit $x$ of $\mathcal L$. Let $q_0$ be the usual quotient map from the multi-scale closure of $\mathcal L$ to WYSIWYG, and let  
\[ 
y = q_0(x). 
\]
We work in a sufficiently small neighborhood of the fiber $q_0^{-1}(y)$ in the multi-scale closure of $\mathcal L$, where the \emph{normalization} of this neighborhood is denoted by $\mathcal{L}_y^{\mathrm{MS}}$, and for simplicity we still write $q_0$ for the quotient map from the normalized $\mathcal{L}_y^{\mathrm{MS}}$ to WYSIWYG. 
\par 
For $x_1, x_2\in \mathcal{L}_y^{\mathrm{MS}}$, if $q_0(x_1) = q_0(x_2) = y'$ and if $x_1, x_2$ belong to the same connected component of $q_0^{-1}(y')$, then we define $x_1\sim x_2$. Let $\mathcal L_y^{\mathrm{WYS}}$ be the quotient of $\mathcal L_y^{\mathrm{MS}}$ modulo this equivalence relation, and let 
\[ 
q\colon \mathcal L_y^{\mathrm{MS}} \to 
\mathcal L_y^{\mathrm{WYS}} 
\]
be the resulting quotient map. In other words, an element of $\mathcal L_y^{\mathrm{WYS}}$ consists of a usual WYSIWYG limit, together with a {\em choice of connected component} of its preimage in the normalized multi-scale closure. Retaining the information of a connected fiber component is for utilizing the {\em Stein factorization}, which will become clear in the proof of Proposition~\ref{prop:local} below. 

By Lemma~\ref{lem:lowerexact}, the (perturbed) period coordinates of \cite[Section 9]{BCGGM2} have no logarithmic terms on the  collision locus. Hence the actual relative period map extends holomorphically to
\[
\Psi\colon \mathcal L^{\mathrm{MS}}_y\to V_n
\] 
near $q_0^{-1}(y)$, and $\Psi$ is locally biholomorphic on the interior by Proposition~\ref{prop:translation}.

Let $\mathcal{L}_y^{\mathrm{MET}}$ denote a sufficiently small neighborhood in the metric completion $\mathcal{L}^{\mathrm{MET}}$ along the limits approached by families in $\mathcal L$ whose multi-scale preimages converge to a point of $q_0^{-1}(y)$.  

In what follows we will compare the two completions $\mathcal{L}_y^{\mathrm{WYS}}$ and $\mathcal{L}_y^{\mathrm{MET}}$. 

\begin{prop}[WYSIWYG and the metric completion]
\label{prop:local}
$\mathcal L_y^{\mathrm{WYS}}$ can be given a normal complex-analytic structure such that locally: 
\begin{enumerate}[label=\textup{(\roman*)}]
\item The actual relative period map $\Psi$ descends to a finite holomorphic map
\[
\Phi\colon \mathcal L_y^{\mathrm{WYS}}\to V_n.
\]
\item There is a natural identification 
\[
\mathcal L_y^{\mathrm{WYS}}\cong \mathcal L_y^{\mathrm{MET}}.
\]
\item The map $\Phi$ is locally biholomorphic away from the collision locus.
\end{enumerate}
\end{prop}

\begin{proof}
Since the quotient map $q$ retains the top-level differential, the map $\Psi$ is constant on every fiber of $q$. Conversely, we will show that $q$ is constant on every connected component of a fiber of $\Psi$.

Let $F$ be a connected component of a fiber of $\Psi$. On each fixed boundary stratum, the fixed absolute and actual relative periods determine the top level differential locally. Hence $q$ is locally constant on the intersection of $F$ with each such stratum. Further degeneration occurring in deeper levels does not affect $q$. So the only remaining possibility for $q$ to vary along $F$ is that a level passage opens and a lower-level component $C$ joins the top level. By Lemma~\ref{lem:lowerexact}, its differential is exact, say, $\eta_C=df_C$, where $f_C$ has at least two distinct critical points at marked zeros or holomorphic nodes on $C$. 
But then $C$ carries a non-zero relative period after scaling, say over the path $\gamma$ joining these two points. By extending the endpoints of $\gamma$ through the lower levels if necessary, we may choose $\gamma$ to represent a relative path between two marked zeros. Hence, when the level passage opens, the corresponding actual relative period has leading term
\[
t_C\int_\gamma\eta_C + o(t_C)
\]
where $\int _\gamma \eta_C\neq 0$ and $t_C$ is the scaling parameter for the level of $C$. It therefore cannot remain zero, which contradicts the assumption that $F$ lies in the fiber of $\Psi$. Hence $q$ is constant on $F$. Since $\Psi$ is constant on each fiber of $q$, this shows that the two maps $q$ and $\Psi$ have the same connected fiber components.

Consider the Stein factorization of $\Psi$. The preceding comparison of the connected components of the fibers of $\Psi$ and $q$ identifies the middle term of its Stein factorization with $\mathcal{L}_y^{\mathrm{WYS}}$: 
\begin{center}
\begin{tikzcd}
\centering
\mathcal{L}_y^{\mathrm{MS}} \arrow[rd, "q"] \arrow[rr, "\Psi"] & & V_n \\
&\mathcal{L}_y^{\mathrm{WYS}} \arrow[ru, "\Phi"] &
\end{tikzcd}
\end{center}
where $q$ has connected fibers and $\Phi$ is finite. Since $\mathcal{L}_y^{\mathrm{MS}}$ is normal, $\mathcal{L}_y^{\mathrm{WYS}}$ admits a compatible normal, complex-analytic structure for which $\Phi$ is finite and holomorphic.  

The interior leaf $\mathcal L$ is dense in $\mathcal{L}_y^{\mathrm{MS}}$, where the relative period map is a local isometry for the flat metric. Locally, a finite holomorphic map has only finitely many analytic branches, and normalization separates these branches. The flat metric pulled back from $V_n$ completes each such branch by adding its limiting point, while keeping distinct branches separate. Hence we have the desired identification
\[
\mathcal{L}_y^{\mathrm{WYS}}\cong \mathcal{L}_y^{\mathrm{MET}}.
\]
\par 
Away from the collision locus, no lower-level component is contracted, so $\Phi$ agrees locally with the relative-period map and is locally biholomorphic. \end{proof}

\begin{rmk}
In genus zero, the absolute period character is determined by the residues, so Proposition~\ref{prop:local} applies to every connected component of an isoresidual fiber at each limit of collision type.
\end{rmk}

\begin{rmk}
The same argument applies after imposing prescribed linear conditions on the relative periods, provided that these equations cut out a local analytic submanifold on which the restricted period map has constant rank (i.e. a locally linear submanifold in the ambient space). One
simply replaces $V_n$ by the corresponding affine linear subspace and defines the local WYSIWYG completion using the contraction from the multi-scale compactification of the linear submanifold (see \cite[Section 3]{cms}).
\end{rmk}

\begin{rmk}\label{rmk:zero-limit}
If the differentials parameterized in the isoperiodic leaf $\mathcal{L}$ are {\em exact}, global scaling by a constant preserves the isoperiodic leaf, and so the 
metric completion can also contain the identically-zero limit obtained from
$$ (X,\omega)\mapsto (X, t\omega), \quad   t\to 0.
$$
Nevertheless, since this special limit is identically zero, it is irrelevant to Proposition~\ref{prop:local}. 
\end{rmk}

\section{Polyhedral period domains}\label{sec:Polyhedral}

In the next two sections, we analyze the residual map in genus zero through the lens of flat geometry. Our objective here is to construct a period mapping that will realize an isoresidual fiber as a union of (unbounded) polyhedra glued together along polyhedral subregions of faces. The stratified fibration structure of the residual map will emerge naturally from the combinatorics: we will see that the isomorphism types of the polyhedra are locally constant along a given flat in the residue space.

\subsection{Main results}
As usual, $\mu = (a_1, \dots, a_n, -b_1, \dots, -b_p)$ is a partition satisfying $\sum_{i = 1}^n a_i - \sum_{j=1}^p b_j = -2$; we assume $n+p \ge 3$ throughout. $\mathcal{H}(\mu)$ denotes the stratum of meromorphic differentials on $\mathbb{CP}^{1}$ with zeros and poles of orders prescribed by $\mu$, and $\res \colon \mathcal{H}(\mu) \to \calR_{p}$ denotes the residual map. 
\par
Before giving the statement of the main results, we first formulate a key genericity condition. In this section and the next, it is more natural to consider periods as opposed to residues, and in particular to rescale residues $\lambda_i$ to the associated periods $2 \pi \mathrm{i} \lambda_i$.

\begin{defn}[Vertically generic]\label{def:vgen}
   A configuration of residues $\lambda = (\lambda_1, \dots, \lambda_p) \in \calR_p$ (or a differential $\omega \in \res^{-1}(\lambda)$) is said to be {\em vertically generic} if for all subsets $I \subset [p]$ such that $\Rea(2 \pi \mathrm{i} \lambda_I) = 0$, then moreover $\lambda_I = 0$, where  $\lambda_{I}=\sum\limits_{i \in I} \lambda_{i}$.
\par
   A subset $U \subset \calR_p$ is likewise said to be vertically generic if each $\lambda \in U$ is. 
\end{defn}

Throughout the next two sections, $B\subset \calR_p$ will denote an open flat (cf. Definition \ref{defn:resonanceintro}), and $U \subset B$ will be a connected vertically generic subset. We will see that such $U$ have precisely the properties we are looking for: we will be able to construct a family of period domains whose combinatorics are constant along such $U$.\\

\para{Structure of the period domain} The main result of Sections \ref{sec:Polyhedral} and \ref{sec:Period} is Theorem \ref{thm:stratfib}, which gives another proof of Theorem \ref{thm:MAIN} for strata of differentials on $\mathbb{CP}^1$.  This is proved via a study of a certain period map, defined on a subset $\res^{-1}(U)$ for $U \subset B \subset \calR_p$ as above. The target $\calP^U$ of the period map is a family over $U$ of what we call ``almost-polyhedral complexes'' $\calP^\lambda$ for $\lambda \in U$. A given $\calP^\lambda$ is a complex constructed by identifying polyhedra along {\em polyhedral subdomains} of their facets (by contrast, an ordinary ``polyhedral complex'' would identify polyhedra along entire facets). The constituent polyhedra $\calP^\lambda_\Gamma \subset \calP^\lambda$ are indexed by the data of a {\em critical graph} $\Gamma$, which tracks the combinatorial structure of the vertical foliations of the differentials (see Section \ref{SS:critgraph}). \\

\para{Statement of main results} In order to prove the main result Theorem \ref{thm:stratfib}, we will establish three intermediate results about these period domains and the associated period maps, including Theorem \ref{lem:Pproduct} as stated in the Introduction.\\

\noindent{\textbf{Theorem \ref{lem:Pproduct}} (Local combinatorics are invariant along a flat)\textbf{.}} \textit{
   Let $\Gamma$ be any critical graph, let $B \subset \calR_p$ be an open flat, and let $U\subset B$ be connected and vertically generic. Then the combinatorial isomorphism type of $\calP_\Gamma^\lambda$ is independent of $\lambda \in U$. Moreover, the equations defining $\calP_\Gamma^\lambda$ are all of the form $\sum \epsilon_i \Rea(s_i) \ge \Rea(2 \pi \mathrm{i}\lambda_I)$ and $\sum \epsilon_i\Ima(s_i) \ge 0$, where the variables $s_i$ denote periods of a certain finite set of distinguished saddle connections and $\epsilon_i \in \{-1,0,1\}$.}

\begin{thm}[Global homeomorphism type is invariant along a flat]\label{prop:pUproduct}
 With $U \subset B \subset \calR_p$ as above, there is a homeomorphism
 \[
 \calP^U \cong U \times \calP^\lambda
 \]
 for any $\lambda \in U$. 
\end{thm}

\begin{thm}[Period map is a homeomorphism]\label{thm:periodmaphomeo}
    Let $B \subset \calR_p$ be an open flat, and let $U \subset B$ be a connected vertically generic subset. Then there is a period map
    \[
    \per^U \colon \res^{-1}(U) \cong \calP^U
    \]
    which is a homeomorphism.
\end{thm}

\begin{thm}[$\res$ is a stratified fibration]\label{thm:stratfib}
   The restriction of $\res$ to any open flat $B \subset \calR_p$ realizes $\res^{-1}(B)$ as a locally trivial fibration over $B$. That is, $\res \colon \mathcal{H}(\mu) \to \calR_p$ is a stratified fibration relative to the resonance arrangement $\mathcal{A}_p \subset \calR_p$.
\end{thm}

Theorem \ref{thm:stratfib} is an essentially immediate corollary of Theorems \ref{prop:pUproduct} and \ref{thm:periodmaphomeo}; these results are established in Sections \ref{SS:periodmap} and \ref{SS:perioddomain}, respectively. Here we give the proof of Theorem \ref{thm:stratfib}, contingent on these results.

\begin{proof}[Proof of Theorem \ref{thm:stratfib}, given Theorems \ref{prop:pUproduct} and \ref{thm:periodmaphomeo}]
    Let $\lambda \in B$ be given. There are finitely many partial sums $\lambda_I \in \C$, and hence finitely many directions $\theta$ for which $e^{i \theta} \lambda$ fails to be vertically generic. $\res$ is equivariant with respect to the actions of $\C^*$ on $\mathcal{H}(\mu)$ and $\calR_p$, and so after adjusting by some $e^{i\theta}$, we can assume that $\lambda$ is vertically generic. As vertical genericity is an open condition, there exists a connected vertically generic polyhedral region $U \subset B$ containing $\lambda$. Combining Theorem \ref{thm:periodmaphomeo} and Theorem \ref{prop:pUproduct}, there is a homeomorphism $\res^{-1}(U) \cong U \times \calP^{\lambda}$, exhibiting the local product structure of $\res$. 
\end{proof}

\para{Organization of sections} The proofs of Theorems~\ref{lem:Pproduct},~\ref{prop:pUproduct}, and~\ref{thm:periodmaphomeo} are split between Sections~\ref{sec:Polyhedral} and \ref{sec:Period}. Section~\ref{sec:Polyhedral} is devoted to constructing and understanding the structure of the polyhedra $\calP^\lambda_\Gamma$, culminating in the proof of Theorem~\ref{lem:Pproduct}; Section \ref{sec:Period} then turns to more global matters. We begin in Section \ref{SS:critgraph} with a discussion of {\em critical graphs}, which encode the various possible topological types of vertical foliations of differentials in the stratum. In Section~\ref{SS:perioddomain}, we construct the ``local'' period domain $\calP^U_\Gamma$ associated to differentials with a fixed critical graph. In Section \ref{SS:isotype}, we show that the combinatorial structure of $\calP_\Gamma^U$ depends only on the signs of the partial sums $\Rea(\lambda_I)$ - this is the key ingredient needed in the proof of Theorem~\ref{prop:pUproduct}.

\subsection{Critical graphs}\label{SS:critgraph}
We begin our analysis of $\res$ by introducing the key combinatorial/topological bookkeeping tool: {\em critical graphs}. These encode the structure of the vertical foliation of a meromorphic differential, and are used to construct local coordinates for the period domains.

The work of the previous section shows that isoresidual fibers admit a metric completion given by the WYSIWYG compactification. We will likewise see such a structure naturally emerge from our perspective. Note that in genus zero, any curve underlying a differential in the WYSIWYG boundary is a ``tree of $\mathbb{CP}^1$'s'' (the $\mathbb{CP}^1$'s being the top-level components of the multi-scale differential, joint by the contraction images of the lower-level components). It will be convenient to give this class of curves (and the associated differentials) a name.

\begin{defn}[Arboreal curve, arboreal differential]
\label{def:arboreal}
    An {\em arboreal curve} $X = \bigcup X_i$ consists of components $X_i \cong \mathbb{CP}^{1}$ joined together at ordinary $m$-pointed singularities, such that $X$ is simply connected. An {\em arboreal differential} is a family of differentials $\omega_i$ on $X_i$ that is regular (has no pole) at each $p \in \operatorname{Sing}(X) \cap X_i$. By convention, we call every point $p \in \operatorname{Sing}(X)$ a {\em zero} of $\omega$, even when some $\omega_i$ on $X_i$ only has an ordinary marked point at $p$. 
\end{defn}

To define our desired period map on $\mathcal{H}(\mu)$, we first organize the differentials in this stratum by the topological type of their vertical foliation. This is described by an associated {\em critical graph}. For this purpose, we recall the standard terminology that a leaf of the vertical foliation is said to be {\em critical} if at least one end terminates at a zero or ordinary marked point, and is {\em regular} otherwise.

\begin{defn}[Critical graph]
    Let $X = \bigcup X_i$ be an arboreal curve, and let $\omega$ be an arboreal differential on $X$ with zeros at points $z_1, \dots, z_n \in X$ and poles at $w_1, \dots, w_p \in X$. Define the {\em critical graph} for $\omega$, written $\Gamma_\omega$, as the embedded graph in $X$ with vertex set $\{z_1, \dots, z_n, w_1, \dots, w_p\} \subset X$ consisting of the zeros and poles of $\omega$, and with edges given as the critical leaves of the vertical foliation associated to $\omega$. 
\end{defn}

{\em A priori}, it is not clear that the critical graph is well-defined: does every critical leaf connect a zero to a zero or pole? The next two results (Lemmas \ref{lemma:vnrcrit} and  \ref{lemma:vccs}) together show that when the configuration of residues $\lambda$ is vertically generic, this is true, and that moreover the critical graph determines the $1$-skeleton of a cell structure on $X$.

\begin{lem}\label{lemma:vnrcrit}
   Let $\omega$ be an arboreal differential, and suppose that $\omega$ is vertically generic. Then every critical leaf of the vertical foliation terminates at a zero or pole of $\omega$. 
\end{lem}
\begin{proof}
    The statement holds for an arboreal differential if and only if it holds for every meromorphic differential $\omega_i$ on each $X_i$, so we proceed under the assumption that $X = \mathbb{CP}^{1}$ and $\omega$ is an ordinary meromorphic differential.
    
    Following \cite[Proposition 5.5]{tahar}, we cut $\omega$ along all vertical saddle connections
    and examine the resulting decomposition into subsurfaces, known as {\em invariant components}. We first remark that a vertical saddle connection is precisely a critical leaf that terminates at a zero at both ends. According to \cite[Proposition 5.5]{tahar}, there are four possible types of invariant components: 
    \begin{enumerate}[label=(\roman*)]
        \item finite-volume cylinders foliated by vertical leaves,
        \item ``minimal components'' of finite volume, where the dynamics of the vertical foliation is induced by some nontrivial interval exchange map,
        \item infinite-volume cylinders bounding a simple pole (necessarily of purely real residue),
        \item ``free components'' where the generic leaves connect one pole to another (possibly itself), and the non-generic leaves are critical leaves connecting a zero to a pole.
    \end{enumerate}
    We claim that if $\omega$ is vertically generic, only invariant components of type (iv) can occur, from which it will follow that every critical leaf terminates at a pole or a zero. 

    We note first that regions of type (ii) cannot exist on $\mathbb{CP}^{1}$ - this is a consequence of the Poincar\'e--Bendixson theorem as applied to the vertical flow.
    To rule out components of types (i) and (iii), we employ the following principle: on $\mathbb{CP}^{1}$, the period of every oriented simple closed curve $c$ is a sum $2 \pi \mathrm{i}\sum_{i \in I} \lambda_i$ of the residues of the poles enclosed by $c$ and lying to the left. This immediately shows that regions of type (i) and (iii) cannot exist, since in either case, the core of the cylinder would have a period that is purely imaginary and nonzero. Since $\omega$ is vertically generic, if $\Rea(2 \pi \mathrm{i}\lambda_I) = 0$, then moreover $\lambda_I = 0$, and so no such curve can exist.
\end{proof}

\begin{lem}\label{lemma:vccs}
    Let $X$ be an arboreal curve, and let $\omega$ be a vertically generic arboreal differential on $X$. Then the critical graph $\Gamma_\omega \subset X$ forms the $1$-skeleton of a CW structure on $X$. There are two types of $2$-cells: 
   \begin{enumerate}[label=(\roman*)]
       \item {\bf Strips}, bounded on both the left and right by critical leaves,
       \item {\bf Half-spaces}, bounded on exactly one side by critical leaves.
   \end{enumerate}
\end{lem}
\begin{proof}
    We must show that each component of $X \setminus \Gamma_\omega$ is a topological $2$-disk. As in Lemma~\ref{lemma:vnrcrit}, it suffices to analyze the case of $X = \mathbb{CP}^{1}$. Following the analysis in the proof of Lemma~\ref{lemma:vnrcrit}, cutting along all vertical saddle connections decomposes $\mathbb{CP}^{1}$ into one or more invariant components of type (iv).
    Let $p$ be a point in the interior of $\mathbb{CP}^{1} \setminus \Gamma_\omega$. Consider a horizontal segment $h$ through $p$, and extend $h$ as far as possible to the left and right, until it meets a critical leaf (or, exceptionally, a zero). If $h$ extends infinitely far in both directions, then $\omega = dz$ in the stratum $\mathcal{H}(-2)$, but here $n+p=1 <3$, so this case is excluded by hypothesis. We claim that the region of $\mathbb{CP}^{1} \setminus \Gamma_\omega$ containing $p$ is homeomorphic to $h \times \R$, for $h$ such a maximal horizontal transversal. This is clear: every regular leaf in a region of type (iv) is homeomorphic to $\R$ and connects one pole to another, crossing $h$ exactly once. 
\end{proof}

\begin{defn}[Topological type]\label{def:toptype}
    Let $X = \mathbb{CP}^{1}$, and let $\omega, \omega' \in \mathcal{H}(\mu)$ be meromorphic differentials with zeros $(z_1, \dots, z_n)$ and $(z_1', \dots, z_n')$ and poles $(w_1, \dots, w_p)$ and $(w_1', \dots, w_p')$, respectively. Such $\omega, \omega'$ are said to have the same {\em topological type} if their critical graphs are isomorphic as embedded graphs in an oriented surface with ordered vertex set. In detail - they have the same topological type if there is an orientation-preserving homeomorphism $f \colon \mathbb{CP}^{1} \to \mathbb{CP}^{1}$ such that $f(z_i) = z_i'$ and $f(w_j) = w_j'$ for $1 \le i \le n$ and $1 \le j \le p$, and such that $f$ restricts to an isomorphism $f \colon \Gamma_\omega \to \Gamma_{\omega'}$.
    General arboreal differentials $\omega, \omega'$ on an arboreal curves $X$ are said to have the same {\em topological type} if the restriction $\omega_i, \omega'_i$ to each component $X_i \subset X$ has the same topological type in the above sense.
    If $\Gamma$ is the common isomorphism class of $\Gamma_\omega \cong \Gamma_{\omega'}$, we say that $\omega, \omega'$ have {\em topological type $\Gamma$}.  
\end{defn}

\para{Distinguished saddles} The decomposition of $X$ given by the vertical foliation distinguishes a certain finite subset of the saddles of $\omega$. These will be used to furnish distinguished period coordinates on the set of differentials of given topological type.

\begin{defn}[Distinguished saddles]\label{def:distsaddle}
    Let $\omega$ be an arboreal differential. A saddle connection $\sigma$ is {\em distinguished} (tacitly, with respect to the vertical foliation, or the underlying topological type) if its interior does not cross any critical leaf of the vertical foliation. There are two kinds of distinguished saddles: {\em transverse saddles}, whose interior points transversely intersect only regular leaves of the vertical foliation, and {\em vertical saddles}, i.e. critical leaves of the vertical foliation that connect two zeros.  
\end{defn}

There is a second graph associated to an arboreal differential - the {\em distinguished saddle graph} - that will be of paramount importance in the following section, where it will be used to describe systems of coordinates on the period domains. Note that if $\omega, \omega'$ have the same topological type, then their distinguished saddles are in natural bijection. We can therefore speak of the distinguished saddles associated to a topological type, not merely to a specific differential. Similarly, in a given critical graph, the vertices corresponding to zeros are well-defined at the level of isomorphism classes (as are the vertices corresponding to poles); we write $Z(\Gamma)$ for this subset of the vertices.

\begin{defn}[Distinguished saddle graph]
    Let $\Gamma$ be a topological type of arboreal differentials. The {\em distinguished saddle graph} $\Sigma_\Gamma$ is the graph with vertices given by the zeros of $\Gamma$, and with edges given by all distinguished saddles. 
\end{defn}

As an initial indication of the importance of the distinguished saddle graph, we establish the lemma below, which will be used in the following subsection. As a general piece of graph-theoretic terminology, let $\Gamma$ be a connected graph, and let $\calT$ be a spanning tree. Any edge $e$ not in $\calT$ then {\em determines a loop} in $\Gamma$ by connecting the endpoints of $e$ via the unique path in $\calT$.

\begin{lem}\label{lem:saddlebasis}
    Let $X = \bigcup X_i$ be an arboreal curve, and let $\Gamma$ be a critical graph on $X$ with set of zero vertices $Z(\Gamma)$. Then the distinguished saddles for $\Gamma$ span $H_1(X,Z(\Gamma); \Z)$, and a basis for $H_1(X, Z(\Gamma);\Z)$ is given by the set of saddles $\upsilon \in \calT$ for $\calT$ any spanning tree of the distinguished saddle graph $\Sigma_{\Gamma}$.
\end{lem}
\begin{proof}
    Let $\omega$ be any differential of topological type $\Gamma$. Let $\gamma \colon I \to X$ define a relative $1$-cycle. Then $\gamma$ can be homotoped to a union of saddle connections for $\omega$, each of which is either vertical, transverse, or else intersects some critical leaf on its interior. Any saddle of the latter type can be further homotoped into a union of transverse saddles, thereby proving the first claim. 

    To prove the second claim, let $\upsilon \in \Sigma_{\Gamma} \setminus \calT$ be given. Then $\upsilon$ determines a loop in $\Sigma_{\Gamma}$ (necessarily bounding some disk in some $X_i \subset X$) whose remaining edges are drawn from $\calT$, thereby expressing $\upsilon$ as a sum of cycles in $\calT$. The cycles in $\calT$ are linearly independent in $H_1(X, Z(\Gamma);\Z)$, as can be seen e.g. by constructing a ``dual basis'' in $H_1(X\setminus Z(\Gamma);\Z)$ of simple closed curves crossing exactly one edge of $\calT$ (since $\calT \subset X$ is an embedded tree, such a dual curve exists for every edge of $\calT$). Such cycles are dual to the elements of $\calT$ under the intersection pairing $H_1(X\setminus Z(\Gamma);\Z) \times H_1(X;Z(\Gamma);\Z) \to \Z$.
\end{proof}

\begin{ex}[Running example, I]\label{ex:run1}
    Here we introduce a running example designed to illustrate the ideas of the section. Let $\mu = (2, 1, 1, -2, -2, -2)$, and let $\lambda = (\lambda_1, \lambda_2, \lambda_3)$ be the configuration of residues shown in Figure \ref{fig:residues}. Figure \ref{fig:runex1} then depicts two differentials $\omega, \omega' \in \res^{-1}(\lambda) \subset \mathcal{H}(\mu)$, as well as the associated critical graphs and sets of distinguished saddles.
    \begin{figure}
\centering
		\labellist
        \small
        \pinlabel $2\pi~i\lambda_1$ at 10 18
        \pinlabel $2\pi~i\lambda_2$ at 45 18
        \pinlabel $2\pi~i\lambda_3$ at 25 0
		\endlabellist
\includegraphics[scale=1.5]{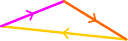}
\caption{The three residues $\lambda = (\lambda_1, \lambda_2, \lambda_3)$ (depicted as periods) used throughout the running example. The value of $\lambda$ is consistent across figures up to a positive real rescaling.}
\label{fig:residues}
\end{figure}

\begin{figure}
\centering
		\labellist
                \pinlabel $\omega$ [tl] at 435.08 504.52
\pinlabel $\Gamma$ [l] at 333.15 360.49
\pinlabel $\omega'$ [tl] at 442.16 173.06
\pinlabel $\Gamma'$ [l] at 339.81 5.35
        \small
\pinlabel $\sigma_1$ [b] at 51.02 566.88
\pinlabel $\sigma_2$ [bl] at 83.61 562.62
\pinlabel $\sigma_3$ [bl] at 214.00 566.88
\pinlabel $\sigma_4$ [b] at 252.26 559.79
\pinlabel $\sigma_1$ [b] at 58.10 240.33
\pinlabel $\tau_2$ [l] at 106.01 231.16
\pinlabel $\sigma_3$ [b] at 199.82 239.66
\pinlabel $\sigma_4$ [br] at 248.29 230.13
\pinlabel $\tau_2$ [r] at 392.56 225.49
\pinlabel $\tau_2$ [br] at 128.68 107.71
\pinlabel to~$w_1$~(res.~$\lambda_1$) at 75 640
\pinlabel to~$w_2$~(res.~$\lambda_2$) at 240 640
\pinlabel to~$w_3$~(res.~$\lambda_3$) at 400 640
\pinlabel to~$w_2$ at 51 490
\pinlabel to~$w_3$ at 90 490
\pinlabel to~$w_1$ at 215 490
\pinlabel to~$w_3$ at 255 490
\tiny
\pinlabel $w_3$ at 160 460
\pinlabel $w_2$ at 267 388
\pinlabel $w_1$ at 375 400
\pinlabel $\sigma_1$ at 295 370
\pinlabel $\sigma_2$ at 125 460
\pinlabel $\sigma_3$ at 275 413
\pinlabel $\sigma_4$ at 213 423
		\endlabellist
\includegraphics[width=0.88\textwidth]{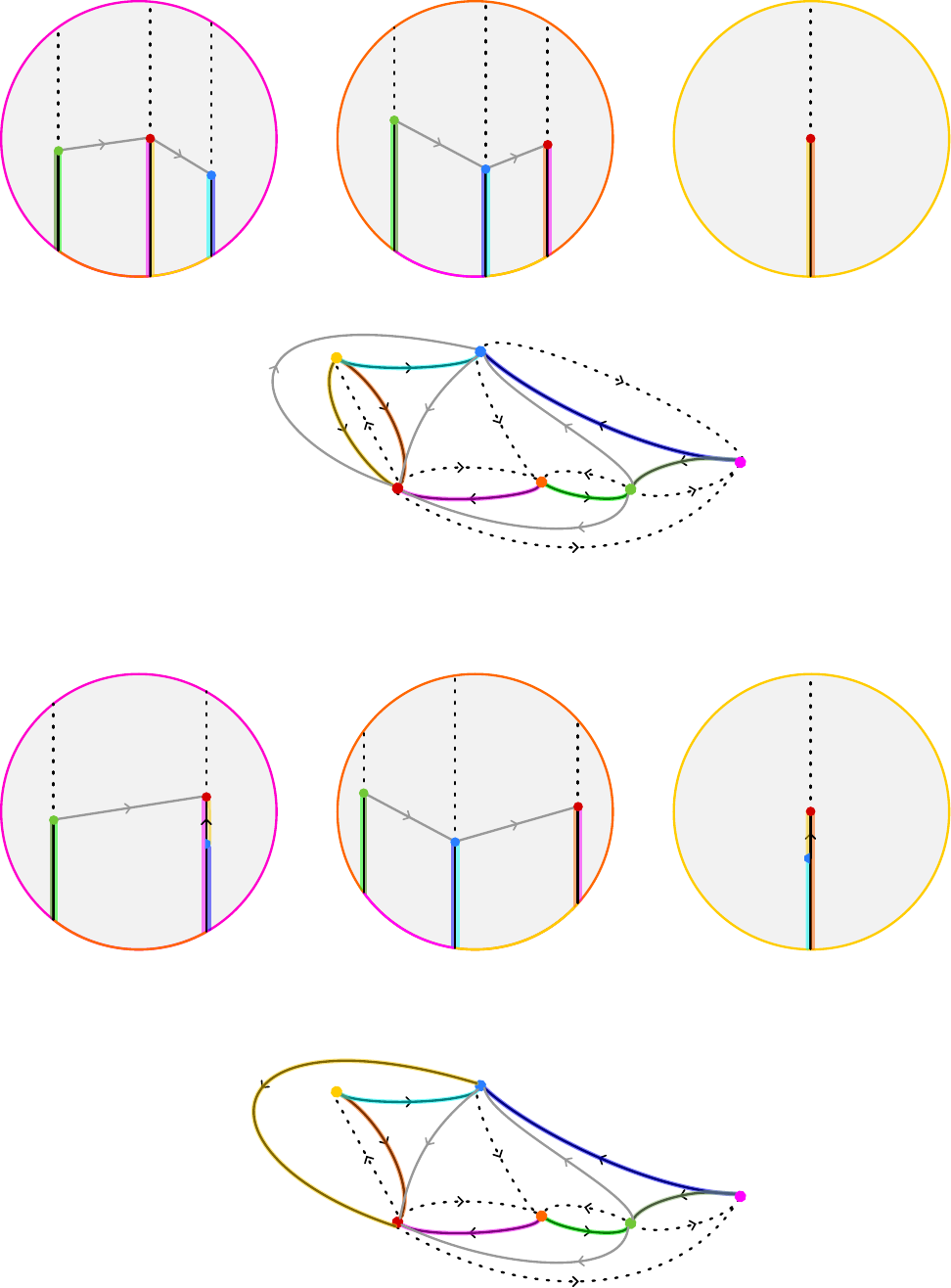}
\caption{This figure is best viewed in full color. From top to bottom: (1) The differential $\omega$. Each of the three shaded circles and their contents should be understood as complements of neighborhoods of each of the three poles (which portion of each circle limits onto which pole is indicated by color and by label). The three zeros are indicated by the colored dots in the middle of the circles. The slits running down from each zero have gluing instructions indicated by the colors on either side. Dotted lines indicate critical leaves for which the left and right half-spaces have been identified in the figure. There are four transverse saddles and no vertical saddles. (2) The critical graph $\Gamma = \Gamma_\omega$, enhanced with the transverse saddles $\sigma_i$. (3) The differential $\omega'$, which has one vertical saddle and three transverse saddles. (4) The critical graph $\Gamma' = \Gamma_{\omega'}$ and the associated transverse saddles.}
\label{fig:runex1}
\end{figure}
\end{ex}

\subsection{Period domains: construction (fixed topological type)}\label{SS:perioddomain}
The work of the previous section leads to natural polyhedra that parametrize the periods of the set of all differentials of a given topological type $\Gamma$. Briefly stated, the idea is to equip all of the distinguished saddles with a canonical orientation - transverse saddles run left-to-right (so that $\Rea(s_i) > 0$, where $s_i \in \C$ is the period of such a saddle $\sigma_i$), and vertical saddles run bottom-to-top (so that $t_j \in i\R_+$, where $t_j$ is the period of a vertical saddle $\tau_j$). Then there is a period mapping for differentials of type $\Gamma$ valued in the product of right half-spaces and positive imaginary axes, a {\em positive cone} in the ambient vector space. But the period map satisfies additional constraints coming from the homology relations among distinguished saddles, and so is valued in a certain {\em cycle subspace}. Intersecting this affine subspace with the positive cone leads to a polyhedral period domain.

Before constructing these polyhedra, we first lay out our terminology and conventions concerning polyhedra in the form most suitable to our situation at hand (where the polyhedron naturally lies in some subspace of an ambient space, and the equations defining both the subspace and half-spaces are not assumed to be independent). 

\begin{defn}[Polyhedron, isomorphism of polyhedra]\label{def:polyhedron}
    Fix an ambient Euclidean space $\R^d$. Let $I$ be a finite index set. For $i \in I$, let $\alpha_i \in (\R^d)^*$ be a functional, and let $a_i \in \R$; define the {\em half-space
    \[
    H_i^+ = \{x \in \R^d \mid \alpha_i(x) \ge a_i\} \subset \R^d
    \]} 
    as well as the {\em hyperplane} $H_i$ on which equality $\alpha_i(x) = a_i$ holds. For some subset $J \subset I$, define the {\em face space} 
    \[
    H_J = \bigcap_{j \in J} H_j.
    \]

    Now let $V \subset \R^d$ be an affine subspace. A {\em polyhedron} $P \subset V$ is defined as the intersection
    \[
    P = \bigcap_{i \in I} H_i^+ \cap V
    \]
    of some set of half-spaces in $\R^d$ with the subspace $V$. A {\em face} of $P$ is the intersection $H_J \cap P$ of some face space with $P$; a {\em facet} is a face of codimension one, defined by a single hyperplane $H_i$.
    
    The {\em face lattice} $\Lambda_{P}$ is the lattice of faces under inclusion. Polyhedra $P \subset V$ and $Q \subset W$ are {\em combinatorially isomorphic} if there is an isomorphism $\Lambda_P \cong \Lambda_Q$. In such a situation, we write $P \approx Q$.
    
    As a general remark, we will be somewhat lax with this terminology, by allowing the inequalities defining certain faces to be strict, still referring to the resulting object as a polyhedron.
\end{defn}

\para{Spherical polytopes} The polyhedra giving our period domains are unbounded. For technical purposes (specifically, to discuss the isomorphism type in Section \ref{SS:isotype} and triangulations in Section \ref{SS:perioddomainglobal}) it will be useful to compactify these as {\em spherical polytopes}, defined as the intersection of half-spaces in positive real projective spaces. Here, an unbounded polyhedron can be compactified by {\em homogenization}. For the reader's convenience, we recall the general setup (cf. \cite[Section 1.5]{ziegler}). 

\begin{defn}[Spherical polytope, homogenization, spherical compactification]\label{def:sphericalpoly}
    Let $\R^d$ be a real vector space, and let $\PP^+\R^d = (\R^d \setminus \{0\})/\R_+ \cong S^{d-1}$ be the positive projectivization. A {\em half-space} $H_\alpha^+ \subset \PP^+\R^d$ is determined by the equation $\alpha(x) \ge 0$, where $\alpha \in (\R^d)^*$ is any nonzero functional; as before we let $H_\alpha$ denote the projective hyperplane on which equality holds. As above, we will also notate these simply by $H_i^{(+)}$ when the defining functional $\alpha_i$ is unimportant. A {\em spherical polytope} is the intersection 
    \[
    P = \bigcap_{i = 1}^k H_{\alpha_i}^+ \cap \PP^+V \subset \PP^+\R^d
    \]
    of finitely many half-spaces with a projective subspace $\PP^+V \subset \PP^+\R^d$. The face lattice of a spherical polytope is defined as before, as is the notion of isomorphism of spherical polytopes.

    Now suppose $\alpha_i(x) \ge a_i$ determines a half-space $H_i^+$ in $\R^d$ in the sense of Definition \ref{def:polyhedron}. The {\em homogenization} of $H_i^+$ is the half-space in $\R\times\R^d$ given by $\alpha_i(x) \ge a_i x_0$, where $x_0$ is a distinguished coordinate in the first factor. Define the homogenization of an affine subspace $V \le \R^d$ similarly.
    
   Likewise, let $P \subset \R^d$ be a polyhedron (as in Definition \ref{def:polyhedron}), defined by the intersection of half-spaces $H_i^+ = \{x \in \R^d \mid \alpha_i(x) \ge a_i\}$ with some affine subspace $V \subset \R^d$. The {\em spherical compactification} $\overline{P} \subset \PP^+(\R \times \R^d)$ is defined to be the spherical polytope determined by the homogenizations of all defining half-spaces and affine subspaces, along with the additional half-space constraint $H_\infty^+$ given by $x_0 \ge 0$. The original $P$ can be recovered by intersecting $\overline{P}$ with the affine hyperplane $x_0 = 1$. 
\end{defn}

We can now define the period domain associated to the set of differentials of given topological type. Following the discussion above, we will work in the compactified setting.

\begin{defn}[Period space, cycle equation, cycle space]\label{def:cyclespace}
    Let $\Gamma$ be a critical graph, and let $\sigma_1, \dots, \sigma_h, \tau_1, \dots, \tau_v$ be the transverse and vertical saddles that are distinguished on $\Gamma$. Define the {\em period space} for $\Gamma$ to be the trivial real vector bundle over $\calR_p$
    \[
    E_\Gamma \coloneqq \calR_p \times V_\Gamma,
    \]
    where
    \[
    V_\Gamma \coloneqq \R \times \C^h\times i\R^v.
    \]
    We endow $V_\Gamma$ with coordinates as follows. Let  $x_0$ be a coordinate on the first factor, used for the purposes of compactification. Endow $\C^h$ with coordinates $s_1, \dots, s_h$ in bijection with the transverse saddles $\{\sigma_i\}$, and likewise endow $i\R^h$ with coordinates $t_1, \dots, t_v \in i\R$ in bijection with $\{\tau_j\}$ (as periods of vertical saddles, they are naturally valued on the imaginary axis).
    
    Now let $U \subset \calR_p$ be vertically generic. We next define the {\em cycle space} \[
    E_\Gamma^U \le E_\Gamma
    \]
    as follows: choose any spanning tree $\calT$ for $\Sigma_\Gamma$. Then each non-edge $e \in \Sigma_\Gamma \setminus \calT$ determines a loop $\gamma_e$ in $\Sigma_\Gamma$ which can be used to define a {\em cycle equation}, as follows. Topologically, when $\Gamma = \Gamma_\omega$ is realized by a differential $\omega$ with residues in $U$, such $\gamma_e$ is represented as a simple closed curve in $\mathbb{CP}^{1}$ enclosing some subset of the poles. After orienting $\gamma_e$, the residue theorem gives that the periods $s_i$ (resp. $t_j$) of the distinguished saddles $\sigma_i$ (resp. $\tau_j$) are related by a {\em cycle equation} of the form
    \begin{equation}\label{eqn:cycle}
        \sum \epsilon_{e,i}s_i + \sum \epsilon_{e,j}t_j = 2 \pi \mathrm{i} \lambda_{I_e},
    \end{equation}
    where each coefficient $\epsilon_{e,k}$ is valued in $\{-1,0,1\}$ and records whether a given saddle is involved in the loop (and if so, with which orientation), and $I_e \subset [p]$ is the set of poles lying to the left of $\gamma_e$. Note that this equation is in general inhomogeneous - to homogenize for the purposes of projectivization, insert a factor of the homogenizing variable $x_0$ on the right-hand side.

    For $\lambda \in U$, we define the subspace 
    \[
    V_\Gamma^\lambda \le V_\Gamma
    \]
    as the solutions to the cycle equations associated to the non-edges in some spanning tree $\calT$ for $\Sigma_\Gamma$; note that by Lemma~\ref{lem:saddlebasis}, $V_\Gamma^\lambda$ is independent of the choice of $\calT$. We then define the {\em cycle space} 
    \[
    E_\Gamma^U \le E_\Gamma
    \]
     as the subbundle of $E_\Gamma$ with fiber $V_\Gamma^\lambda \le V_\Gamma$ over $\lambda \in U$. 
\end{defn}

\begin{rmk}\label{rmk:cycleeqn}
    Note that $V_\Gamma$ is a {\em real} vector space, so that the real and imaginary parts of the cycle equations \eqref{eqn:cycle} technically need to be considered separately, noting that $x_0$ is purely real and each $t_j$ is purely imaginary. 
\end{rmk}

The period domain associated to a critical graph $\Gamma$ can now be defined as a family of polyhedra inside the cycle space $E_\Gamma^U$. We will also define three bordifications/compactifications - for a discussion of the use cases of each, see Remark \ref{rmk:usecases} below.

\begin{defn}[Period domain (fixed topological type)]\label{defn:perioddomainfixedtype}
    Continuing in the setting of Definition \ref{def:cyclespace}, let $C_\Gamma  \subset E_\Gamma$ be the family of cones defined by the strict inequalities $x_0 > 0, \Rea(s_i) > 0, \Ima(t_j) > 0$, and let $\overline{C_\Gamma}$ be the family of closed cones defined by the same conditions. Let also $\overline{C_\Gamma}^{pr}$ be the semi-open cone allowing $\Rea(s_i) = 0$ but requiring the remaining inequalities to be strict.
    
    For any subset $U \subset \calR_p$, the {\em period domain} $\calP_\Gamma^U$ is then defined as
    \[
    \calP_\Gamma^U \coloneqq C_\Gamma \cap E_\Gamma^U,
    \]
    while the {\em compactified period domain} is defined as
    \[
    \overline{\calP_\Gamma^U}^+ \coloneqq \overline{C_\Gamma} \cap E_\Gamma^U.
    \]
    The {\em proper bordification} is given by
    \[
    \overline{\calP_\Gamma^U}^{pr} \coloneqq \overline{C_\Gamma}^{pr} \cap E_\Gamma^U.
    \]
    Finally, define the {\em metric completion} of the period domain as the closed polyhedron in the affine setting, via
    \[
    \overline{\calP_\Gamma^U} \coloneqq \overline{\calP_\Gamma^U}^+ \cap \{x_0 = 1\}.
    \]
\end{defn}

\begin{rmk}[Where is each variant used?]\label{rmk:usecases}
    Let us comment briefly on the purpose of each of the constructions above. The open polyhedron $\calP_\Gamma^U$ is designed to capture the periods of every differential of topological type $\Gamma$ - see Lemma~\ref{lem:periodhomeolocal}. As mentioned above, the compactified period domain $\overline{\calP_\Gamma^U}^+$ is used in settings where having access to all possible faces is useful - this occurs in Theorem \ref{lem:Pproduct} and Theorem \ref{prop:pUproduct}. The proper bordification allows faces given by the inequalities $\Rea(s_i)$ to degenerate, but not the faces given by $\Ima(t_j)$. As seen in Section \ref{SS:degenerating}, degenerations of the first type correspond to degenerating to a differential of different topological type, still contained in the original stratum, while degenerations of the second type necessarily leave the stratum. Thus the proper bordification is relevant for the construction of the actual period domain (and not some compactification). Finally, the metric completion $\overline{\calP_\Gamma^U}$ allows for both types of degeneration, but excludes the faces at infinity present in the compactification $\overline{\calP_\Gamma^U}^+$. This is the largest subset of $\overline{\calP_\Gamma^U}^+$ on which we know how to give a ``modular interpretation'' of the limiting points, as arboreal differentials.
\end{rmk}

\begin{ex}[Running example, II]\label{ex:run2}
    Returning to the running example of Example \ref{ex:run1}, here we describe the cycle spaces for the topological types $\Gamma$ and $\Gamma'$. For $\Gamma$, we choose the spanning tree $\calT = \{\sigma_2, \sigma_3\}$. There are two distinguished saddles not in $\calT$, namely $\sigma_1$ and $\sigma_4$. By examining the critical graph $\Gamma$ shown in Figure \ref{fig:runex1}, the (affine) cycle equation for $\sigma_1$ is found to be $s_1 +s_2 - s_3 = 2 \pi \mathrm{i}\lambda_1$, and the one for $\sigma_4$ is $s_2 + s_4 = 2 \pi \mathrm{i}(\lambda_1 + \lambda_2)$.

    For $\Gamma'$, we choose the spanning tree $\calT' = \{\tau_2, \sigma_3\}$. Then the distinguished saddles not in $\calT'$ are $\sigma_1$ and $\sigma_4$. The cycle equations here can be worked out from the principle that $\omega'$ is a deformation of $\omega$, under which $\sigma_2$ on $\omega$ limits to $\tau_2$ on $\omega'$ with the opposite orientation (the remaining saddles on $\omega$ limit onto their counterparts on $\omega'$). Thus the cycle equations on $\omega'$ are simply given by substituting $-t_2$ for $s_2$.
\end{ex}

\subsection{Period domains: polyhedral structure}\label{SS:isotype}
The goal of this section is to prove the first of our key results, Theorem \ref{lem:Pproduct}. This shows that the combinatorial type of $\calP^\lambda_\Gamma$ depends only on the signs of the partial sums $\Rea(2 \pi \mathrm{i}\lambda_I)$ (recall our standing convention that $U \subset B \subset \calR_p$ is a connected vertically generic subset of some open flat; in particular, the signs of any such partial sum are constant along $U$). 

This will require a more careful investigation of the relationship between the graph of distinguished saddles and coordinates on the cycle spaces, as mediated by well-chosen spanning trees.

\begin{defn}[Adapted spanning tree]\label{def:Cadapted}
    Let $\Gamma$ be a critical graph, and let $C = \{\sigma_1, \dots, \sigma_a, \tau_1, \dots, \tau_b\} $ be some set of distinguished saddles ($a$ and/or $b$ may be zero). A spanning tree $\calT$ for the graph of distinguished saddles $\Sigma_\Gamma$ is {\em adapted to $C$} if it satisfies the following properties:
    \begin{enumerate}[label=(\Roman*)]
        \item Any $\tau_j \in C$ not in $\calT$ determines a loop in $\Sigma_\Gamma$ consisting entirely of vertical saddles $\tau_\ell \in C$,
        \item Any $\tau_j$ not in $\calT$ determines a loop in $\Sigma_\Gamma$ consisting entirely of vertical saddles $\tau_\ell$,
        \item Any $\sigma_i \in C$ not in $\calT$ determines a loop in $\Sigma_\Gamma$ consisting entirely of vertical saddles and transverse saddles $\sigma_k \in C$.
    \end{enumerate}
\end{defn}

\begin{lem}\label{lem:adaptedexists}
    Let $\Gamma$ be a critical graph, and let $C = \{\sigma_1, \dots, \sigma_a, \tau_1, \dots, \tau_b\} $ be any set of distinguished saddles. Then a spanning tree $\calT$ for $\Sigma_\Gamma$ adapted to $C$ exists.
\end{lem}
\begin{proof}
    Running from $\tau_1$ to $\tau_b$, adjoin $\tau_i$ to $\calT$ if and only if the result does not form a loop. If the resulting graph does not yet form a spanning tree, then proceed through the remaining vertical saddles in any order, again adjoining $\tau_i$
    if and only if the result does not for a loop. If a spanning tree has not yet been constructed at this stage, proceed through
    $\sigma_1, \dots, \sigma_a \in C$, adjoining $\sigma_i$ if and only if the result does not form a loop. Finally if the resulting graph does not yet form a spanning tree, adjoin additional edges arbitrarily to create a spanning tree $\calT$. It follows immediately from the construction that such $\calT$ is adapted to $C$. 
\end{proof}

We will apply these notions to the problem of parameterizing the cycle space $V_\Gamma^\lambda \le V_\Gamma$. It will be helpful to establish some terminology. We say a set of variables $\calB \subset \{x_0, s_i, t_j\}$ is a {\em set of coordinates} for $V_\Gamma^\lambda$ if all remaining such variables depend $\R$-linearly on these (treating $\Rea(s_i), \Ima(s_i)$ as separate variables in keeping with Remark \ref{rmk:cycleeqn}); any such variable of the latter type is said to be {\em dependent}.

\begin{lem}\label{lem:adaptedcoordinates} 
    Let $\Gamma$ be a critical graph, let $\lambda \in \calR_p$ be vertically generic, and let $C = \{\sigma_1, \dots, \sigma_a, \tau_1, \dots, \tau_b\} $ be any set of distinguished saddles. Then there is a set of coordinates $\calB_C$ for $V_\Gamma^\lambda$ with the properties listed below. We say such a set of coordinates is {\em adapted to $C$}.
     \begin{enumerate}[label=(\Roman*)]
        \item If $t_j \in C$ is dependent, there is an expression of the form
        \begin{equation} \label{eq:specialtcycle}
        t_j = \sum_{t_\ell \in C}\epsilon_{j\ell} t_\ell,
        \end{equation}
        expressing $t_j$ in terms of coordinates drawn from $C$ alone (here and below, any symbol of the form $\epsilon_{pq}$ is valued in the set $\{-1,0,1\}$),
        \item If $t_j$ is dependent, there is an expression of the form
        \begin{equation}\label{eq:tcycle}
        t_j = \sum \epsilon_{j\ell} t_\ell
        \end{equation}
        expressing $t_j$ in terms of vertical saddle coordinates alone,
        \item If $s_i \in C$, then either both $\Rea(s_i), \Ima(s_i)$ are dependent or neither is; if both are dependent, then there is an expression
            \begin{equation}\label{eq:scycle}
                s_i = 2 \pi \mathrm{i}\lambda_{I_i} x_0 + \sum_{s_k \in C} \epsilon_{ik}s_k + \sum \epsilon_{i\ell}t_\ell
            \end{equation}
        for some subset $I_i \subset [p]$ of the poles; note that the real part of this expression is independent of $t_\ell$.
    \item For any $s_i \not\in C$, either both $\Rea(s_i), \Ima(s_i)$ are dependent or neither is. If both are dependent, there is an expression as in \eqref{eq:scycle} with $s_k$ now drawn from $\calB_C$, not just $C$. 
    \end{enumerate}
\end{lem}
\begin{proof}
    The space $V_\Gamma^\lambda$ is defined by the cycle equations \eqref{eqn:cycle} (indexed by graph-theoretic cycles in the distinguished cycle graph $\Sigma_\Gamma$), where each cycle equation determines a pair of real equations as in Remark \ref{rmk:cycleeqn}.
    Let $\calT$ be a spanning tree for $\Sigma_\Gamma$ adapted to $C$, as in Lemma~\ref{lem:adaptedexists}. Define $\calB_C$ to consist of $x_0$ and those $s_i,t_j\in \calT$.

    Having chosen a spanning tree $\calT$, the cycle equations are in bijection with the edges not contained in $\calT$. For an edge $t_j \not \in \calT$, since $\calT$ is adapted to $C$, property (II) in Definition \ref{def:Cadapted} implies that the corresponding equation is of the form
    \[
    t_j = 2 \pi \mathrm{i}\lambda_{I_j}x_0 + \sum_{t_\ell \in \calT} \epsilon_{j\ell} t_\ell
    \]
    for some partial sum $\lambda_{I_j}$ of the residues. Since $\Rea(t_k) = 0$ for all $k$, it follows that $\Rea(2 \pi \mathrm{i}\lambda_{I_j}) = 0$, and hence by the vertical genericity assumption, $\lambda_{I_j} = 0$. If moreover $t_j \in C$, then by property (I) of an adapted spanning tree, the variables on the right-hand side are all drawn from $C$ as required.  Likewise, any $s_i \not \in \calT$ determines an equation of type \eqref{eq:scycle}, and by property (III), if $s_i \in C$, then all the $s_k$ terms on the right-hand side are drawn from $C$ as well.
    We have thus shown that every dependent variable can be expressed in terms of the chosen coordinates, giving the desired parameterization of $V_\Gamma^\lambda$. 
\end{proof}

In preparation for Theorem \ref{lem:Pproduct}, we next use this theory of adapted coordinates to make some preliminary observations about the face structure of $\overline{\calP_\Gamma^\lambda}^+$. Recall (e.g. \cite[Defn 1.11]{ziegler}) that the {\em lineality space} of a spherical polytope $P \subset \PP^+V$ is the largest projective subspace $\PP^+W \subset P$ contained in $P$. Define 
\[
W_\Gamma \le V_\Gamma
\]
as the span of the variables $\{\Ima(s_i)\}$, and 
\[
W_\Gamma^\lambda \coloneqq W_\Gamma \cap V_\Gamma^\lambda \le V^\lambda_\Gamma.
\]
Observe that the constraints imposed by the cone $\overline{C_\Gamma}$ are independent of $\Ima(s_i)$, so that $\PP^+W_\Gamma^\lambda \subset \overline{\calP_\Gamma^\lambda}^+$, i.e. the lineality space of $\overline{\calP_\Gamma^\lambda}^+$ is nonempty. In the lemma below, recall that $n$ denotes the number of zeros of $\Gamma$, equivalently the number of vertices of $\Sigma_\Gamma$.

\begin{lem}\label{lem:Plineality}
$\operatorname{codim}(\PP^+W_\Gamma^\lambda \subset \PP^+V_\Gamma^\lambda) = n$ for any $\lambda \in U$.
\end{lem}
\begin{proof}
    Let $\calB = \{x_0, s_k, t_\ell\}$ be a set of coordinates for $V_\Gamma^\lambda$ adapted to $C = \emptyset$ as in Lemma~\ref{lem:adaptedcoordinates}. Then the subset of coordinates $\Ima(s_k)$ restrict to give coordinates on $W_\Gamma^\lambda$. There are $n$ remaining coordinates: $x_0$, along with $n-1$ additional coordinates given by either $\Rea(s_k)$ or $t_j$, in bijection with the $n-1$ edges of the adapted spanning tree.
\end{proof}

In a spherical polytope $P \subset \PP^+V$ with nontrivial lineality space $\PP^+L$, the analogue of a vertex is played by a face space $\PP^+W$ containing $\PP^+L$ with codimension one. We call these {\em minimal faces}.
To prove Theorem \ref{lem:Pproduct}, we must characterize when a set of hyperplanes intersect in a minimal face. The possible hyperplanes (or their defining functionals) are drawn from the set 
\[
\calH_\Gamma\coloneqq \{x_0, \Rea(s_1), \dots, \Rea(s_h), \Ima(t_1), \dots, \Ima(t_v)\}.
\]
Functionals of the latter two types ($\Rea(s_i)$ or $\Ima(t_j)$) have an {\em associated edge} in the graph $\Sigma_\Gamma$, either $\sigma_i$ or $\tau_j$ as appropriate. For a subset $C \subset \calH_\Gamma$, we define a subgraph 
\[
\Sigma_{\Gamma, C}\subset \Sigma_\Gamma
\] 
as the subgraph of all such associated edges.

We are now ready to prove Theorem \ref{lem:Pproduct}; recall that this asserts that the combinatorial isomorphism type of $\calP^\lambda_\Gamma$ is independent of $\lambda \in U$, for $U \subset B$ a connected vertically generic subset of an open flat $B \subset \calR_p$. The assertions in Theorem \ref{lem:Pproduct} about the structure of the equations have already been established in the course of constructing $\calP^\lambda_\Gamma$.

\begin{proof}[Proof of Theorem \ref{lem:Pproduct}]
  $\calP_\Gamma^\lambda$ is given as the interior of the spherical compactification $\overline{\calP_\Gamma^\lambda}^+$; we will prove that the combinatorial isomorphism type of the latter is independent of $\lambda \in U$. 
  The faces of $\overline{\calP_\Gamma^\lambda}^+$ are given by intersecting some collection of half-spaces $C \subset \calH_\Gamma$ with the subspace $\PP^+V_\Gamma^\lambda$. To show that the combinatorial isomorphism type of $\overline{\calP_\Gamma^\lambda}^+$ is constant for $\lambda \in U$, we will show that whether such a face is empty is independent of $\lambda \in U$. Throughout, we will tacitly use the fact that since $\lambda \in U$ is connected and vertically generic, the sign of any expression $\Rea(2 \pi \mathrm{i} \lambda_I)$ is constant for $\lambda \in U$. 

    Any spherical polytope is determined by its set of minimal-dimensional faces. We will show that whether such an intersection lies in $\overline{\calP_\Gamma^\lambda}^+$ is independent of $\lambda$. Let $C \subset \calH_\Gamma$ index a set of functionals on $V_\Gamma$ (equivalently, a set of hyperplanes), and suppose that it determines a minimal face, i.e. that 
    \[
    \dim\left(\bigcap_{\alpha \in C} H_\alpha \cap V_\Gamma^\lambda \pmod {W_\Gamma^\lambda}\right)=1.
    \]
    For simplicity, denote the subspace $\bigcap_{\alpha \in C} H_\alpha \cap V_\Gamma^\lambda \pmod {W_\Gamma^\lambda}$ of $V_\Gamma^\lambda/W_\Gamma^\lambda$ by $F_C$.
    
    Let $\calB_C = \{x_0, s_k, t_\ell\}$ be a set of coordinates for $V_\Gamma^\lambda$ adapted to $C$ as in Lemma~\ref{lem:adaptedcoordinates}. Then $\dim(F_C) = 1$ if and only if the conditions $\alpha \equiv 0\ (\alpha \in C)$ together define a subspace of dimension one in the span $\R\calB_C$ of $\calB_C$. A condition $\alpha \equiv 0$ for $\alpha \in C$ imposes a constraint in $\R \calB_C$ in one of two ways: either {\em directly}, if $\alpha \in \calB_C$, or {\em indirectly}, if $\alpha \not \in \calB_C$. 

    As $\calB_C$ is adapted to $C$, equation \eqref{eq:specialtcycle} shows that for every $t_j \in C \setminus \calB_C$, both sides of the associated cycle equation vanish, so that no indirect constraint is imposed by such a functional. A functional $\Rea(s_i) \in C \setminus \calB_C$ has the real part of the associated cycle equation \eqref{eq:scycle} reduce to $0 = \pm \Rea(2 \pi \mathrm{i}\lambda_{I_i}) x_0$. This imposes an indirect constraint only if $x_0 \not \in C$ and $\Rea(2 \pi \mathrm{i}\lambda_{I_i}) \ne 0$. In this case, we can exchange $\Rea(s_i)$ for $x_0$ in $C$ without affecting $F_C$, removing all indirect constraints in the process.

    We conclude that when $F_C$ is a minimal face, $\calB_C \setminus C = \{\alpha\}$ is a singleton. Suppose first that $\alpha = x_0$. Then all equations \eqref{eq:tcycle} reduce to $t_j = 0$, and so the half-space condition $\Ima(t_j) \ge 0$ is necessarily satisfied. A cycle equation expressing some $\Rea(s_i)$ in terms of $\calB_C$ reduces to $\Rea(s_i) = \pm \Rea(2 \pi \mathrm{i}\lambda_{I_i})x_0$; the sign of $\Rea(s_i)$ is therefore constant for $\lambda \in U$ as claimed.

    Next suppose that $\alpha = \Ima(t_\ell)$. Here, all expressions for $\Rea(s_i)$ reduce to $\Rea(s_i) = 0$. For $t_j \not \in C$, \eqref{eq:tcycle} reduces to $t_j = \epsilon_{j\ell} t_\ell$, with $\epsilon_{j\ell} \in \{-1,0,1\}$ determined purely by the topological information encoded by $\Gamma$. We again conclude that whether the associated face satisfies all the half-space conditions is independent of $\lambda \in U$.

    Finally suppose that $\alpha = \Rea(s_k)$. As in the first case, all remaining $t_j = 0$, and similarly to the second case, all remaining $\Rea(s_i)$ admit expressions of the form $\Rea(s_i) = \epsilon_{ik} \Rea(s_k)$, again independent of $\lambda \in U$. 
\end{proof}

\begin{ex}[Running example, III]\label{ex:run3}
    Here we give a concrete description of the polytope $\calP_\Gamma^\lambda$. Continuing from Example \ref{ex:run2}, we see that $s_2, s_3$ are an adapted set of coordinates for $C = \emptyset$, and we solve the cycle equations obtained in Example \ref{ex:run2} for the remaining periods:
    \[
    s_1 = -s_2 + s_3 +2 \pi \mathrm{i}\lambda_1 \quad s_4 = -s_2 - 2 \pi \mathrm{i}\lambda_3.
    \]
    Thus, in $\C^2 = \C\langle s_2, s_3 \rangle$, the polyhedron $\calP_{\Gamma}^\lambda$ is given by the four half-space equations
    \begin{align*}
        -\Rea(s_2) + \Rea(s_3) + \Rea(2 \pi \mathrm{i}\lambda_1) > 0 &\quad & \Rea(s_2) > 0\\
        \Rea(s_3) > 0 &\quad & \Rea(s_2) + \Rea(2 \pi \mathrm{i}\lambda_3) < 0.
    \end{align*}
    The fact that these are independent of $\Ima(s_2), \Ima(s_3)$ is a reflection of the nonempty lineality space (Lemma \ref{lem:Plineality}).
    With reference to Figure \ref{fig:residues}, we see that $\Rea(2 \pi \mathrm{i}\lambda_1) > 0$, while $\Rea(2 \pi \mathrm{i}\lambda_3) < 0$. We can therefore depict $\calP_\Gamma^\lambda$ as in Figure \ref{fig:polytope}.
    Anticipating Section \ref{SS:degenerating}, the polytope $\calP_{\Gamma'}^\lambda$ appears in the boundary of $\calP_\Gamma^\lambda$, as the region given by $\Ima(t_2) = -\Ima(s_2) > 0$ and $\Rea(t_2) = -\Rea(s_2) = 0$, i.e. as the lower half of the facet given by $\Rea(s_2) = 0$. 

    \begin{figure}
\centering
		\labellist
        \small
        \pinlabel $\Rea(s_2)$ at 143 51
        \pinlabel $\Rea(s_3)$ at 9 146
        \tiny
        \pinlabel $\Rea(s_1)>0$ at 57 80
        \pinlabel $\Rea(s_2)>0$ at 23 100
        \pinlabel $\Rea(s_3)>0$ at 23 55
        \pinlabel $\Rea(s_4)>0$ at 80 130
		\endlabellist
\includegraphics[scale=1]{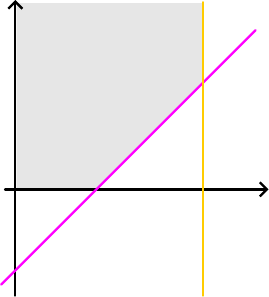}
\caption{The polytope $\calP_{\Gamma}^\lambda$ of our running example (after quotienting by the lineality space). In keeping with Theorem \ref{lem:Pproduct}, note that the isomorphism type of $\calP^\lambda_\Gamma$ is completely determined by the conditions $\Rea(\lambda_3) < 0 < \Rea(2 \pi \mathrm{i}\lambda_1) < \Rea(-2 \pi \mathrm{i}\lambda_3) = \Rea(2 \pi \mathrm{i}\lambda_1) + \Rea(2 \pi \mathrm{i}\lambda_2)$.}
\label{fig:polytope}
\end{figure}
\end{ex}

\section{The period map}\label{sec:Period}

We continue our analysis of the flat geometry of the isoresidual fiber in genus zero. In the previous section we constructed the polyhedron $\calP^\lambda_\Gamma$ associated to a critical graph $\Gamma$, culminating in the proof of Theorem \ref{lem:Pproduct} concerning the invariance of its combinatorial type as $\lambda$ varies in a connected non-resonant subset of an open flat. We turn our sights now to the remaining results Theorem \ref{prop:pUproduct} and Theorem \ref{thm:periodmaphomeo}. There are two tasks to accomplish: we must give the {\em global} construction of the period domain $\calP^\lambda$, by assembling the polyhedra $\calP^\lambda_\Gamma$, and we must construct the period map (at both the local and global levels).

The section is organized as follows. The local period map (for a fixed critical graph) is described in Section \ref{SS:localperiod}. In Section \ref{SS:degenerating}, we show how to identify points on the boundary of a polyhedron $\calP^\lambda_\Gamma$ with limit differentials, leading to our perspective on the modular interpretation of the metric completion. In Section \ref{SS:perioddomainglobal}, we use this to assemble the local period domains $\calP^\lambda_\Gamma$ into the global object $\calP^\lambda$, and prove Theorem \ref{prop:pUproduct}. Finally in Section \ref{SS:periodmap}, we define the global period map and prove it gives a homeomorphism $\calF_\lambda(\mu) \cong \calP^\lambda$ (Theorem \ref{thm:periodmaphomeo}).

\subsection{Period mapping (fixed topological type)}\label{SS:localperiod}

Here we give the construction of the period mapping for a set of differentials of fixed topological type, and prove that this is a homeomorphism.

\begin{defn}[Period map (fixed topological type)]
    Let $\Gamma$ be a critical graph, and let  $\res^{-1}(U)_{\Gamma}$ denote the subset of $\res^{-1}(U)$ consisting of differentials $\omega$ of topological type $\Gamma$. Choose orderings $\sigma_1, \dots, \sigma_h, \tau_1, \dots, \tau_v$ of the distinguished saddles in $\Gamma$.
    The {\em period map} over $\Gamma$,
    \[
    \per^U_{\Gamma} \colon \res^{-1}(U)_{\Gamma} \to \calP^U_{\Gamma} \subset E_\Gamma^U
    \]
    is then defined by 
    \[
    \per^U_\Gamma(\omega) = (\res(\omega), (1,s_1, \dots, s_h, t_1, \dots, t_v))
    \]
    where $s_i$ (resp. $t_j$) is the period of the distinguished saddle $\sigma_i$ (resp. $\tau_j$) (recall we work with the homogenized cycle equations and so set the first coordinate $x_0 = 1$). 
\end{defn}

\begin{lem}\label{lem:periodhomeolocal}
    $\per^U_\Gamma$ induces a homeomorphism $\res^{-1}(U)_\Gamma \cong \calP_\Gamma^U$.
\end{lem}
\begin{proof}
    The period of a given saddle varies continuously under deformation, showing that $\per^U_\Gamma$ is continuous. We will exhibit a continuous inverse. This is a direct construction. We must build a set of strips and half-spaces realizing the desired periods of the distinguished saddles. Once these are constructed, the gluing information encoded by $\Gamma$ gives the instructions for how to assemble them into an arboreal differential on $X$, and the construction will be complete.

    Let $S$ be a strip, and let $w \in \calP_\Gamma^U$ be any point. Suppose that there are $k_L$ vertical saddles appearing on the left side of $S$, and $k_R$ vertical saddles on the right (either of these quantities may be zero); enumerate these running from bottom to top as $\tau_1^L, \dots, \tau_{k_L}^L$ and $\tau_1^R, \dots, \tau_{k_R}^R$, respectively. 
    The geometries of the left and right sides of $S$ are then fixed by the associated periods $\{t_i^L\}$ and $\{t_i^R\}$ encoded in $w \in \calP_\Gamma^U$. On the left side, running from bottom to top, there is a critical leaf running from a pole up to the bottom-most zero appearing on the left side, followed by the $k_L$ saddles $\tau_1^L, \dots, \tau_{k_L}^L$, followed by a critical leaf running up from the topmost zero to the pole. The right side of $S$ admits an analogous description. 

    The transverse saddles contained in $S$ then admit the following description: there is one transverse saddle for each pair of a vertex on the left and right side, for a total of $(k_L+1)(k_R+1)$. Choosing any pair of vertices on the left and right sides, $S$ is then completely determined by setting the period of the chosen saddle to the corresponding coordinate of $w$. The cycle equations imply that the period of any other transverse saddle is determined by this value and the periods of the vertical saddles, so that the geometry of the constructed strip is independent of which transverse saddle is selected.

    The construction for half-spaces is similar. Running from bottom to top, the boundary of a half-space consists of a critical leaf running from the pole up to the bottom-most critical point, followed by a sequence of vertical saddles, and finally a critical leaf running from the top-most critical point up to the pole. The periods (i.e. vertical lengths) of the vertical saddles are prescribed by the period information encoded by $w \in \calP_\Gamma^U$. 
\end{proof}

\subsection{Extending to the boundary}\label{SS:degenerating}

Our next objective is to describe an extension of the inverse $(\per^U_\Gamma)^{-1}$ to the metric completion $\overline{\calP_\Gamma^U}$ (recall Definition \ref{defn:perioddomainfixedtype}). The reconstruction process described in the proof of Lemma \ref{lem:periodhomeolocal} just above suggests that points in the boundary can likewise be incarnated as arboreal differentials, but where the critical graph undergoes some type of degeneration that remains to be better understood. There are two types of facets of $\overline{\calP_\Gamma^U}$: those given by $\Rea(s_i) = 0$, and those given by $\Ima(t_j) = 0$. Approaching each of these faces has a clear geometric interpretation. The second type is easy to describe - approaching a facet defined by $\Ima(t_j) = 0$ causes the associated vertical saddle $\tau_j$ to shrink to a point, merging the zeros at either end and possibly causing the underlying curve to degenerate if multiple homologous saddles collapse at once (see the discussion of {\em linking} just below).\\

\para{Collapsing strips} Approaching a facet defined by $\Rea(s_i) = 0$ corresponds topologically to a degeneration we call ``collapsing a strip'' - we pause here to give a brief discussion. Under a deformation of differentials limiting onto such a facet, the strip containing the associated saddle $\sigma_i$ {\em collapses}, degenerating the transverse saddles within to new vertical saddles. Examples of this can be seen in the degenerations of $\omega$ onto $\omega'$ and $\omega''$ in Figures \ref{fig:runex1}, \ref{fig:runex2}, respectively - in each of these examples, the strip containing $\sigma_2$ collapses, and the transverse saddle $\sigma_2$ degenerates to the vertical saddle $\tau_2$. A key feature of this type of degeneration is that the limiting topological type depends on the value of $\Ima(s_i)$. The degenerations shown in Figures \ref{fig:runex1}, \ref{fig:runex2} are onto the same face ($\Rea(s_2) = 0$) of the same polytope $\calP^\lambda_\Gamma$; the difference is governed by the sign of $\Ima(s_2)$.\\

\para{Linking} There is one further subtlety to understand - in the presence of resonance, there can be multiple saddles constrained to share the same period. This can force several strips to collapse at once when approaching a facet given by $\Rea(s_i) = 0$, and can cause several vertical saddles to shrink simultaneously when approaching a facet given by $\Rea(t_j) = 0$. This can cause the underlying curve to degenerate; it is the source of our need to consider arboreal differentials.

\begin{defn}[Linked saddles, strips]
    Let $\Gamma$ be a critical graph. Let $\upsilon, \upsilon'$ be distinguished saddles (either transverse or vertical), and suppose that each of $\upsilon, \upsilon'$ has two distinct endpoints (i.e. neither is a loop in the distinguished saddle graph $\Sigma_\Gamma$). Let their periods be denoted $u,u'$, respectively, and view these as functions on $\res^{-1}(U)_\Gamma$. Such saddles are said to be {\em linked} if $u \equiv \pm u'$ holds identically on  $\res^{-1}(U)_\Gamma$.
    Strips $S, S'$ determined by $\Gamma$ are said to be {\em linked} if there are transverse saddles $\sigma \subset S$ and $\sigma' \subset S'$ that are linked in the above sense.
     
\end{defn}
\begin{rmk}
    We require the endpoints of $\upsilon$ to be distinct in the above definition in order to exclude a ``degenerate'' form of linking. If $\upsilon$ connects some zero to itself, then its period is an absolute period, given as the sum of residues of the poles lying to one side. If there happens to be a second such $\upsilon'$, and if the corresponding sums of residues happen to be equal throughout $U$ (e.g. if $U =\{\lambda\}$ is a singleton), then $\upsilon$ and $\upsilon'$ have equal periods throughout $\calP_\Gamma^U$. On the other hand, such behavior does not factor into our analysis below, which is concerned with when the periods of saddles degenerate (transition into being purely imaginary, or zero). {\em Absolute} periods do not degenerate on the sets $U$ we consider, as $U$ is a vertically generic subset of an open flat by hypothesis.
\end{rmk}

The next lemma characterizes when linkage occurs, and shows that it is an effect caused by resonance.

\begin{lem}\label{lem:whenequal}
    Let $\Gamma$ be a critical graph. Saddles $\upsilon, \upsilon'$ are linked if and only if they connect the same zeros, and moreover the sum of the residues enclosed by $\upsilon \cup \upsilon'$ on either side vanishes.
\end{lem}
\begin{proof}
    The `if' assertion is an immediate consequence of the residue theorem. For the reverse implication, suppose that the endpoints of $\upsilon, \upsilon'$ consist of at least three points (if the endpoints are the same but the sum of the enclosed residues is nonzero, the result is again clear from the residue theorem). Then $\upsilon$ and $\upsilon'$ determine distinct elements of $H_1(X,Z(\Gamma);\Z)$, and so can be completed to a basis. The periods of $\upsilon, \upsilon'$ can then be perturbed independently while remaining in $\calP_\Gamma^U$ (as can be seen by an appeal to the theory of adapted coordinates as in Lemma \ref{lem:adaptedcoordinates}).
\end{proof}

This example shows how shrinking linked vertical saddles can lead to a degeneration of the underlying curve, necessitating the introduction of arboreal differentials.

\begin{ex}\label{ex:degenerating}
    Let $\mu = (1,1, -2,-2)$, and let $\lambda = (0,0)$. Consider a differential $\omega$ obtained by taking two copies of $(\mathbb{CP}^{1}, dz)$ and joining them along a vertical slit consisting of two vertical saddles $\tau_1, \tau_2$. Then $\tau_1$ and $\tau_2$ are linked. Shrinking the periods $t_1 = t_2$ down to $t_1 = t_2 = 0$ yields an arboreal differential on $X = \mathbb{CP}^{1} \cup \mathbb{CP}^{1}$, restricting to the original $dz$ on each component.
\end{ex}

We can now understand how differentials degenerate as they approach a facet of their period domain.

\begin{lem}\label{lem:sfacet}
    Let $F$ be a facet of $\overline{\calP_\Gamma^U}$ defined by the condition $\Rea(s_k)=0$ for some transverse saddle period $s_k$ crossing the strip $S$. Then $F$ decomposes as a union of polyhedra of the form $\overline{\calP_{\Gamma'}^U}$, where $\Gamma'$ ranges over the set of critical graphs obtained from $\Gamma$ by collapsing all strips linked with $S$. 
\end{lem}
\begin{proof}
    To understand $F$ defined as in the statement of the lemma, we will consider degenerations of the form $w \colon [0,1) \to \overline{\calP_\Gamma^U}$, with $w(t) \in \calP_\Gamma^U$ for $t > 0$ and $w(0)$ contained in the relative interior of the facet of $\overline{\calP_\Gamma^U}$ defined by $\Rea(s_k) = 0$. 
    
    Treating each period $s_i(t), t_j(t)$ as a function of $t$, we have $\Rea(s_i(0)) = 0$ if and only if the associated saddle lies in a strip linked with $S$ (the ``only if'' direction follows from the assumption that $w(0)$ is contained in the relative interior of the face). Geometrically, this corresponds to each of the strips linked with $S$ collapsing into a chain of vertical saddles as in the discussion at the beginning of the section. 
    
    Now, further suppose $w$ is chosen so that $\Ima(s_i(t)) \ne 0$ for all transverse saddles contained in all strips linked with $S$. As $t \to 0$, the transverse saddles therefore limit onto vertical saddles in a new critical graph $\Gamma'$. Each distinguished saddle for $\Gamma'$ can be identified with a corresponding distinguished saddle for $\Gamma$, though the distinguished orientation may change under this identification, when $\Rea(s_i(t)) \to 0$ and $\Ima(s_i(0)) < 0$. After accounting for this, this gives an embedding $\overline{\calP_{\Gamma'}^U} \into F \subset \overline{\calP_\Gamma^U}$.
    This shows that $F$ decomposes as a union of polyhedra of the form $\overline{\calP_{\Gamma'}^U}$, where $\Gamma'$ ranges over the possible critical graphs produced by such degenerations. 
    \end{proof}

    \begin{rmk}\label{rmk:partialorder}
        Although we will not pursue this in detail, let us comment in passing on how to enumerate the set of possible critical graphs produced by collapsing a given strip. We assume for simplicity that the collapsing strip is not linked with any other; the interested reader is invited to work out the modifications necessary in the linked case. The key observation is that in a strip undergoing collapse, the set of zeros appearing on the left and right sides must be {\em distinct}, for if not, there would appear a purely imaginary absolute period in the limit, which is prohibited by the vertical genericity assumption. In particular, there is a {\em partial order} on the set of zeros of $\Gamma$, with $z_1 < z_2$ whenever there is a chain of vertical saddles beginning at $z_1$ and proceeding upwards to $z_2$. 
        The set of degenerations of a strip can then be seen to be in bijection with the set of {\em interleavings} of the totally-ordered subsets of zeros on the left and right (an interleaving of two disjoint totally-ordered subsets of a poset is a refinement of the partial order that induces a total order on the union). 
    \end{rmk}

We have a corresponding description of the second kind of facet. Here the analysis is simpler, since there is a {\em unique} critical graph appearing in the limit.

\begin{lem}\label{lem:tfacet}
    Let $F$ be a facet of $\overline{\calP_\Gamma^U}$ defined by the condition $\Ima(t_\ell) = 0$. Then $F$ can be identified with the polyhedron $\overline{\calP_{\Gamma'}^U}$, where $\Gamma'$ is the critical graph obtained from $\Gamma$ by collapsing all saddles linked with $\tau_\ell$.
\end{lem}
\begin{proof}
    As in the previous lemma, we understand $F$ by way of a degeneration $w \colon [0,1) \to \overline{\calP_\Gamma^U}$ with $w(t) \in \calP_\Gamma^U$ for $t > 0$ and $w(0)$ in the relative interior of the facet of $\overline{\calP_\Gamma^U}$ defined by $\Ima(t_\ell) = 0$. Under such a degeneration, the corresponding vertical saddle shrinks, until at $t = 0$ the (distinct) zeros at the endpoints have been identified. Moreover, all saddles linked with $\tau_\ell$ are also collapsing. There is a unique resulting critical graph $\Gamma'$, obtained from $\Gamma$ by contracting the corresponding edges. Viewing $\Gamma$ as the $1$-skeleton of a cell structure on $X$, this degeneration gives a cell structure on the quotient of $X$ by the union of the saddles linked with $\tau_\ell$. Note that such a quotient is an arboreal curve.
    As in Lemma \ref{lem:sfacet}, the distinguished saddles in $\Gamma'$ can be identified with corresponding distinguished saddles in $\Gamma$, leading to an embedding $\overline{\calP_{\Gamma'}^U} \into \overline{\calP_{\Gamma}^U}$ as claimed.
\end{proof}

It will be useful to introduce a partial ordering on the set of topological types that appear as degenerations of differentials in $\res^{-1}(U)$ - we call this {\em subordination}.

\begin{defn}[Subordinate]
    The topological types $\Gamma'$ appearing on the boundary of a given $\overline{\calP_\Gamma^U}$ are said to be {\em subordinate} to $\Gamma$; we write $\Gamma' \prec \Gamma$. 
\end{defn}

For given $U \subset B \subset \calR_p$, where $B$ is an open flat and $U$ is a vertically generic subset as usual, we write 
    \[
    \calG(U)
    \]
    to denote the set of topological types $\Gamma$ realized as $\Gamma = \Gamma_\omega$ for $\omega \in \res^{-1}(U)$, and
    \[
    \overline{\calG(U)}
    \]
    for the set of $\Gamma'$ subordinate to some $\Gamma \in \calG(U)$.

The next result summarizes the work of this subsection.
\begin{lem}\label{lem:Pfaces}
    Let $\Gamma, \Gamma' \in \overline{\calG(U)}$ be critical graphs. If $\Gamma' \prec \Gamma$, then there is an inclusion
    \[
    i_{\Gamma'}^\Gamma \colon \overline{\calP_{\Gamma'}^U} \into \overline{\calP_\Gamma^U}
    \]
    with image some polyhedral subset of a face of $\overline{\calP_\Gamma^U}$ defined by the usual constraints $\Rea(s_i) \ge 0, \Ima(t_j) \ge 0$ as well as additional constraints of the form $\pm \Ima(s_i) \ge 0$. Conversely, every face of $\overline{\calP_\Gamma^U}$ decomposes as a union of one or more subfaces of the form $i_{\Gamma'}^\Gamma(\overline{\calP_{\Gamma'}^U})$ for $\Gamma' \prec \Gamma$.
\end{lem}
\begin{proof}
    The work of Lemmas \ref{lem:sfacet}, \ref{lem:tfacet} together prove the claim for the facets of $\overline{\calP_\Gamma^U}$. As any face of codimension $d$ is a facet of a face of codimension $d-1$, successive applications of these results show how to understand an arbitrary face of $\overline{\calP_\Gamma^U}$ as a union of polyhedra $\overline{\calP_{\Gamma'}^U}$ for subordinate $\Gamma' \prec \Gamma$. Under a degeneration from topological type $\Gamma$ to $\Gamma'$, some set of transverse saddles limit onto vertical saddles, and some collection of vertical saddles collapse to points. After reconciling the orientation conventions for transverse and vertical saddles, the upward-pointing condition $\Ima(t_j) \ge 0$ in the limit vertical saddle may correspond to the constraint $\Ima(s_i) \le 0$ on the original transverse saddle, accounting for the possible sign in the additional constraint $\pm \Ima(s_i) \ge 0$.  
\end{proof}
\begin{rmk}
    There is an analogous inclusion map (still notated $i_{\Gamma'}^{\Gamma}$) defined for the compactification. Likewise, there is an inclusion map for the proper bordification, so long as $\Gamma' \in \calG(U)$ and not merely $\overline{\calG(U)}$. 
\end{rmk}

\begin{ex}[Running example, IV]\label{ex:run4}
    Continuing from Examples \ref{ex:run1}, \ref{ex:run2}, \ref{ex:run3} let us view $\omega'$ as being obtained as a deformation of $\omega$, e.g. by swinging $\sigma_2$ down until it degenerates into the vertical saddle $\tau_2$. Using $s_2$ and $s_3$ as coordinates and freezing $s_3$ along this deformation, the periods $s_1, s_4$ of $\sigma_1, \sigma_4$ vary continuously as prescribed by the associated cycle equations. Conversely, if periods $s_1, t_2, s_3, s_4$ are chosen subject to the cycle equations for $\Gamma'$, such values can be used to construct a differential in the boundary of $\calP_\Gamma^\lambda$, by identifying $\tau_2$ as being limited onto by $\sigma_2$ but with orientation reversed.

    The remaining portion of the face $\Rea(s_2) = 0$ in $\calP_\Gamma^\lambda$ corresponds to deformations onto differentials of topological type $\Gamma''$ as shown in Figure \ref{fig:runex2}, where $\sigma_2$ instead limits onto a vertical saddle pointing up. While the vertical saddle in $\Gamma'$ appears on the boundary of two half-space regions, the vertical saddle in $\Gamma''$ appears on the boundary of two strip regions. These strips then have multiple transverse saddles contained within. The inclusion $\calP_{\Gamma''}^\lambda \into \calP_{\Gamma}^\lambda$ identifies the periods $s_1^1,t_2,s_3, s_4^2$ in $\calP_{\Gamma''}^\lambda$ with $s_1, s_2, s_3, s_4$ in $\calP_{\Gamma}^\lambda$, respectively.

\begin{figure}
\centering
		\labellist
    \pinlabel $\omega$ at 435 330
    \pinlabel $\omega''$ at 435 110
    \pinlabel $\Gamma''$ at 380 20
    \Small
    \pinlabel $\sigma_1$ [b] at 49.60 405.32
\pinlabel $\sigma_2$ [b] at 87.30 400.50
\pinlabel $\sigma_3$ [b] at 214.00 405.88
\pinlabel $\sigma_4$ [b] at 250.84 398.23
\pinlabel $\sigma_1^1$ [t] at 49.60 172.33
\pinlabel $\sigma_1^2$ [b] at 48.18 189.90
\pinlabel $\tau_2$ [r] at 71.43 188.20
\pinlabel $\sigma_3$ [b] at 206.48 188.77
\pinlabel $\tau_2$ [l] at 219.47 181.71
\pinlabel $\sigma_4^1$ [t] at 244.32 165.81
\pinlabel $\sigma_4^2$ [b] at 244.89 185.37
		\endlabellist
\includegraphics[width=0.95\textwidth]{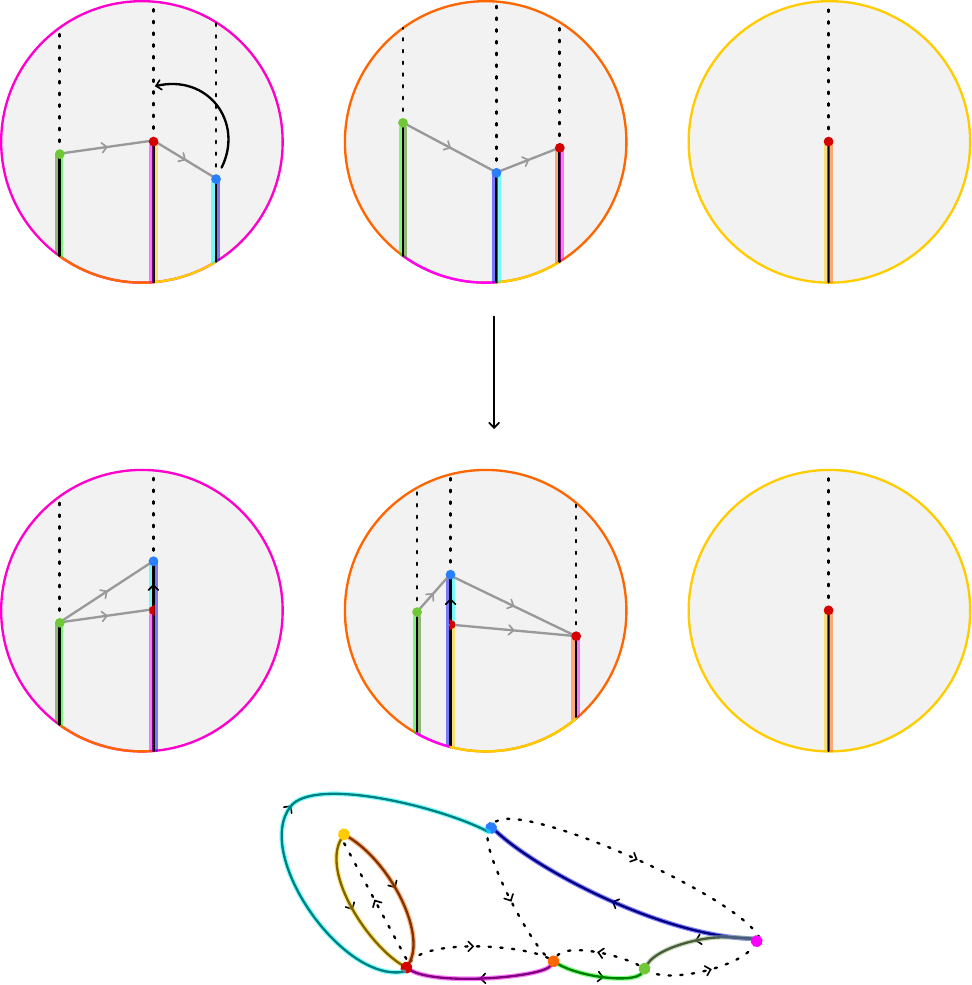}
\caption{Deforming $\omega$ onto $\omega''$ in the boundary of $\calP_\Gamma^\lambda$. }
\label{fig:runex2}
\end{figure}
\end{ex}

\subsection{Period domain: global construction}\label{SS:perioddomainglobal}

Here we come to the definition of the global period domain $\calP^U$ defined over a connected vertically generic subset $U \subset B \subset \calR_p$, by gluing together the polyhedra $\calP_\Gamma^U$ associated to given topological types along the inclusion maps $i_{\Gamma'}^\Gamma$ of Lemma \ref{lem:Pfaces}. From here, we give a proof of one of the main results of the section, Theorem \ref{prop:pUproduct}, which asserts that $\calP^U$ is homeomorphic to a product $\calP^U \cong U \times \calP^\lambda$ for any $\lambda \in U$. 

\begin{defn}[Period domain $\calP^U$]
    Recall from Definition \ref{defn:perioddomainfixedtype} that the proper bordification $\overline{\calP^{U}_{\Gamma}}^{pr}$ is obtained from $\calP_\Gamma^U$ by including faces of the form $\Rea(s_i) = 0$ (corresponding to degenerations that collapse strips, but do not shrink vertical saddles). 
    Define the {\em period domain}
    \begin{equation}\label{eq:PU}
        \calP^U \coloneqq \left.\bigsqcup_{\Gamma \in \calG(U)} \overline{\calP^{U}_{\Gamma}}^{pr}\right/ \sim,
    \end{equation}
    with a given $\overline{\calP^{U}_{\Gamma'}}^{pr}$ identified with its image in each $\overline{\calP^U_\Gamma}^{pr}$ for which $\Gamma' \prec \Gamma$ as in Lemma \ref{lem:Pfaces}.
    Likewise, define the {\em completed period domain}
    \begin{equation}\label{eq:PUcomplete}
        \overline{\calP^U} \coloneqq \left.\bigsqcup_{\Gamma \in \overline{\calG(U)}} \overline{\calP^{U}_{\Gamma}}\right/ \sim,
    \end{equation}
    noting that the union now ranges over the larger set $\overline{\calG(U)}$ of all topological types realized under degeneration.
    Define also the {\em compactified period domain}
    \begin{equation} \label{eq:PUbar}
        \overline{\calP^U}^+ \coloneqq \left.\bigsqcup_{\Gamma \in \overline{\calG(U)}} \overline{\calP^{U}_{\Gamma}}^+\right/ \sim,
    \end{equation}
    where the union ranges over the same set as before. 
      
\end{defn}

We are now in striking distance of the proof of Theorem \ref{prop:pUproduct}. It remains only to make some remarks about the ``face lattice'' of $\overline{\calP^U}^+$.
Since a general inclusion $\overline{\calP_{\Gamma'}^U}^+ \into \overline{\calP_{\Gamma}^U}^+$ only identifies $\overline{\calP_{\Gamma'}^U}^+$ with a polyhedral {\em subface} of $\overline{\calP_{\Gamma}^U}^+$, there is some care required in defining this notion.

\begin{defn}[Face lattice of the (compactified) period domain]
    A (nonempty) face $F$ of some $\overline{\calP^U_\Gamma}^+$ is said to be {\em imprimitive} if there is some inclusion $i_{\Gamma'}^\Gamma \colon \overline{\calP_{\Gamma'}^U}^+ \into \overline{\calP_\Gamma^U}^+$ and a face $F'$ of $\overline{\calP_{\Gamma'}^U}^+$ with image $i_{\Gamma'}^\Gamma(F') \subset F$ a proper subset of $F$; otherwise a face is said to be {\em primitive}. The {\em face lattice} $\overline{\Lambda^U}^+$ of $\overline{\calP^U}^+$ is the lattice of primitive faces of the constituent polytopes $\overline{\calP^U_\Gamma}^+$ under inclusion, modulo the identifications of \eqref{eq:PUbar}.  
\end{defn}

\begin{lem}\label{lem:facelattice}
Every minimal element of $\overline{\Lambda^U}^+$ is a vertex, i.e. a face of some $\overline{\calP_\Gamma^U}^+$ of dimension zero.    
\end{lem}
\begin{proof}
    In general, let $P$ be a spherical polytope with face lattice $\Lambda_P$ and lineality space $\PP^+W \subset P$ of dimension $d \ge -1$ (interpreting $d=-1$ as the case where the lineality space is empty). If $d=-1$ then the minimal elements of $\Lambda_P$ are the vertices of $P$. Otherwise every face of $P$ contains $\PP^+W \cong S^d$. In the case $d = 0$, we can view $\PP^+W \cong S^0$ as a union of two vertices, while for $d > 0$, the face lattice $\Lambda_P$ contains a unique minimal element $\PP^+W \cong S^d$. 

    Specializing now to our setting, we consider the complex $\overline{\calP^U}^+$, as assembled from the spherical polytopes $\overline{\calP_\Gamma^U}^+$. While by Lemma \ref{lem:Plineality}, the lineality space of a general $\overline{\calP_\Gamma^U}^+$ will be nonempty, we claim that these are always imprimitive. 
    
    To see this, recall from Lemma \ref{lem:Plineality} that the (unprojectivized) lineality space $W_\Gamma^\lambda$ is the span of the variables $\{\Ima(s_i)\}$ contained in $V_\Gamma^\lambda$. This is nonempty if and only if $\Gamma$ contains at least one strip region. In this case, there is at least one strip that can be collapsed as in Lemma \ref{lem:sfacet} (as can be seen by choosing coordinates adapted to some transverse saddle $\sigma_i$ as in Lemma \ref{lem:adaptedcoordinates}). Let $\Gamma'$ be a critical graph obtained by some such collapse; for simplicity suppose the degeneration arises by the transverse saddle $\sigma_i$ in $\Gamma$ degenerating onto a vertical saddle $\tau_i$ in $\Gamma'$, with the canonical orientations of $\sigma_i$ and $\tau_i$ agreeing. Let $s_i, t_i$ denote the respective periods.
    
    Such a degeneration induces an inclusion $i_{\Gamma'}^\Gamma \colon \overline{\calP_{\Gamma'}^U}^+ \into \overline{\calP_{\Gamma}^U}^+$ with image contained in the facet $\Rea(s_i) = 0$ of $\overline{\calP_{\Gamma}^U}^+$. This containment is proper, as the image does not contain any point with $\Ima(s_i) < 0$. This shows that the face $\Rea(s_i)=0$ is imprimitive, and shows moreover that the lineality space $\PP^+W_\Gamma^\lambda$ is imprimitive, as the inclusion $i_{\Gamma'}^\Gamma$ subdivides this according to the sign of $\Ima(s_i)$.
\end{proof}

\begin{proof}[Proof of Theorem \ref{prop:pUproduct}]
    The overall idea is quite straightforward, and can be summarized in one paragraph. Theorem \ref{lem:Pproduct} implies that $\overline{\calP^U}^+$ forms a family of almost-polyhedral complexes $\overline{\calP^\lambda}^+$ for $\lambda \in U$ whose combinatorial structure is independent of $\lambda \in U$. We endow each polyhedral piece $\overline{\calP^\lambda}^+$ with a triangulation that is independent of $\lambda \in U$, giving a canonical family of identifications as $\lambda \in U$ varies, via the piecewise-linear homeomorphism that is linear on each simplex. Assembling these across $\lambda \in U$ yields a homeomorphism $\overline{\calP^U}^+ \cong U \times \overline{\calP^{\lambda}}^+$ for any fixed $\lambda \in U$; this restricts to give the required homeomorphism $\calP^U \cong U \times \calP^\lambda$.

  This triangulation process requires some further comment.
    We will make use of the following inductive principle: for any spherical polytope $P$ of dimension $d$, suppose the boundary $\partial P$ has been given the structure of a simplicial complex of dimension $d-1$. Then $P$ can be decomposed into $d$-simplices by taking the cones on the $(d-1)$-simplices in $\partial P$ with apex at the barycenter of $P$. 
    
    By Lemma \ref{lem:facelattice}, every minimal element of the face lattice $\overline{\Lambda^\lambda}^+$ is zero-dimensional and hence is canonically triangulated. Each face of $\overline{\Lambda^\lambda}^+$ corresponds to some spherical polytope $\overline{\calP_\Gamma^\lambda}^+$, and hence its boundary is a union of spherical polytopes of lower dimension, corresponding to not necessarily primitive faces. However, by Lemma \ref{lem:Pfaces}, each imprimitive face is a union of primitive faces. Therefore, $\overline{\calP^\lambda}^+$ can be triangulated by invoking the above principle.
    
    We now make the following three observations, which will conclude the argument.
   Firstly, this simplicial structure depends only on the combinatorics of the polytopes $\overline{\calP^\lambda_\Gamma}^+$, so that the isomorphism type of $\overline{\calP^\lambda}^+$ as a simplicial complex is independent of $\lambda$. Secondly, by construction, each face of each $\overline{\calP^\lambda}^+$ determines a {\em subcomplex} of $\overline{\calP^\lambda}^+$, so that a simplicial homeomorphism $\overline{\calP^\lambda}^+ \cong \overline{\calP^{\lambda'}}^+$ restricts to a homeomorphism $\calP^\lambda \cong \calP^{\lambda'}$. Finally, note that the barycenter of a polytope varies continuously as a function of the vertices, so that this family of homeomorphisms assembles into global homeomorphisms
   \[
   \overline{\calP^U}^+ \cong U \times \overline{\calP^{\lambda}}^+ \quad \mbox{and} \quad \calP^U \cong U \times \calP^{\lambda}
   \]
   as claimed.
\end{proof}

\subsection{Period mappings (global construction)}\label{SS:periodmap}
Finally in this section, we define the global period map $\per^U \colon \res^{-1}(U) \to \calP^U$ and prove that it is a homeomorphism (Theorem \ref{thm:periodmaphomeo}). As always, $B\subset \calR_p$ denotes an open flat, and $U \subset B$ is a connected vertically generic subset.

\begin{defn}[Period map]
    $\res^{-1}(U)$ admits a stratification
    \[
    \res^{-1}(U) = \bigsqcup_{\Gamma \in \calG(U)} \res^{-1}(U)_{\Gamma}
    \]
    according to topological type. Using this, define the {\em period map}
    \[
    \per^U \colon \res^{-1}(U) \to \calP^U
    \]
    by defining
    \[
    \per^U |_{\res^{-1}(U)_{\Gamma}} \coloneqq \per^U_{\Gamma},
    \]
    identifying the target $\calP^U_{\Gamma}$ of $\per^U_{\Gamma}$ with its image in $\calP^U$.
\end{defn}

\begin{proof}[Proof of Theorem \ref{thm:periodmaphomeo}]
    By Lemma \ref{lem:periodhomeolocal}, the restriction of $\per^U$ to any subset $\res^{-1}(U)_\Gamma \subset \res^{-1}(U)$ gives a homeomorphism $\res^{-1}(U)_\Gamma \cong \calP_\Gamma^U$. By Lemma \ref{lem:Pfaces}, this extends to a homeomorphism
    \[
    \overline{\res^{-1}(U)_\Gamma} \to \overline{\calP_\Gamma^U}^{pr},
    \]
    where $\overline{\res^{-1}(U)_\Gamma}$ denotes the closure of $\res^{-1}(U)_\Gamma$ in $\res^{-1}(U)$. These homeomorphisms are compatible with the identifications used to construct $\calP^U$ from the disjoint union of the polyhedral pieces $\overline{\calP_\Gamma^U}^{pr}$ as $\Gamma$ ranges over the set of realized topological types $\calG(U)$. This shows that the various homeomorphisms $\per^U_\Gamma \colon \res^{-1}(U)_\Gamma \to \calP_\Gamma^U$ assemble into a global homeomorphism $\per^U \colon \res^{-1}(U) \to \calP^U$.
\end{proof}

\section{Estimating Betti numbers}\label{sec:Betti} 

Here we turn to the question of estimating the Betti numbers $\operatorname{Betti}_i(\mathcal{F}_\lambda(\mu))$ of an isoresidual fiber $\mathcal{F}_\lambda(\mu)$ in a moduli space $\mathcal{H}(\mu)$ of differentials of genus zero, as well as the Betti numbers of the metric completion $\overline{\calF}_\lambda(\mu)$. In this section, we fix $\mu = (a_1, \dots, a_n, -b_1, \dots, -b_p)$ for $n + p \ge 3$ and set $B\coloneqq \sum_{j =1}^p b_j$.

\subsection{Guiding ideas}\label{sub:Ideas}

While the estimate we obtain below is rather crude, it has an important structural feature: it separates out the contributions coming from the purely {\em topological} constants $n$ and $p$ from the contributions due to the {\em geometric} constants, i.e. the orders $a_i,b_j$ of the zeros and poles. The dimension (over $\C$) of $\mathcal{F}_\lambda(\mu)$ is given by $n-1$ (and is therefore a purely topological invariant). It is natural to ask how the Betti numbers grow when the dimension (or indeed both constants $n$ and $p$) are fixed, and the orders of the zeros and poles grow. This is answered by Theorem~\ref{thm:BettiEstimate}, which bounds the Betti numbers in terms of a polynomial in $B$ (the total orders of the poles) of degree at most $n+p-3$. 
\par
To prove this result, we may assume that $\lambda$ is vertically generic, so that by Theorem \ref{thm:periodmaphomeo}, $\mathcal{F}_\lambda(\mu) \cong \calP^\lambda$. We can compute the homology of $\mathcal{F}_\lambda(\mu)$ via a cell structure {\em dual} to the polyhedral structure of $\calP^\lambda$, so that a polyhedron of real codimension $i$ (relative to the dimension of a top-dimensional polyhedron) corresponds to a cell of dimension $i$. As the Betti number is bounded by the number of cells in the corresponding dimension, it suffices to bound the number of critical graphs that appear among differentials in $\mathcal{F}_\lambda(\mu)$. 
\par
We first require some basic information about the structure of a critical graph.

\begin{lem}\label{lem:gammacount}
Let $\Gamma$ be the critical graph of a differential in an isoresidual fiber $\mathcal{F}_{\lambda}(\mu)$, and suppose that $\Gamma$ has $r$ vertical saddle connections. Then the number of vertices and edges of $\Gamma$ are given by
        \[
        |v(\Gamma)| = n+p, \quad |e(\Gamma)| = 2B + 2n - 4 - r.
        \]
\end{lem}

\begin{proof}
The count $|v(\Gamma)| = n+p$ is clear from the definition of $\Gamma$. To determine $|e(\Gamma)|$, observe that each edge of $\Gamma$ has at least one endpoint at a vertex corresponding to a zero, with {\em both} endpoints at a zero if and only if it is a vertical saddle connection. A zero of order $a_i$ has $2a_i+2$ vertical critical leaves incident, for a total of $2A + 2n = 2B+2n - 4$ half-edges. A half-edge is double-counted if and only if it is a vertical saddle, so that $|e(\Gamma)| = 2B + 2n -4 -r$ as claimed.
\end{proof}

\subsection{Weighted critical graphs}

    To separate out the contributions to the number of critical graphs coming from $n,p$ alone from those coming from $B$, we introduce the notions of a {\em weighted critical graph} and a {\em reduced critical graph}. 
    
\begin{figure}
\centering
		\labellist
\pinlabel $\Gamma$ [tl] at 168.93 141.44
\pinlabel $\overline{\Gamma}$ [tl] at 181.97 13.89
\tiny
\pinlabel 2 [bl] at 33.73 62.64
\pinlabel 1 [b] at 64.62 40.82
\pinlabel 1 [bl] at 123.58 38.83
\pinlabel 1 [tl] at 166.38 38.26
\pinlabel 1 [bl] at 117.63 70.86
\pinlabel 0 [b] at 41.38 92.12
\pinlabel 0 [tr] at 80.21 65.19
\pinlabel 0 [b] at 121.03 2.83
		\endlabellist
\includegraphics[scale=1]{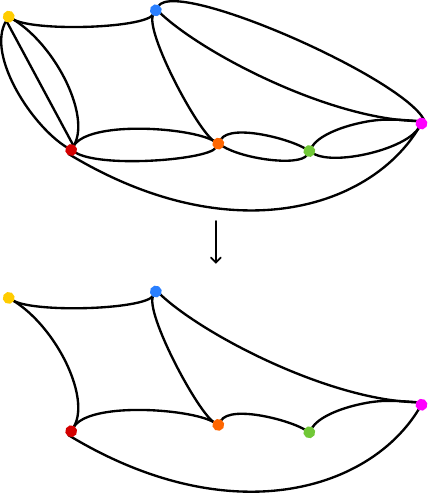}
\caption{Passing from the critical graph $\Gamma$ to the reduced critical graph $\overline{\Gamma}$.}
\label{fig:reducedcrit}
\end{figure}

    \begin{defn}[Weighted, reduced critical graph]
        Let $\Gamma$ be a critical graph. The {\em weighted critical graph} is a graph $\overline{\Gamma}$ together with a weight function $w \colon e(\overline{\Gamma}) \to \Z_{\ge 0}$, obtained from $\Gamma$ by the following process.

        Begin by assigning a weight of zero to each edge of $\Gamma$. The complement of $\Gamma$ has some number of bigon regions, corresponding to half-space regions in a neighborhood of a pole, that are bounded by a union of exactly two critical leaves. To construct $\overline{\Gamma}$, successively {\em collapse} the two edges bounding a bigon region into a single edge, and increase the weight of that edge by $1$. See Figure \ref{fig:reducedcrit}.

        The {\em reduced critical graph} $\overline \Gamma$ associated to $\Gamma$ is just the graph $\overline \Gamma$ formed by the above process, without the weight function.
    \end{defn}

\begin{rmk}\label{rmk:weightcrit}
    Here are some simple but important observations about the structure of a weighted/reduced critical graph.
    \begin{enumerate}
        \item The data encoded by a critical graph and a weighted critical graph are equivalent: given $(\overline \Gamma, w)$, one can uniquely reconstruct $\Gamma$ by splitting each weighted edge into $w(e) + 1$ edges bounding bigons in between.
        \item Recall that the vertices of $\Gamma$ correspond to zeros and poles of a differential, and by convention are globally labeled. Thus the vertices of $\overline{\Gamma}$ are likewise labeled, and we consider two weighted (or reduced) critical graphs to be distinct if the labels assigned to the vertex sets are different, even if the (weighted) graphs themselves are otherwise isomorphic.
        \item All edges of $(\overline \Gamma, w)$ corresponding to vertical saddles have weight zero, since each edge of a bigon necessarily runs from a zero to a pole and so all edges of positive weight are of this form. In particular, the graph $\overline{\Gamma}'$ obtained from $\overline \Gamma$ by deleting all vertical saddle edges is {\em bipartite}, since each remaining edge connects a zero vertex to a pole vertex.
    \end{enumerate}
\end{rmk}

In a critical graph $\Gamma$, the order of a zero is encoded by the valence at the corresponding vertex (it is $2a_i + 2$), and the order of a pole is encoded by the number of half-space regions adjacent to the corresponding vertex (there are $2b_i - 2$). The weight structure in a weighted critical graph likewise allows these invariants to be recovered. For zeros this is straightforward, but for poles, care must be taken to account for so-called {\em special half-spaces}. A half-space is {\em special} if there is some vertical saddle connection appearing on the boundary (and is {\em ordinary} otherwise). The complementary region in $\Gamma$ for a special half-space is not a bigon, and so does not get collapsed into a weighted edge. These observations can be summarized as follows.

\begin{lem}[Weight conditions]\label{lem:weights}
Let $\Gamma$ be a critical graph for the moduli space $\mathcal{H}(\mu)$, and let $(\overline{\Gamma}, w)$ be the associated weighted critical graph. Extend $w$ to a function on the vertices of $\overline{\Gamma}$ by taking $w(v)$ to be the sum of the weights of edges incident to $v$, and let $\val_{\overline{\Gamma}}(v)$ denote the valence of $v$ as a vertex of $\overline \Gamma$. Then $w$ satisfies the following {\em weight conditions} at each vertex:
\begin{itemize}
    \item At a vertex $v$ corresponding to a zero of order $a_i$, 
    \[
    w(v) + \val_{\overline{\Gamma}}(v) = 2a_i + 2.
    \]
    \item At a vertex $v$ corresponding to a pole of order $b_j$, let $s_j$ denote the number of special half-spaces at the pole. Then
    \[
    w(v) + s_j = 2 b_j - 2.
    \]
\end{itemize}
Moreover, the total weight is given by
\[
\sum_{e \in \overline \Gamma} w(e)= 2B-2p - s,
\]
where $s = \sum s_j$ denotes the total number of special half-spaces.
\end{lem}
\begin{proof}
    The total weight is just given by the number of regions collapsed in the process of forming $(\overline{\Gamma},w)$, i.e. the number of ordinary half-spaces, which is easily counted to be $2B-2p - s$. To see that the weight conditions hold, observe that the quantities on the left-hand sides of both equations are unchanged by the process of collapsing and weighting edges. They are trivially satisfied for $\Gamma$ itself (with identically zero weighting function), and so also hold for $(\overline \Gamma, w)$. 
\end{proof}

Our next objective is to estimate the number of possible weight functions on a given reduced critical graph $\overline \Gamma$, that produce a potential critical graph for $\mathcal{H}(\mu)$, i.e. that at a minimum satisfy the weight conditions of Lemma \ref{lem:weights}. To do so, we momentarily consider real-valued weight functions satisfying the weight conditions, which therefore form an affine space. To determine its dimension, we first analyze the case when $\Gamma$ has no vertical saddles in Lemma \ref{lem:weightspacer0}, and then parlay this into the general case in Lemma \ref{lem:weightspace}. We must first formulate and prove one graph-theoretic result, Lemma \ref{lem:Aker}.

Let $\Gamma$ be an arbitrary graph, and let $\R^{v(\Gamma)}$ denote the space of real-valued functions on the vertices of $\Gamma$, and define $\R^{e(\Gamma)}$ similarly as the space of real-valued functions on the edges; we equip each with the canonical basis of vertices and edges, respectively. Let $A \colon \R^{e(\Gamma)} \to \R^{v(\Gamma)}$ be the unsigned adjacency map, so that if $e$ is an edge connecting vertices $u,v$, then $A(e) = u + v$. The adjoint $A^{\ast} \colon \R^{v(\Gamma)} \to \R^{e(\Gamma)}$ then has $A^{\ast}(v) = \sum_{v \in e} e$.
\begin{lem}\label{lem:Aker}
    In the above setting, 
    \[
    \dim(\ker(A^*)) = b,
    \]
    where $b$ is the number of bipartite components of $\Gamma$.
\end{lem}
\begin{proof}
    A function $y = \sum y_v v \in \R^{v(\Gamma)}$ satisfies $A^*(y) = 0$ if and only if for each edge $e$, say connecting $u$ to $v$, the condition $y_u + y_v = 0$ holds. Now observe that if $u$ and $v$ are in the same component of $\Gamma$, then $y_u \ne 0$ if and only if $y_v \ne 0$, so that the support of $y$ must be some union of connected components of $\Gamma$, and that moreover each such component must be bipartite (since each edge must connect a vertex $v$ with $y_v > 0$ to a vertex $u$ with $y_u < 0$). From this it is easy to deduce the claimed dimension of $\ker(A^*)$.
\end{proof}

\begin{lem}\label{lem:weightspacer0}
    Let $\Gamma$ be a critical graph for the moduli space $\mathcal{H}(\mu)$ with no vertical saddles, and let $\overline \Gamma$ be the associated reduced critical graph. Let $W(\overline \Gamma, \mu) \subset \R^{|e(\overline \Gamma)|}$ be the affine space of weight functions that satisfy the weight conditions of Lemma \ref{lem:weights}. Then $\dim(W(\overline \Gamma, \mu)) = n+p-3$.
\end{lem}
\begin{proof}
    A weight function $w \in \R^{e(\overline \Gamma)}$ satisfies the weight conditions of Lemma \ref{lem:weights} if and only if it satisfies the equation $A(w) = y$, where $y(v) = 2a_i+2-\val_{\overline{\Gamma}}(v)$ if $v$ is a vertex corresponding to a zero of order $a_i$, and $y(v) = 2b_j - 2 - s_j$ if $v$ is a vertex corresponding to a pole of order $b_j$ abutted by $s_j$ special half-spaces; since $\Gamma$ contains no vertical saddles, there are no special half-spaces to consider.
    As $\overline \Gamma$ is assumed to come from a critical graph $\Gamma$, the weight conditions are consistent, and so the dimension of $W(\overline \Gamma, \mu)$ is therefore computed as
    \[
    \dim(W(\overline \Gamma,\mu)) = \dim(\ker(A)) = \dim(\operatorname{coker}(A^*)) = |e(\overline \Gamma)|- |v(\overline \Gamma)| + \dim(\ker(A^*)).
    \]

    By Remark \ref{rmk:weightcrit}.3, $\overline \Gamma$ is bipartite (and connected, since $\Gamma$ is), and hence $\dim(\ker(A^*)) = 1$ by Lemma \ref{lem:Aker}.
    To prove the result, we must show that $|e(\overline \Gamma)|-|v(\overline \Gamma)| = n + p - 4$.
    Firstly, $|v(\overline \Gamma)| = n+p$. To determine $|e(\overline \Gamma)|$, recall from Lemma \ref{lem:gammacount} that $|e(\Gamma)| = 2B + 2n - 4 -r$ (of course, here $r = 0)$. In passing to $\overline{\Gamma}$, we lose $2B-2p-s$ edges corresponding to the collapsed ordinary half-spaces (and here $s = 0$), so that
    \[
    |e(\overline \Gamma)| = (2B + 2n - 4) - (2B-2p) = 2n + 2p - 4 
    \]
    and so
    \[
    |e(\overline \Gamma)| - |v(\overline \Gamma)| = n + p - 4
    \]
    as required.
\end{proof}

\begin{figure}
\centering
		\labellist
        \pinlabel $\Gamma$ [tl] at 70.58 186.79
\pinlabel $\Gamma_0$ [tl] at 200.67 190.19
\pinlabel $\overline{\Gamma}'$ [tl] at 70.58 24.66
\pinlabel $\overline{\Gamma_0}'$ [tl] at 200.67 24.66
\pinlabel S [c] at 28.34 225.33
\pinlabel $e_1$ [tr] at 43 250.28
\pinlabel $e_2$ [br] at 44.22 204
		\endlabellist
\includegraphics[scale=1]{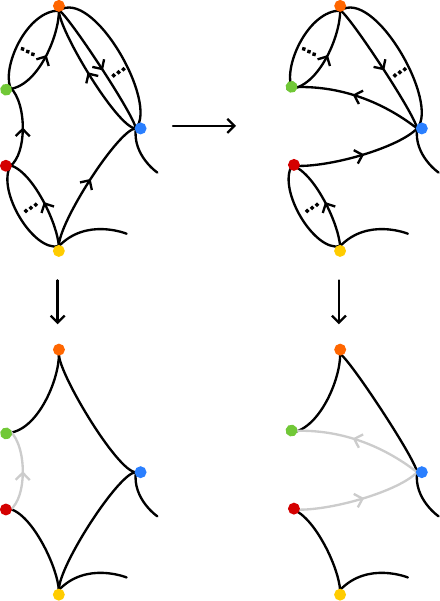}
\caption{An example of the graphs $\Gamma, \Gamma_0, \overline{\Gamma}'$, and $\overline{\Gamma_0}'$ - as $S$ collapses, the zero vertex on the right merges in between the two on the left, destroying edges $e_1, e_2$ in the process. The vertical saddles removed from $\overline{\Gamma}'$ and $\overline{\Gamma_0}'$ are shown in light gray for reference, but are not edges in the graphs on the bottom row (though they are edges in $\overline{\Gamma}$ and $\overline{\Gamma_0}$). Here, $e_2$ bounds no bigons and so disappears when passing from $\overline{\Gamma}'$ to $\overline{\Gamma_0}'$, while $e_1$ does bound some bigon and so survives to $\overline{\Gamma_0}'$.}
\label{fig:stripcollapse}
\end{figure}

\begin{lem}\label{lem:weightspace}
    Let $\Gamma$ be a critical graph for the moduli space $\mathcal{H}(\mu)$, and let $\overline \Gamma$ be the associated reduced critical graph. Let $W(\overline \Gamma, \mu) \subset \R^{|e(\overline \Gamma)|}$ be the affine space of weight functions that satisfy the weight conditions of Lemma \ref{lem:weights} and moreover satisfy $w(e) = 0$ for any edge of $\overline\Gamma$ coming from a vertical saddle connection. Then $\dim(W(\overline \Gamma,\mu)) \le n+p-3$.
\end{lem}
\begin{proof}
    Let $\overline{\Gamma}'$ be the graph obtained from $\overline \Gamma$ by deleting every edge coming from a vertical saddle connection. Then, following Remark \ref{rmk:weightcrit}.3, $W(\overline \Gamma, \mu)$ can be canonically identified with $W(\overline \Gamma', \mu)$, the space of real-valued weight functions on $\overline \Gamma'$ satisfying the same weight conditions. 
    
    Let $\Gamma_0$ be a critical graph obtained from $\Gamma$ by collapsing a strip, as in Section \ref{SS:degenerating}, and define $\overline{\Gamma_0}$ and $\overline{\Gamma_0}'$ as above. We will show that there is an embedding of affine spaces 
    \[
    W(\overline{\Gamma_0}, \mu) \into W(\overline{\Gamma}, \mu).
    \]
    Since every critical graph $\Gamma$ can be obtained from a critical graph $\tilde \Gamma$ with no vertical saddles by collapsing some sequence of strips, the claim will then follow from Lemma \ref{lem:weightspacer0}.

    We claim that $\overline{\Gamma_0}'$ naturally embeds as a subgraph of $\overline{\Gamma}'$. As Figure \ref{fig:stripcollapse} shows, when a strip $S$ in $\Gamma$ collapses and creates $\Gamma_0$, two of the four infinite ray edges of $S$, say $e_1$ and $e_2$, are deleted, and the vertices on either side of $S$ are joined by some sequence of vertical saddles as in Remark \ref{rmk:partialorder}. In passing to $\overline{\Gamma}$ and $\overline{\Gamma_0}$, if $k_i$ bigons merge onto $e_i$ in $\overline{\Gamma}$, then $k_i -1$ bigons merge onto $e_i$ in $\overline{\Gamma_0}$. Finally when passing to $\overline{\Gamma}'$ and $\overline{\Gamma_0}'$, all vertical saddle edges are deleted. The result is that $\overline{\Gamma_0}'$ is obtained from $\overline{\Gamma}'$ by deleting $e_1$ or $e_2$ only if there are $k_i= 0$ bigons merging onto $e_i$ in $\overline{\Gamma}$. 

    This moreover shows how to construct the embedding $W(\overline{\Gamma_0}, \mu) \into W(\overline{\Gamma}, \mu)$. Given a weight function $w_0 \in W(\overline{\Gamma_0}, \mu)$, define $w \in \R^{e(\overline{\Gamma}')}$ by taking $w(e) = w_0(e)$ for $e \ne e_1, e_2$. For $i = 1,2$, if $e_i \in \overline{\Gamma_0}'$, define $w(e_i) = w_0(e_i) +1$, and if $e_i \not \in \overline{\Gamma_0}'$, define $w(e_i) = 0$. Then it is easy to verify that $w$ satisfies the weight conditions defining $W(\overline{\Gamma}, \mu)$: at the endpoint $z_i$ of $e_i$ corresponding to a zero, $w(z_i) = w_0(z_i) + 1$, but 
    \[
    \val_{\overline{\Gamma_0}}(z_i) = \val_{\overline{\Gamma}}(z_i) + 1,
    \]
    since there is a new vertical saddle in $\overline{\Gamma_0}$ at $z_i$. At any other zero vertex, both the total weight and valence are unchanged. The analysis at a pole is similar: the loss in weight is offset by the addition of a new special half-space. 
 \end{proof}

\subsection{Counting unweighted critical graphs}

We are now in striking distance of the proof of Theorem~\ref{thm:BettiEstimate}. We must prove just one more preparatory result, which shows that the number of unweighted reduced critical graphs appearing for a given moduli space $\mathcal{H}(\mu)$ depends only on $n,p$.

\begin{lem}\label{lem:Cnp}
Let $C_{\mu,\lambda}$ denote the number of (unweighted) reduced critical graphs $\overline \Gamma$ obtained from critical graphs $\Gamma$ appearing in $\overline{\mathcal{F}_\lambda(\mu)}$. Then 
    \[
    C_{\mu, \lambda} \le C_{n,p},
    \]
    where $C_{n,p}$ is a constant depending only on $n,p$.
\end{lem}

\begin{proof}
Following the computation of $|e(\overline{\Gamma})|$ in the proof of Lemma \ref{lem:weightspacer0}, we see that any $\overline{\Gamma}$ is a graph with at most $n+p$ vertices and at most $2n+2p + s - r - 4$ edges, where $r$ is the number of vertical saddles and $s$ is the number of special half-spaces (as we consider $\overline{\calF_\lambda}(\mu)$ here and not merely $\calF_\lambda(\mu)$, the number of zeros in $\overline{\Gamma}$ may be less than $n$). Note that every special half-space contains some vertical saddle in its boundary, and conversely a given vertical half-space can appear in the boundary of at most two special half-spaces, so that $s \le 2r$. Thus
\[
|e(\overline{\Gamma})| \le 2n + 2p -4 + r.
\]
We claim that $r$ is bounded as a function of $n,p$ alone. To see this, we first observe that since $\Gamma$ is embedded in $\mathbb{CP}^{1}$ (or more generally in an arboreal curve $X$ given topologically as a wedge of $2$-spheres), the number of isotopy classes of edges corresponding to vertical saddles is bounded in terms of $n, p$ alone (since their interiors must be pairwise disjoint). We moreover claim that no two vertical saddle edges can form a bigon. Any such pair of vertical saddle edges forming a bigon would give a closed loop whose period would be nonzero and purely imaginary, but the bigon region contains no poles and hence the period would be zero, a contradiction.
\end{proof}

\begin{proof}[Proof of Theorem~\ref{thm:BettiEstimate}]
As discussed in Section \ref{sub:Ideas}, it suffices to estimate the number of cells of codimension $i$ in the period domain $\calP^\lambda$ and its metric completion $\overline{\calP^\lambda}$.  
We first show that in the range $0 \le i \le 2n-2$, the number of cells in codimension $i$ (and hence the corresponding Betti numbers) of the completion $\overline{\calP^\lambda}$ satisfies the claimed bound. The corresponding result for $\calP^\lambda$ then follows, since it embeds as a subcomplex of the completion.
The cells of $\overline{\calP^\lambda}$ are indexed by critical graphs appearing for differentials limited onto from $\calH(\mu)$ and hence have at most $n$ zeros and at most $p$ poles. By Lemma~\ref{lem:Cnp}, the number of possible unweighted reduced critical graphs that appear is a constant $C_{n,p}$ depending only on $n,p$, as claimed; this constant will change depending on whether we are considering critical graphs appearing only in $\calP^\lambda$, or its completion $\overline{\calP^\lambda}$, but in either case it depends on $n,p$ alone.

It remains to show that for a given unweighted reduced critical graph $\overline \Gamma$, the number of associated critical graphs $\Gamma$ (i.e. the number of weightings) with the correct orders of zeros and poles, is bounded above by $(2B-2p+1)^{n+p-3}$. 
\par
By Lemma~\ref{lem:weights}, the total weight of any $\Gamma$ reducing to $\overline{\Gamma}$ is $2B-2p - s$, which may be trivially bounded above by $2B-2p$. By Lemma~\ref{lem:weightspace}, the space of real-valued weight functions $w \colon e(\overline{\Gamma}) \to \R$ satisfying the weight conditions at each vertex has dimension at most $n+p-3$. As each weight is an integer between $0$ and the total weight, it follows that there are at most $2B-2p+1$ possible values that can be assigned independently to each of $n+p-3$ edges, for a maximum total of $(2B-2p+1)^{n+p-3}$ possible weight functions.
\par

We next show that the dual cell structure for $\calP^\lambda$ has no cells in dimension above $n-1$, i.e. that $\calP^\lambda$ has no polyhedra of codimension $i \ge n$. This will establish the claimed vanishing $\operatorname{Betti}_i(\calF_\lambda(\mu)) = 0$ for $i \ge n$. This follows from the theory of adapted coordinates (cf. Lemma \ref{lem:adaptedcoordinates}): any polyhedron $\calP^\lambda_\Gamma$ has coordinates coming from the periods of a set of $n-1$ distinguished saddles that form a spanning tree for the graph $\Sigma_\Gamma$ of distinguished saddles. Each transverse saddle among this spanning tree contributes two parameters (its real and imaginary parts), while each vertical saddle contributes only one. The codimension of $\calP^\lambda_\Gamma$ is therefore equal to the number of vertical saddles among this set, and so is at most $n-1$. The vanishing of Betti numbers of the completion $\overline{\calP^\lambda}$ above dimension $2n-2$ is a triviality, since $\overline{\calP^\lambda}$ is a union of polyhedra of real dimension $2n-2$.
\end{proof}

\section{Isoresidual Hurwitz fibers}
\label{sec:isoresidual-hurwitz-fibers}

We present an application of our main results, using here a genus-zero interpretation of certain isoresidual fibers as Hurwitz-type spaces. This gives a useful way to translate the stratified fibration theorem into a statement about families of rational maps. It also provides a conceptual bridge between the isoresidual fibration and the stratified polynomial spaces studied in \cite{SaBraid, SaBraid2, SaHuxBraid}.

We use the standard residue convention. Thus, if a local coordinate satisfies $z=w^m$, then
\[
\frac{dz}{z} = m\frac{dw}{w}
\]
has residue $m$ at $w=0$.  We say that a residue vector is \emph{integral} if its coordinates all lie in $\ZZ$.

Consider a stratum of differentials $\mathcal{H}(\mu)$ with $\mu = (a_1,\ldots, a_n, -1^p)$, $a_{1},\dots,a_{n} \geq 1$ and $\sum_{i=1}^n a_i = p-2$ parameterizing differentials on $\mathbb{CP}^{1}$ with all poles being simple. Given its residual map $\res\colon \calH(\mu)\to \calR_p$, for $\lambda \in \calR_{p}$, the corresponding isoresidual fiber is $\mathcal{F}_\lambda(\mu) = \res^{-1}(\lambda)$. 
\par
Since all poles are simple,  $\mathcal{F}_\lambda(\mu)$ is empty whenever some $\lambda_j = 0$.

\begin{defn}[Isoresidual Hurwitz fiber]
\label{def:IH-fiber}
Let $\lambda\in \calR_p$ be integral and suppose that all coordinates of $\lambda$ are nonzero.  The \emph{isoresidual Hurwitz fiber} associated to $(\mu,\lambda)$ is
\[
  \cIH_\mu(\lambda)\coloneqq \mathcal{F}_\lambda(\mu),
\]
viewed through the identification with rational maps given in Proposition~\ref{prop:df-f-dictionary} below.  Equivalently, it is the moduli space of rational maps $f\colon \mathbb{CP}^1\to\mathbb{CP}^1$, modulo automorphisms of the source and multiplication of the target coordinate by $\CC^*$, such that
\[
  \frac{df}{f}\in \calH(\mu), 
  \quad  \res\left(\frac{df}{f}\right) = \lambda.
\]
Here the branch points of $f$ away from $0$ and $\infty$ are {\em not} fixed and are allowed to collide.
\end{defn}

We write $\lambda^+ = (\alpha_1,\ldots,\alpha_r)$ and $\lambda^- = -(\beta_1,\ldots,\beta_s)$
for the positive and negative parts of the integer vector $\lambda$, and let 
\[
d(\lambda) = \sum_{i=1}^r\alpha_i = \sum_{j=1}^s\beta_j.
\]
This is the {\em degree} of the rational map associated to a point of $\cIH_\mu(\lambda)$.

\begin{prop}[Isoresidual fiber with nonzero integral residues]
\label{prop:df-f-dictionary}
Let $\lambda\in \calR_p$ be integral with all coordinates nonzero.  Then every differential in $\calF_\lambda(\mu)$ is uniquely of the form $df/f$, up to multiplying $f$ by a constant in $\CC^*$.  Under this identification, $\cIH_\mu(\lambda)$ parametrizes rational maps $f\colon \mathbb{CP}^1\to\mathbb{CP}^1$ such that:
\begin{enumerate}[label=\textup{(\roman*)}]
\item the ramification profile over $0$ is $(\alpha_1,\ldots,\alpha_r)$;
\item the ramification profile over $\infty$ is $(\beta_1,\ldots,\beta_s)$;
\item away from $0$ and $\infty$, the ramification indices are $a_1+1,\ldots, a_n+1$.
\end{enumerate}
The ramification points in {\rm (iii)} are labeled by the zeros of $df/f$, while the branch points in the target are allowed to move and to coincide.
\end{prop}

\begin{proof}
Let $\omega\in \calF_\lambda(\mu)$, let $p_j$ denote the labeled pole with residue $\lambda_j$, and let $D_\lambda = \sum_{j=1}^p  \lambda_j p_j$. Note that $D_\lambda$ is a degree-zero divisor on $\mathbb{CP}^1$, hence it is principal. Therefore there exists a rational function $f$, unique up to multiplication by $\CC^*$, such that
\[
 {\mathrm{div}}(f) = D_\lambda.
\]
The differential $df/f$ has the same poles and the same residues as $\omega$.  Thus $\omega - df/f$ is holomorphic on $\mathbb{CP}^1$, and hence vanishes.
\par 
The positive and negative coefficients of $D_\lambda$ give the profiles of $f$ over $0$ and $\infty$.  If $q$ is neither a zero nor a pole of $f$ and locally
\[
  f(u) = c + u^e\cdot \mathrm{unit}, \quad c\in \CC^{*},
\]
then
\[ \frac{df}{f} = u^{e-1}\cdot \mathrm{unit}\,du.
\]
Thus a zero of $df/f$ of order $a$ is the same as a ramification point of $f$ of index $a+1$ away from $0$ and $\infty$.  This shows the desired correspondence.
\end{proof}

\begin{rmk}
The stratified polynomial spaces of \cite{SaBraid} are special cases of Definition~\ref{def:IH-fiber}.  If
\[
  f(z) = \prod_{i=1}^d (z - z_i)
\]
is a squarefree polynomial of degree $d$, i.e. the $z_i$ are all distinct, then $df/f$ has residue $1$ at each root and residue $-d$ at infinity. If the critical points of $f$ have multiplicity profile $\kappa$, this gives the isoresidual fiber with residue vector
\[
  (1,\ldots,1,-d)
\]
in the stratum $\calH(\kappa,-1^{d+1})$, modulo the isomorphisms of the source.  In this case the translation surface structure associated with $df/f$ is used in \cite{SaBraid} to study the braid monodromy of the polynomial stratum.  This is the same as the one appearing here, but for a particular residue vector.
\end{rmk}

We now combine Proposition~\ref{prop:df-f-dictionary} with Theorem~\ref{thm:MAIN}. Recall that an {\em open flat} of the resonance arrangement is determined by the collection of partial sums
\[
 \lambda_I = \sum_{i\in I}\lambda_i
\]
that vanish identically, with all other partial sums required to be nonzero.

\begin{cor}[Homeomorphic isoresidual Hurwitz fibers]
\label{cor:Hurwitz-homeomorphisms}
Let $B\subset \calR_p$ be an open flat of the resonance arrangement, and let $\lambda,\lambda'\in B$ be integral residue vectors with all coordinates nonzero. Then the corresponding two isoresidual Hurwitz fibers are homeomorphic:
\[
\cIH_\mu(\lambda) \cong  \cIH_\mu(\lambda').
\]
\end{cor}

\begin{proof}
By Theorem~\ref{thm:MAIN}, the residual map is locally topologically trivial over the open flat $B$. Since the complement of a complex hyperplane arrangement in a flat is connected, the fibers over any two points of $B$ are homeomorphic (possibly both empty).  If the two points are integral, Proposition~\ref{prop:df-f-dictionary} identifies these two fibers with the corresponding isoresidual Hurwitz fibers.
\end{proof}

\begin{ex}[An explicit Hurwitz homeomorphism]
If $\lambda' = m\lambda$ with $m\in \ZZ_{>0}$, they clearly belong to the same open flat, and one can see the above homeomorphism explicitly.  Indeed, the map
\[
  f\mapsto f^m
\]
sends $\cIH_\mu(\lambda)$ to $\cIH_\mu(m\lambda)$.  Conversely, if a rational function $g$ represents a point of $\cIH_\mu(m\lambda)$, then
\[
  {\mathrm{div}}(g) = mD_\lambda.
\]
Since $D_\lambda$ is principal on $\mathbb{CP}^1$, we can write $D_\lambda = {\mathrm {div}}(f)$, and then $g=cf^m$ for some $c\in\CC^*$.  Hence the inverse is obtained by dividing the differential $dg/g$ by $m$. 
\end{ex}

\begin{rmk}
In general, for two non-proportional integral residue vectors in the same open flat, there is usually no analogous formula as above on their  rational functions. If one fixes the marked poles $p_j$ and only changes the coefficients of the principal divisor $\sum \lambda_j p_j$, the zero type of $df/f$ will generally change.  The homeomorphism in Corollary~\ref{cor:Hurwitz-homeomorphisms} therefore comes from the isoresidual fibration: the poles and zeros move simultaneously so that the zero orders remain fixed.
\end{rmk}

\begin{rmk}
Note that Corollary~\ref{cor:Hurwitz-homeomorphisms} (unexpectedly) compares Hurwitz spaces that parameterize branched covers of different degrees. Nevertheless, it concerns the version of Hurwitz spaces in which the branch points away from $0$ and $\infty$ can move. In contrast, we do not claim that the classical Hurwitz spaces with fixed distinct branch points have the same topology. 
\end{rmk}

The same idea also translates known results on the image of the residual map into a criterion for the \emph{Hurwitz realization problem}. For an integral vector $\lambda \in \mathbb{Z}^{p}$, let $\lambda_{\mathrm{prim}}$ be the primitive integer vector that spans the same ray as $\lambda$.  Let 
\[
  d_{\mathrm{prim}}(\lambda) = \sum_{(\lambda_{\mathrm{prim}})_j>0}(\lambda_{\mathrm{prim}})_j, \quad A_{\max}(\mu) = \max\{a_1,\ldots,a_n\}.
\]

\begin{prop}[Isoresidual Hurwitz fiber with simple poles]
\label{prop:Hurwitz-realization}
Assume $p\ge 3$ and all poles are simple.  The isoresidual Hurwitz fiber $\cIH_\mu(\lambda)$ is nonempty if and only if
\[
  d_{\mathrm{prim}}(\lambda) > A_{\max}(\mu).
\]
Equivalently, the Hurwitz realization problem with profiles $\lambda^+$ over $0$, $\lambda^-$ over $\infty$, and ramification indices $a_i+1$ away from $0$ and $\infty$, is realizable precisely under the above inequality.
\end{prop}

\begin{proof}
This is the translation of the residue realization result  \cite[Th\'eor\`eme 1.2 (ii)]{getaab} in the case of genus zero and all poles being simple.  In this case the complement of the image of the residual map is the union of the $\CC^*$-lines spanned by primitive integral residue vectors whose positive degree is at most the largest zero order.  Proposition~\ref{prop:df-f-dictionary} translates that statement into the asserted Hurwitz realization criterion.

One can also prove the ``only if'' part directly as follows. Let $f\colon \mathbb{CP}^1\to \mathbb{CP}^1$ 
be a branched cover of degree $d$, and suppose that two branch points, say $0$ and $\infty$, have ramification profiles
\[
\alpha = (\alpha_1,\ldots,\alpha_r),\quad
\beta = (\beta_1,\ldots,\beta_s),
\]
all of whose parts are divisible by an integer $m\geq 1$.  Then
\[
(f) = \sum_i \alpha_i p_i-\sum_j\beta_j q_j
\]
is divisible by $m$ as a divisor.  Since every  divisor of degree zero on $\mathbb{CP}^1$ is principal, we can write $f=c h^m$ for some rational function
$h$ of degree $d/m$ with a nonzero scalar $c$.  Hence $f$ factors through the $m^{th}$-power map $z\mapsto z^m$.  As this map is unramified over $\CC^{*}$, the
ramification indices of $f$ away from $0$ and $\infty$ must already occur for $h$, and are therefore at most $d/m$. Thus any remaining branch profile containing a part larger than $d/m$ gives a
non-realizable Hurwitz profile.  
\end{proof}

\begin{ex}[An empty isoresidual Hurwitz fiber]
Let
\[
  \mu = (3,1,-1^6), \quad \lambda = (m, m, m, -m, -m, -m).
\]
Then $d_{\mathrm{prim}}(\lambda) = 3 = A_{\max}(\mu)$, and hence $\cIH_\mu(\lambda)$ is empty for all $m\ge 1$ by Proposition~\ref{prop:Hurwitz-realization}.  In terms of the Hurwitz realization problem, for $m\geq 2$, there is no rational map of degree $3m$ with fiber profiles
\[
  (m, m, m) \ \text{over }0, \quad (m, m, m) \ \text{over }\infty,
\]
and with additional ramification indices $4$ and $2$ away from $0$ and $\infty$, even if the corresponding branch points are allowed to move and collide.
\end{ex}

\begin{rmk}
Isoresidual Hurwitz fibers do not correspond to the {\em double ramification cycle} in genus zero. On $\mathbb{CP}^1$, every divisor of degree zero is principal, so the Abel--Jacobi condition defining the double ramification locus holds automatically. The extra condition here is that $df/f$ has zero type $(a_1,\ldots, a_n)$. Thus $\cIH_\mu(\lambda)$ is better understood as a Hurwitz stratum with fixed residues  and fixed profile of critical points. 
\end{rmk}

\section{An isoresidual elliptic surface}
\label{sec:isoresidual-elliptic-surface}

This section brings together three different kinds of geometry associated with the residual map. We consider the case where the projectivized stratum is a surface, and the residual map has fibers of dimension one, though these ideas generalize. First, from the viewpoint of flat geometry, each generic isoresidual curve carries a canonical translation structure defined by a meromorphic differential whose zeros and poles encode distinguished degenerations of the differentials parametrized by the curve. Second, after compactification, these isoresidual curves become the fibers of a rational elliptic surface: the zeros of the differential together give a multisection, while each pole gives a section of the elliptic fibration. Third, these sections, multisections, and fibers of the residual map are algebraic curves on the total surface and therefore contribute directly to its birational geometry.

We make this relationship explicit for the signature $\mu=(2, 2, -2, -2, -2)$. In this example the compactified generic isoresidual curve is an elliptic curve equipped with a meromorphic differential of signature $(4^4, -2^{8})$.  The four zeros form a degree-four multisection whose monodromy records the three resonance points, whereas the eight poles form distinguished sections of the elliptic fibration. The elliptic surface has Mordell--Weil rank two and consequently contains infinitely many sections. Their classes give infinitely many negative curves and hence infinitely many extremal rays of the pseudo-effective cone. The general fiber of the elliptic fibration forms an extremal moving curve that also detects the non-polyhedrality of the  pseudo-effective cone showing that the surface is \emph{not} a Mori dream space.
\par 
The example also serves to illustrate the distinction between the  compactifications considered in this paper: the minimal resolution of the projectivized residual map (see Remark~\ref{rmk:minimal-resolution}), realized by  $X_{\min}$ introduced below, and the finer residue-aligned multi-scale compactification introduced in Section~\ref{sec:RAMS}. Here ``minimal'' refers to resolving the rational residue map, rather than minimal as an elliptic fibration. 

For completeness, we recall the birational terminology used below. For a smooth projective surface $Y$, its {\em Picard number} is the rank of the Néron--Severi group $\rho(Y) = \operatorname{rank}\operatorname{NS}(Y)$. We write $\operatorname{Nef}(Y)$ for the {\em nef cone} of divisor classes having nonnegative intersection with every curve and $\overline{\operatorname{Eff}}(Y)$ for the closure of the cone generated by {\em effective divisors}. On a smooth projective surface, the intersection pairing identifies the nef cone of divisor classes with the closed cone of moving curves. After choosing a zero section of an elliptic fibration $\pi\colon Y\to \mathbb P^1$, its sections form the {\em Mordell--Weil group} $\operatorname{MW}(\pi)$. A positive Mordell--Weil rank produces infinitely many sections and, in the example below, infinitely many negative curves. This obstructs the rational polyhedrality required of a Mori dream space, or equivalently in this setting, the finite generation of its \emph{Cox ring}.

Consider the signature $\mu = (2, 2, -2, -2, -2)$. Let $X_{\min}$ be the normalization of the closure of the graph of the projectivized residual map in $\overline{\mathcal M}_{0,5}\times\PP\mathcal R_3$. Since $\PP\mathcal R_3\cong\PP^1$, projection to the second factor gives a morphism
\[
\PRes\colon X_{\min}\to \PP\mathcal{R}_3\cong \PP^1.
\]
We can construct this map explicitly. On the open stratum, normalize the three poles to $0,1,\infty$ and denote the two labeled zeros by $x$ and $y$. Up to scaling, the corresponding differential is given by
\[
\omega_{x, y} = \frac{(z-x)^2(z-y)^2}{z^2(z-1)^2}dz.
\]
The residues at $0,1,\infty$ are given by
\begin{align*}
r_0 & = 2xy(xy - x - y)  \\
r_1 & = 2(1-x)(1-y)(1-xy) \\
r_\infty & = -r_0 - r_1
\end{align*}
respectively. Hence we observe that the residue map is induced on $\PP^1\times \PP^1$ by the pencil of curves of bidegree $(2, 2)$ generated by 
\begin{align*}
f &= xy(xy - x - y)\\
g &= (1-x)(1-y)(1-xy).
\end{align*}
The base locus consists of the eight points
\begin{align*}
B_\partial &= \{(x,y)\mid x,y\in \{0,1,\infty\}, x\ne y\}\\
B_{\mathrm{res}} &= \{(x,y)\mid xy = 1, x+y = 1\}.
\end{align*}
The six points of $B_\partial$ are boundary points, while the two points of $B_{\mathrm{res}}$ are the residueless differentials.

Let $Z = \Bl_{B_\partial,B_{\mathrm{res}}}(\PP^1\times\PP^1)$ and $H_x, H_y$ denote the two ruling classes and $E_1,\ldots, E_8$ the exceptional curves. Every $E_i$ is a section of the resolved pencil and the strict transform of a member of the pencil has class
\[
F = 2H_x + 2H_y - \sum_{i=1}^8 E_i.
\]
At each of the eight base points, the two curves defined by $f$ and $g$ meet transversely, so each base point is simple. Consequently blowing up these eight points resolves the pencil, every $E_i$ is a section, and
\[
F = -K_Z, \quad F^2=0.
\]

\begin{thm}[The elliptic model]
\label{thm:theellipticmodel}
The resolved residue pencil
\[
\pi\colon Z\to \PP^1
\]
is a rational elliptic surface with three singular fibers, all of Kodaira type $IV$, lying over $r_0 = 0$, $r_1 =0$, and $r_\infty=0$. Further the Mordell--Weil group of this elliptic fibration has rank two. Moreover,
\[
X_{\min}\cong \Bl_{B_\partial,B_{\mathrm{res}}}(\overline{\mathcal{M}}_{0,5})\cong \Bl_{(0,0), (1,1), (\infty,\infty)}(Z).
\]
Hence $X_{\min}$ is obtained by blowing up the triple point of each type $IV$ fiber.
\end{thm}

\begin{proof}
After blowing up the eight base points the pencil is base-point-free, so a general member is smooth by Bertini. Adjunction, together with $F=-K_Z$ and $F^2=0$, shows that it has genus one. Since an exceptional curve over a base point is a section and $Z$ is rational, this gives a rational elliptic surface.

The fibers over $r_0 = 0$, $r_1 = 0$, and $r_\infty = 0$ are $f=0$, $g=0$, and $f+g=0$, respectively. In bihomogeneous coordinates we have
\begin{align*}
f &= X_1Y_1(X_1Y_1 - X_1Y_0 - X_0Y_1),\\
g &= (X_0 - X_1)(Y_0 - Y_1)(X_0Y_0 - X_1Y_1),\\
f+g &= X_0Y_0(X_0Y_0 - X_1Y_0 - X_0Y_1).
\end{align*}
Each of these fibers has two ruling components and one $(1,1)$-component. The ruling components pass through two base points, and the $(1,1)$-component passes through four, so all their strict transforms have self-intersection $-2$. The three components meet with distinct tangent directions at $(0,0)$, $(1,1)$, and $(\infty,\infty)$. 

Recall that a Kodaira fiber of type $IV$ is a reduced fiber consisting of three smooth rational $(-2)$-curves meeting at a single point with three distinct tangent directions. Its topological Euler characteristic is $4$. The factorizations above show that each of the three displayed fibers has precisely this configuration.

The topological Euler characteristic of $Z$ is
\[
\chi(Z) = \chi(\PP^1\times\PP^1) + 8 = 12.
\]
A type $IV$ fiber has Euler characteristic $4$, hence the three fibers exhaust the Euler characteristic and there are no other singular fibers. For a reducible fiber $v$, let $R_v$ denote the root lattice generated by the irreducible components that do not meet the chosen zero section. A type $IV$ fiber has three components. After removing the component met by the zero section, the remaining two $(-2)$-curves have intersection matrix
\[
\begin{pmatrix}
-2 & 1 \\
1 & -2
\end{pmatrix},
\]
the negative of the Cartan matrix of type~$A_2$. Thus $R_v\cong A_2$ and $\operatorname{rank}(R_v)=2$. Consequently, since $\rho(Z)=10$, the Shioda--Tate formula \cite{Shi90} gives
\[
10 = 2+\operatorname{rank}\operatorname{MW}(\pi)+3\cdot \operatorname{rank}(A_2) = 8+\operatorname{rank}\operatorname{MW}(\pi), 
\]
giving $\operatorname{rank}\operatorname{MW}(\pi)=2$.

The closure of the graph of $[f:g]$ is the blowup of its base ideal. Since its base scheme consists of the eight reduced points $B_\partial\cup B_{\mathrm{res}}$, this proves the first asserted description of $X_{\min}$. Finally,
\[
\overline{\mathcal M}_{0,5}\cong
\Bl_{(0,0),(1,1),(\infty,\infty)}(\PP^1\times\PP^1).
\]
The three diagonal points are disjoint from $B_\partial\cup B_{\mathrm{res}}$, so the blowups commute and give the second description.
\end{proof}

The ellipticity of the generic isoresidual curve in this signature also follows from \cite{CGPT}. The preceding minimal resolution is generally coarser than the residue-aligned multi-scale compactification.

Let $\Delta\subset \overline{\mathcal M}_{0,5}$ be the divisor where two zeros lie on a common tail. After normalizing the three top-level poles to $0,1,\infty$, write $c$ for the node and denote by $Q_0, Q_1, Q_\infty$ the points $c = 2,-1,\tfrac 12$, respectively.

\begin{prop}[Comparison with residue alignment]
\label{prop:RAMS-minimal-comparison}
There is a natural proper birational morphism
\[
\varpi\colon\RAMS(\mu)\to X_{\min}.
\]
It is not an isomorphism. At each of $Q_0,Q_1,Q_\infty$, one factor of the local residue-alignment ideal in \eqref{eq:Jres} has the form $(a,q^5)$ and its normalized blow-up is nontrivial.
\end{prop}

\begin{proof}
The forgetful morphism to $\overline{\mathcal M}_{0,5}$ together with the extended residual map on $\RAMS(\mu)$ gives a dominant morphism to the graph closure, and hence factors through its normalization $X_{\min}$. This produces the proper birational morphism $\varpi$.

Normalize the three poles on the top component to $0,1,\infty$, and let $c$ be the attaching point of the zero tail. The top-level differential is
\[
\eta_c = \frac{(z-c)^4}{z^2(z-1)^2}\,dz,
\]
with residues
\[
r_0 = 2c^3(c-2),\quad
r_1 = 2(1-c)^3(1+c),\quad
r_\infty = 4c-2.
\]
Thus the respective residues vanish simply at
\[
Q_0\colon c = 2,\quad Q_1\colon c = -1,\quad Q_\infty\colon c = \tfrac12.
\]
At the node separating the zero tail from the top component the orders are $-6$ and $4$, so the enhancement is $5$. If $q$ is the node-smoothing parameter, the lower-level adjusting parameter is therefore a unit times $q^5$; see \cite[Equation~(11.1)]{BCGGM2}. Taking the corresponding vanishing residue as a local coordinate $a$, the factor $(r_i,h_v)$ in \eqref{eq:Jres} is $(a,q^5)$. This ideal is not principal, so its blowup is nontrivial. Hence $\varpi$ is not an isomorphism: $X_{\min}$ remembers the projective residual direction, whereas $\RAMS(\mu)$ remembers additionally the scale of each individually vanishing residue relative to the lower level.
\end{proof}

The canonical translation structure on the generic isoresidual curves gives a geometric interpretation of the exceptional curves and of the three resonance fibers.

\begin{prop}[Translation singularities and their monodromy]
\label{prop:elliptic-translation-singularities}
Over
\[
B^\circ=\PP\mathcal R_3\setminus\{r_0 r_1 r_\infty=0\},
\quad X^\circ_{\min}=\PRes^{-1}(B^\circ), 
\]
the zero and polar divisors of the canonical translation differential satisfy
\[
Z(\omega) = 4\widetilde\Delta^\circ, \quad P(\omega) = 2\sum_{i=1}^8 E_i^\circ.  
\]
Thus the four zeros form the degree-four multisection $\widetilde\Delta^\circ$, whose permutation monodromy is $A_4$, whereas the eight labeled poles form the sections $E_1^\circ,\ldots,E_8^\circ$.
\end{prop}

\begin{proof}
The canonical translation structure is induced by the differential of the relative period between the two labeled zeros. By \cite[Propositions~4.1 and~4.7]{CGPT}, its zeros occur when these zeros merge, and for the present signature there are
\[
\frac{4!}{(4+2-3)!} = 4
\]
such points, each of order $4$. Their locus is the diagonal $x = y$. The same results give eight poles of order $2$: six boundary degenerations corresponding to the points of $B_\partial$ and two residueless degenerations corresponding to the points of $B_{\mathrm{res}}$. In the minimal resolution these are precisely the exceptional curves $E_1,\ldots,E_8$, where each $E_i$ maps isomorphically to the residual line and is therefore a section. This proves the asserted divisor formulas.

Restricting the residual map to the diagonal $x = y = c$ gives
\[
\phi\colon \PP^1\to \PP\mathcal{R}_3, \quad c\mapsto [c^3(c-2) : (1-c)^3(1+c)]
\]
Its fibers over the resonance points are
\[
\begin{aligned}
\phi^*(r_0=0) & = 3[0] + [2],\\
\phi^*(r_1=0) & = 3[1] + [-1],\\
\phi^*(r_\infty=0) &= 3[\infty] + [\tfrac12].
\end{aligned}
\]
These three points of ramification index $3$ account for the full ramification degree $2\deg(\phi) - 2 = 6$. Thus there is no further branching. The local monodromies are $3$-cycles, and they generate a transitive subgroup of $S_4$. Since they are even permutations, this subgroup is $A_4$.
\end{proof}

\begin{rmk}[Relation with the birational models]
The polar sections $E_i$ are eight distinguished $(-1)$-curves on $X_{\min}$. The strict transform $\widetilde\Delta$ of the zero multisection is also a $(-1)$-curve, since $\Delta^2=2$ on $\PP^1\times\PP^1$ and it passes through the three points $(0,0), (1,1), (\infty,\infty)$ that are blown up in passing from $Z$ to $X_{\min}$. These three points are the triply ramified points $c = 0,1,\infty$ of $\phi$. The remaining points over the resonance locus, namely $c = 2,-1,\tfrac12$, are $Q_0, Q_1, Q_\infty$ from Proposition~\ref{prop:RAMS-minimal-comparison}, where the finer residue-aligned compactification makes the additional local modifications $(a,q^5)$. Thus the zero multisection records both the blowups producing $X_{\min}$ and the further refinement from $X_{\min}$ to $\RAMS(\mu)$.
\end{rmk}

The elliptic model has the following birational consequence.

\begin{cor}
Let $F$ denote the class of a general isoresidual fiber on $X_{\min}$ for $\mu = (2, 2, -2, -2, -2)$. Then
\begin{enumerate}[label=\textup{(\roman*)}]
\item $F$ spans an extremal ray of $\operatorname{Nef}(X_{\min})$, equivalently of the cone of moving curves;
\item $\overline{\operatorname{Eff}}(X_{\min})$ is not polyhedral; and
\item $X_{\min}$ is not a Mori dream space, and its Cox ring is not finitely generated.
\end{enumerate}
Precisely, $\rho(X_{\min})=13$ and
\[
\dim_{\mathbb{R}}\langle [D]\mid D\geq 0,\ F\cdot D = 0 \rangle = 10 < \rho(X_{\min})-1.
\]
Thus $F$ is a non-polyhedral extremal moving ray in the sense of \cite[Definition~3.2 and Theorem~3.3]{mullane2025isoresidual}.
\end{cor}

\begin{proof}
For a smooth projective surface, the intersection pairing identifies the nef cone of divisors with the closed cone of moving curves. Hence the first statement follows from \cite[Theorem~3.4]{mullane2025isoresidual}, since $F$ is the class of a general fiber of a morphism.

The Mordell--Weil group has rank two, so $Z$ contains infinitely many sections. If $S$ is a section, then $K_Z\cdot S = -F\cdot S = -1$, and adjunction gives $S^2 = -1$. A section cannot pass through the triple point of a type $IV$ fiber: at such a point $d\pi = 0$, whereas $\pi\circ S = \operatorname{id}_{\PP^1}$ forces $d\pi\circ dS$ to be the identity. The strict transforms of the sections are therefore infinitely many $(-1)$-curves on $X_{\min}$. Each spans an extremal pseudo-effective ray, proving that $\overline{\operatorname{Eff}}(X_{\min})$ is not polyhedral. Since $X_{\min}$ is smooth, projective and rational, the Mori dream space and Cox-ring assertions follow from \cite[Proposition~2.9]{HK00}.

Finally, blowing up the three triple points increases $\rho(Z) = 10$ to $\rho(X_{\min}) = 13$. Each of the three special fibers on $X_{\min}$ has four irreducible components. An irreducible effective curve has zero intersection with $F$ precisely when it is vertical. The span of the vertical divisor classes consequently has dimension
\[
1 + 3\cdot (4-1) = 10,
\]
one dimension for the common fiber class and three independent component classes for each reducible fiber. This also recovers non-polyhedrality from \cite[Theorem~3.3]{mullane2025isoresidual}.
\end{proof}

We conclude this section by identifying the connected component of the stratum containing the canonical translation structures of the generic isoresidual fibers. This is a special case of Open Problem~9.8 of \cite{CGPT}.

\begin{prop}[Spin parity and rotation number]
\label{prop:elliptic-spin-parity}
For every $[\lambda]\in B^\circ$, the canonical translation structure
\[
(X_{[\lambda]},\omega_\lambda)\in\mathcal H(4^4, -2^8)
\]
has odd spin structure. Equivalently, it belongs to the connected component of rotation number $2$.
\end{prop}

\begin{proof}
Let $z_1,\ldots,z_4$ and $p_1,\ldots,p_8$ denote respectively the zeros and poles of $\omega_\lambda$. By Proposition~\ref{prop:elliptic-translation-singularities},
\[
\operatorname{div}(\omega_\lambda) =
4\sum_{j=1}^4 z_j - 2\sum_{i=1}^8 p_i.
\]
Consequently, the theta-characteristic line bundle
\[
\Theta_\lambda \coloneqq \mathcal O_{X_{[\lambda]}}
\left(
2\sum_{j=1}^4 z_j - \sum_{i=1}^8 p_i
\right)
\]
satisfies
$
\Theta_\lambda^{\otimes2} \cong
K_{X_{[\lambda]}}\cong
\mathcal O_{X_{[\lambda]}}$.
Thus $\Theta_\lambda$ defines a $2$-torsion point in 
$\operatorname{Pic}^0(X_{[\lambda]})~[2]$.

These points vary globally with $[\lambda]$. Indeed, the divisor
\[
2\widetilde\Delta^\circ-\sum_{i=1}^8 E_i^\circ
\]
on $X^\circ_{\min}$ defines a section
\[
\tau\colon
B^\circ\to
\operatorname{Pic}^0(X^\circ_{\min}/B^\circ)~[2],
\quad
\tau([\lambda]) = [\Theta_\lambda].
\]

For completeness, we recall the local monodromy calculation. The monodromy of a Kodaira fiber of type $IV$ can be computed from the standard local {\em Weierstrass model}
\[
y^2 = x^3 + t^2
\]
over a disk with coordinate $t$. After making the base change
$t = s^3$ and setting $x = s^2 X$, $y = s^3 Y$, the equation becomes
\[
Y^2 = X^3+1.
\]
Thus the pulled-back family is isotrivial. If $\zeta = \exp(2\pi {\mathrm{i}}/3)$, analytic continuation once around $t=0$
sends $s$ to $\zeta s$. This induces the order-three automorphism
\[
(X,Y)\mapsto (\zeta X,Y)
\]
of the elliptic curve $Y^2 = X^3 + 1$.

The induced action on the holomorphic differential $dX/Y$ is multiplication by $\zeta$. Hence the two eigenvalues of the action on first cohomology are $\zeta$ and $\zeta^2$, and its characteristic polynomial is
\[
u^2+u+1.
\]
After choosing a symplectic basis of $H_1(X_{[\lambda]},\mathbb Z)$, the monodromy is therefore up to inversion and conjugacy
\[
M= \begin{pmatrix}
0 & 1\\
-1 & -1
\end{pmatrix}.
\]
Notice that $M^3=I$. The blowups of the triple points used to pass from $Z$ to $X_{\min}$ affect only the central fibers, so they do not change the monodromy of the smooth family over the punctured base. Reducing $M$ modulo $2$ gives its action on the $2$-torsion which cyclically permutes the three nonzero elements,
\[
(1,0)\mapsto(0,1)\mapsto(1,1)\mapsto(1,0),
\]
and therefore fixes only the zero element of $(\mathbb Z / 2\mathbb Z)^2$. Since the section $\tau$ is globally defined over $B^\circ$, its value is invariant under this local monodromy. It follows that $\tau$ is the zero section and hence $\Theta_\lambda \cong \mathcal O_{X_{[\lambda]}}$ for every $[\lambda]\in B^\circ$.

Therefore, $h^0(X_{[\lambda]}, \Theta_\lambda) = h^0(X_{[\lambda]},\mathcal O_{X_{[\lambda]}}) = 1$, and the spin structure is odd. Finally, in genus one the algebraic spin parity of a component of rotation number $d$ is $d + 1$ (mod $2$), see \cite[Remark~5.5]{CG22}, hence the odd parity obtained above is equivalent to rotation number $d = 2$.
\end{proof}

\printbibliography
\end{document}